\documentclass[11pt]{amsart}

\usepackage{amssymb}
\usepackage{amsmath}
\usepackage{a4wide}
\usepackage{graphicx}
\newtheorem{theorem}{Theorem}[section]
\newtheorem{proposition}[theorem]{Proposition}
\newtheorem{lemma}[theorem]{Lemma}
\newtheorem{corollary}[theorem]{Corollary}
\newtheorem{remark}[theorem]{Remark}
\newtheorem{definition}[theorem]{Definition}

\makeatletter
 \@addtoreset{equation}{section}
 \makeatother

\newcommand{\fract}[2]{\genfrac{}{}{0pt}{}{\scriptstyle #1}{\scriptstyle #2}}
\newcommand{\trait}{\leavevmode\vrule height 2pt depth -1.6pt width 40pt}
\newcommand{\Subsection}[1]{\subsection{ #1} ${}^{}$}

\def\ds{\displaystyle}

\def\R{{\mathbb R}}
\def\C{{\mathbb C}}
\def\N{{\mathbb N}}

\def\CC{{\mathcal C}}
\def\E{{\mathbb E}}
\def\CD{{\mathcal D}}
\def\CE{{\mathcal E}}

\def\CI{{\mathcal I}}
\def\CJ{{\mathcal J}}
\def\CK{{\mathcal K}}
\def\CH{{\mathcal H}}
\def\CM{{\mathcal M}}
\def\CO{{\mathcal O}}
\def\CP{{\mathcal P}}
\def\CS{{\mathcal S}}
\def\CT{{\mathcal T}}

\def\CS{{\mathcal S}}

\def\diag{\mathop{\rm diag}\nolimits}
\def\im{\mathop{\rm  Im}\nolimits}
\def\re{\mathop{\rm Re}\nolimits}
\def\supp{\mathop{\rm supp}\nolimits}
\def\tr{\mathop{\rm tr}\nolimits}

\def\MS{\mathop{\rm MS}\nolimits}
\def\FS{\mathop{\rm FS}\nolimits}

\def\ker{\mathop{\rm Ker}\nolimits}
\def\ran{\mathop{\rm Ran}\nolimits}
\def\rank{\mathop{\rm Rank}\nolimits}
\def\bra{\langle}
\def\ket{\rangle}
\def\<{\langle}
\def\>{\rangle}

\def\op{\mathop{\rm Op}\nolimits}

\title[Microlocal kernel near an hyperbolic fixed point]{Microlocal kernel of pseudodifferential operators at an hyperbolic fixed point}

\begin{document}

\date{\today}

%\author[J.-F. Bony]{Jean-Fran\c{c}ois Bony}
%\address{Laboratoire MAB, CNRS, Universit\'e de Bordeaux I, Bordeaux, France}
%\address{Laboratoire MAB, CNRS, Universit\'e de Bordeaux I, Bordeaux, France}
%\email{Jean-Francois.Bony@math.u-bordeaux1.fr}
%\author[S. Fujiie]{Setsuro Fujiie}
%\address{Graduate School of Material Science, University of Hyogo, Japan}
%\address{Graduate School of Material Science, University of Hyogo, Japan}
%\email{fujiie@sci.u-hyogo.ac.jp}
%\author[T. Ramond]{Thierry Ramond}
%\address{Math\'ematiques, Universit\'e Paris XI, UMR CNRS 8628, Orsay, France}
%\address{Math\'ematiques, Universit\'e Paris XI, UMR CNRS 8628, Orsay, France}
%\email{thierry.ramond@math.u-psud.fr}
%\author[M. Zerzeri]{Maher Zerzeri}
%\address{D\'epartement de Math\'ematiques, Universit\'e Paris XIII, Villetaneuse, France}
%\address{D\'epartement de Math\'ematiques, Universit\'e Paris XIII, Villetaneuse, France}
%\email{zerzeri@math.univ-paris13.fr}

\author[J.-F. Bony]{Jean-Fran\c{c}ois Bony}
%\address{Laboratoire MAB, CNRS, Universit\'e de Bordeaux I, Bordeaux, France}
\address{Jean-Fran\c{c}ois Bony, Laboratoire MAB, UMR CNRS 5466, Universit\'e de Bordeaux I, Bordeaux, France}
\email{jean-francois.bony@math.u-bordeaux1.fr}
\author[S. Fujiie]{Setsuro Fujiie}
\address{Setsuro Fujiie, Graduate School of Material Science, University of Hyogo, Japan}
%\address{Graduate School of Material Science, University of Hyogo, Japan}
\email{fujiie@sci.u-hyogo.ac.jp}
\author[T. Ramond]{Thierry Ramond}
\address{Thierry Ramond, Math\'ematiques, Universit\'e Paris XI, UMR CNRS 8628, Orsay, France}
%\address{Math\'ematiques, Universit\'e Paris XI, UMR CNRS 8628, Orsay, France}
\email{thierry.ramond@math.u-psud.fr}
\author[M. Zerzeri]{Maher Zerzeri}
\address{Maher Zerzeri, D\'epartement de Math\'ematiques, UMR CNRS 7539,  Institut Galil\'ee, Universit\'e Paris XIII, Villetaneuse, France}
%\address{D\'epartement de Math\'ematiques, Universit\'e Paris XIII, Villetaneuse, France}
\email{zerzeri@math.univ-paris13.fr}

 
\begin{abstract}
We study the microlocal kernel of  $h$-pseudodifferential operators $\op_h(p)-z$, where $z$ belongs to some neighborhood of size $\CO(h)$ of a critical value of its  principal symbol $p_0(x,\xi)$. We suppose that this critical value corresponds to a hyperbolic fixed point of the Hamiltonian flow $H_{p_0}$. 
First we describe  propagation of  singularities at such a hyperbolic fixed point, both in the analytic and in the $\CC^\infty$ category. In both cases, we show that the null solution is  the only element of this microlocal kernel which vanishes on the stable incoming manifold, but for energies $z$ in some discrete set.
For energies $z$ out of this set,  we build the element  of the microlocal kernel  with prescribed data on the incoming manifold. We describe completely the operator which associate the value of this null solution  on the outgoing manifold to the initial data on the incoming one. In particular it appears  to be a semiclassical Fourier integral operator associated to some natural canonical relation.
\end{abstract}

%\begin{altabstract}
%On \'etudie le noyau microlocal d'op\'erateurs $h$-pseudodiff\'erentiels $\op_h(p)-z$, o\`u $z$ appartient \`a un voisinage de taille $\CO(h)$ d'une valeur critique  de son symbole principal $p_0(x,\xi)$. On suppose que celle-ci correspond \`a un point fixe hyperbolique du champ hamiltonien $H_{p_0}$ associ\'e.

%On d\'ecrit d'abord la propagation des singularit\'es au voisinage de ce  point fixe hyperbolique, en cat\'egorie analytique comme en cat\'egorie $\CC^\infty$. Dans ces deux cas, on montre que le noyau microlocal  ne contient pas de solution non-triviale qui s'annule sur la vari\'et\'e stable entrante, sauf pour un certain ensemble discret de valeurs de $z$.

%Pour  $z$ hors de cet ensemble exceptionnel, on construit l'unique \'el\'ement du noyau microlocal  de $\op_h(p)-z$ dont la valeur est donn\'ee sur la vari\'et\'e stable entrante. On d\'ecrit compl\'etement  l'op\'erateur qui associe \`a une telle donn\'ee la valeur de la solution correspondante sur la vari\'et\'e stable sortante. En particulier, nous montrons qu'il s'agit d'un op\'erateur int\'egral de Fourier semiclassique associ\'e \`a une relation canonique naturelle.

%\end{altabstract}

\thanks{\textbf{Acknowledgments:}  We would like to thank   B. Helffer  and J. Sj\"ostrand for many enlightening discussions during the preparation of this work.
}

\maketitle

\section{Introduction}

This paper is devoted to the study of the microlocal solutions near $(0,0)$ to the equation $(P-z)u=0$, where $P=\op_h(p(x,\xi,h))$ is  a self-adjoint $h$-pseudodifferential operator whose principal symbol can be reduced to
\begin{equation}\label{prems}
p_0(x,\xi)=\ds\sum_{j=1}^d\frac{ \lambda_j}{2} (\xi_j^2 - x_j^2) +  \CO((x, \xi )^3),
\end{equation}
for some  real and positive $\lambda_{j}$'s. The energies $z$ are supposed to lie at distance $\CO(h)$ of the critical value $p_0(0,0)=0$.
 
Of course such a situation occurs  for  a Schr\"odinger operator $-h^2\Delta+V$ when the potential $V$ has a non-degenerate local maximum, and the results of this paper might have many applications to quantum theory, allowing precise study of spectral or scattering quantities attached to these Schr\"odinger operators.

In this setting, the Hamiltonian vector field associated to $P$  has an hyperbolic fixed point at $(0,0)$, and the stable/unstable manifold theorem ensures  the existence of a stable incoming manifold $\Lambda_-$, and  of a stable outgoing manifold $\Lambda_+$ in $T^*\R^d$. The manifold $\Lambda_{-}$ (resp. $\Lambda_{+}$) can be described as the union of bicharacteristics $t\mapsto \gamma(t)$ such that $\gamma(t)\to (0,0)$ as $t\to +\infty$  (resp. as $t\to -\infty$). It is therefore a very natural question to ask, if the knowledge of a microlocal solution of the equation $Pu=0$ in $\Lambda_-$  determines the solution on $\Lambda_+$, thus in a whole neighborhood of the fixed point.

In the analytic, one-dimensional case, this problem has been given a complete  answer by B. Helffer and J. Sj\"ostrand in their study of Harper's operator \cite{hesjh3}. Their reduction to a normal form result (on the operator side), has then been used in several works, as for the study of gaps width  for Hill's equation by C. M\"arz \cite{marz} and the third author \cite{tr93}, or the computation of the scattering matrix at barrier tops \cite{tr96,fura}. There is also a series of work by Y. Colin de Verdi\`ere and B. Parisse \cite{cdvpa1,cdvpa2} about the so-called double-well problem where the same ideas are developed in a $\CC^\infty$ setting.

Here we address that question in the $d$-dimensional case, $d>1$.
We want to stress out the fact that  the results by N. Hanges, V. Ivrii or  R. Melrose, concerning propagation of singularities for  operators  with multiple characteristics (see e.g. \cite{ha}), do not apply here, since we are not in the case where the symbol factorizes as $p=p_{1}p_{2}$, with $p_1,p_2$ of principal type. Also, we don't think that a Birkhoff normal form reduction on the classical level can be used to obtain the results we give in this paper. In any case, such a reduction would require a non-resonant assumption on the $\lambda_j$'s, that we don't need here.

First, we  prove some kind of propagation of singularity result, both  in the analytic and in the $\CC^\infty$ category.
In these two categories, we show in Theorem \ref{un1} and Theorem  \ref{un2}  below that, roughly speaking, the null solution is the only microlocal solution of the equation $(P-z)u=0$ defined in a neighborhood of $(0,0)$, which vanishes on the stable incoming manifold $\Lambda_-$. This holds for energies $z$ in any neighborhood of the critical energy $0$ of size $\CO(h)$, that do not belong to some discrete subset $\Gamma(h)$. If $z\in \Gamma(h)$, then purely outgoing solutions exist - that is  solutions which vanish out of $\Lambda_+$.

In the analytic case,  our discussion is strongly related to the study of the resonances generated by a critical point of the principal symbol of a Schr\"odinger operator, and we use the same strategy as J.~Sj\"{o}strand  in \cite{sj87res} (see also \cite{kk} and \cite{sjreshyp}):  Our proof relies on  energy estimates  rather than on a reduction to a normal form.

In the $\CC^\infty$ case, our proof rely also on energy estimates, but these are obtained using quite different ideas from recent works by N. Burq and M. Zworski, S.H. Tang and M. Zworski (see  \cite{buzw} and \cite{tazw}), together with $h$-pseudodifferential calculus in some suitable class of symbols.

Then we turn to existence results in the  $\CC^\infty$ case:
For energies $z$ away from the discrete set $\Gamma(h)$,  we show the existence and give a representation formula for the solution  of $(P-z)u=0$  with given Cauchy data on $\Lambda_-$. 
Our proof relies heavily on  ideas from B. Helffer and J. Sj\"ostrand in \cite{hsmw3}, devoted to the study of the tunnel effect between non-resonant potential wells.
Thanks to this representation formula,  we build a microlocal transition operator, which associates the microlocal value of this solution on $\Lambda_+$ to the data on $\Lambda_-$.
We describe completely this operator (see Theorem \ref{explicit} and Theorem \ref{cho}), which turns out to be a $h$-Fourier Integral Operator associated to the  canonical relation $\Lambda_+ \times \Lambda_-$. 

%As a matter of fact, our initial motivation for this work was the description of the resonances associated  to a homoclinic orbit in a $d$-dimensional setting.
%Using these results we give here,  we are able to define and compute a quantum monodromy operator associated with such  a homoclinic orbit, in the spirit of the one given by J. Sj\"ostrand and M. Zworski in \cite{sjzw} for a closed orbit.
%As in that case, we think that we can extract a lot of information  from this monodromy operator, but we shall give results in that direction elsewhere.

The rest of the paper is organized as follows: In Section \ref{secamr}, we describe precisely our geometrical settings, give our assumptions, and state our results. Section \ref{unan} and Section \ref{uncinf}  are devoted to the proof of Theorem \ref{un1} and  of Theorem \ref{un2}, concerning the propagation of singularities at the hyperbolic fixed point, respectively in the analytic category, and in the $\CC^\infty$ category. Then, in Section \ref{secexis}, we address the question of existence of microlocal solutions of the natural Cauchy problem associated to our geometric setting, and we prove Theorem \ref{exis}. In Section 6 we obtain a precise formula for  that solution which is given in Theorem \ref{explicit} and \ref{cho}. Eventually, we have recalled in a short Appendix the results from $h$-pseudodifferential calculus that we use in Section 4.

%------------------------------
%
%-------------------------------
\section{Assumptions and main results}
\label{secamr}

\Subsection{Microlocal terminology}

Since our results are of microlocal nature, and since we shall constantly use this vocabulary through the paper, we briefly recall from \cite{sjsam} (see also \cite{mabk} and  \cite{de}) the precise meaning of  expressions like "$u=0$ microlocally in $\Omega$".
For $u\in \CS'(\R^d)$, we  denote  $\CT u$ the
Sj\"ostrand-FBI-Bargmann  transform of $u$ given by

\begin{equation}
\CT u(z,h)=c_d(h)\int e^{-(z-y)^2/2h} u(y) dy,
\label{bargmann}
\end{equation}
where $c_d(h)= 2^{-d/2}(\pi h)^{-3d/4}$ is a normalization constant. The function $\CT u$ is an holomorphic  function of $z\in \C^d$, and $\CT$  is isometric from $L^2(\R^d)$ to the Sj\"ostrand space $H_\Phi (\C^d)$, defined by
\begin{equation}
H_\Phi(\C^d)=L^2(e^{-2\Phi (z)/h} dz)\cap \CH(\C^d), \ \Phi(z)=\frac{(\im z)^2}2,
\label{sjostrandspace}
\end{equation}
where $ \CH(\C^d)$ is the space of holomorphic functions on $\C^d$, and $H_\Phi (\C^d)$ is endowed with the norm
\begin{equation}
\Vert f\Vert_{H_{\Phi}}=\left(\int \vert f(z,h)\vert^2 e^{-2\Phi(z)/h} dz\right )^{1/2}.
\label{normHphi}
\end{equation}
To the transform $\CT$, one also associates a canonical map $\kappa_\CT: T^*(\R^d)\to \C^d$ defined by
\begin{equation}
\kappa_\CT(x,\xi)=(x-i\xi, \xi).
\label{kappaT}
\end{equation}

We shall say that a family $(u_{h})_{h} \in \CS'(\R^d)$ is a tempered semiclassical distribution if there exists $N_{0}>0$ such that $h^{-N_{0}} u_{h}$ is bounded in $\CS'(\R^d)$. Such a  tempered semiclassical distribution $u\in \CS'(\R^d)$ is said to be analytically microlocally 0 in $\Omega$, an open subset of $T^*(\R^d)$, when there exists a constant $\varepsilon>0$ such that, 
\begin{equation}
\Vert \CT u \Vert_{H_{\Phi}(\Omega')}= \CO(e^{-\varepsilon/h}) \mbox{ as } h\to 0,
\label{egmicroloc}
\end{equation}
where $\Omega'=\Pi_1\kappa_\CT(\Omega)=\{x-i\xi, (x,\xi)\in \Omega\}$. The closed subset  of $T^*\R^d$ where  $u=(u_{h})_{h} $ is not analytically microlocally equal to 0 is called the microsupport of $u$, and we denote it by   $\MS(u)$.

In the $\CC^\infty$ category, one says that $u\in \CS'(\R^d)$ is  microlocally 0 in $\Omega$ when 
$\Vert \CT u\Vert_{H_{\Phi}(\Omega')}= \CO(h^\infty)$.
As a matter of fact, in this $\CC^\infty$ setting, we shall use $L^2$ norms instead of the above $H_{\Phi}$ norm, and it will be more convenient to use another version the FBI transform:  We set, for $z=x-i\xi$,
\begin{eqnarray}
\nonumber
\CT' u(x,\xi,h)&=&c_d(h)e^{-\xi^2/2h}\int e^{-(z-y)^2/2h} u(y) dy\qquad
\\
&=& c_d(h)\int e^{i(x-y)\xi/h -(x-y)^2/2h}u(y)dy.
\label{mfbi}
\end{eqnarray}
Then $\CT' u$ is a $\CC^\infty$ function on $\R^{2d}$, and $u\in \CS'(\R^d)$ is microlocally 0 in $\Omega$ if and only if $\Vert \CT' u\Vert_{L^2(\Omega)}= \CO(h^\infty)$.
The closed set of points where $u$ is not microlocally 0 is called  the frequency set of $u$, and we shall denote it by $\FS(u)$.

\Subsection{The geometrical setting}
\label{geoset}

We consider, microlocally near $(0,0) \in T^{*} \R^{d}$, a $h$-pseudodifferential operator
\begin{equation}  \label{hh1}
P = \op_{h} (p (x, \xi ,h) ),
\end{equation}
with symbol $p (x, \xi ,h) \in \CS_{h}^{0} (1)$ (see the Appendix \ref{pdo} for notations and a short review of $h$-pseudodifferential calculus). We assume that $p$ is real valued and
\begin{equation}
p (x, \xi ,h) \sim \sum_{j=0}^{\infty} p_{j} (x, \xi )h^{j},
\end{equation}
where the principal symbol satisfies, up to a symplectic change of variables,
\begin{equation}  \label{symbolpschro}
p_{0} (x,\xi)=\xi^2 - \frac14\sum_{j=1}^d \lambda_{j}^2x_{j}^2 + \CO((x, \xi)^3) , 
\end{equation}
in a neighborhood of $(0,0)$ in $T^* \R^d$. Here we have ordered the $\lambda_{j}$ such that
\begin{equation} \label{hh2}
0<\lambda_{1}\leq\lambda_{2}\leq \dots\leq \lambda_{d}.
\end{equation}
Since we work microlocally near $(0,0)$, we will assume that $p$ has compact support.

As usual, we denote by
\begin{equation}
H_{p} = \frac{\partial p_{0}}{\partial {\xi}} \frac{\partial}{\partial {x}} - \frac{\partial p_{0}}{\partial {x}} \frac{\partial}{\partial {\xi}}
\end{equation}
the Hamiltonian field of $p_{0} (x, \xi)$. In the $(x,\xi)$ coordinates, the linearized vector field $F_{p}$ of $H_{p}$ at $(0,0)$ is
\begin{equation}
F_{p}=d_{(0,0)}H_{p}=\left(
\begin{array}{cc}
0&2I\\
\frac12{L^2}&0
\end{array}
\right ),
\label{linearise}
\end{equation}
where $L$ is the $d\times d$ matrix defined as $L=\diag(\lambda_{1},\ldots,\lambda_{d})$. Then, the spectrum of $F_{p}$ is $\sigma (F_{p}) = \{ - \lambda_{d}, \ldots , - \lambda_{1}, \lambda_{1}, \ldots , \lambda_{d} \}$.
Associated to the hyperbolic fixed point, we have therefore a natural decomposition of $T_{(0,0)}(T^*\R^d)=\R^{2d}$ in a direct sum of two linear subspaces $\Lambda^0_{+}$ and $\Lambda^0_{-}$, of dimension $d$, associated respectively to the positive and negative eigenvalues of $F_{p}$.  These spaces $\Lambda^0_{\pm}$ are given by
\begin{equation}
\Lambda^0_{\pm}=\{(x,\xi)\in \R^{2d},\; \xi_{j}=\pm\frac{\lambda_{j}}2 x_{j}, \ j=1, \ldots,  d\}.
\label{equationlambda0}
\end{equation}

\begin{figure}
\begin{center}
\begin{picture}(0,0)%
\includegraphics{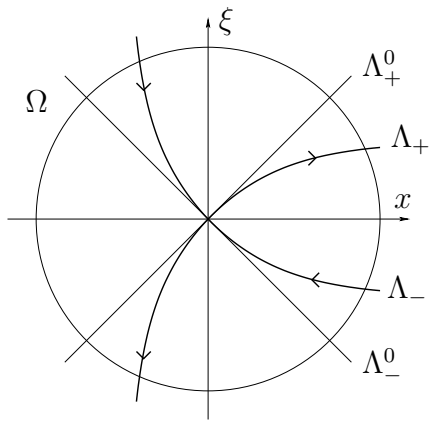}%
\end{picture}%
\setlength{\unitlength}{1184sp}%
\begingroup\makeatletter\ifx\SetFigFont\undefined%
\gdef\SetFigFont#1#2#3#4#5{%
  \reset@font\fontsize{#1}{#2pt}%
  \fontfamily{#3}\fontseries{#4}\fontshape{#5}%
  \selectfont}%
\fi\endgroup%
\begin{picture}(8424,8424)(1789,-8173)
\put(9226,-7186){\makebox(0,0)[lb]{\smash{{\SetFigFont{12}{14.4}{\rmdefault}{\mddefault}{\updefault}$\Lambda_{-}^ {0}$}}}}
\put(2176,-1711){\makebox(0,0)[lb]{\smash{{\SetFigFont{12}{14.4}{\rmdefault}{\mddefault}{\updefault}$\Omega$}}}}
\put(6226, 14){\makebox(0,0)[lb]{\smash{{\SetFigFont{12}{14.4}{\rmdefault}{\mddefault}{\updefault}$\xi$}}}}
\put(9901,-3736){\makebox(0,0)[lb]{\smash{{\SetFigFont{12}{14.4}{\rmdefault}{\mddefault}{\updefault}$x$}}}}
\put(9751,-5611){\makebox(0,0)[lb]{\smash{{\SetFigFont{12}{14.4}{\rmdefault}{\mddefault}{\updefault}$\Lambda_{-}$}}}}
\put(9826,-2461){\makebox(0,0)[lb]{\smash{{\SetFigFont{12}{14.4}{\rmdefault}{\mddefault}{\updefault}$\Lambda_{+}$}}}}
\put(9226,-961){\makebox(0,0)[lb]{\smash{{\SetFigFont{12}{14.4}{\rmdefault}{\mddefault}{\updefault}$\Lambda_{+}^ {0}$}}}}
\end{picture}%
\end{center}
\caption{The geometry at the singular point.}
\label{figsom}
\end{figure}

The stable/unstable manifold theorem gives us the existence of two
smooth Lagrangian manifolds $\Lambda_{+}$ and $\Lambda_{-}$, defined in a vicinity $\Omega$ of $(0,0)$, which are invariant under the $H_{p}$ flow, and whose tangent space at (0,0) are precisely $\Lambda^0_{+}$ and $\Lambda^0_{-}$.  In particular, we see that these manifolds can be written as
\begin{equation}
\Lambda_\pm=\{(x,\xi)\in T^*\R^d, \  \xi=\nabla\varphi_\pm(x)\},
\label{phi+-}
\end{equation}
for some smooth functions $\varphi_+$ and $\varphi_-$, which can be chosen so that 
\begin{equation}
\varphi_\pm(x)=\pm \sum_{j=1}^d\frac{\lambda_j}4x_j^2+ \CO(x^3).
\label{phi+-en0}
\end{equation}
Notice that if $P$ were a Schr\"{o}dinger operator, that is $p (x,\xi) = \xi^{2}
+ V(x)$, we would have $\varphi_{+} (x) = - \varphi_{-} (x)$.

We shall say that $\Lambda_{+}$ is the outgoing Lagrangian manifold, as $\Lambda_{-}$ will be referred to as the incoming
Lagrangian manifold associated to the hyperbolic fixed point. 
Indeed $\Lambda_{+}$ (resp. $\Lambda_{-}$) can be
characterized as the set of points $(x,\xi)\in \Omega$ such that $\exp tH_{p}(x,\xi)\to (0,0)$ as $t\to -\infty$ (resp.
as $t\to +\infty$).

\Subsection{Main results}
\label{mr}

Let $\Omega$ be a small neighborhood of $(0,0) \in T^{*} \R^{d}$. For $\varepsilon >0$ small enough, we  set $S = \Lambda_{-} \cap \{ (x, \xi ); \ \vert x \vert = \varepsilon\} \subset \Omega$. For $U \subset \Omega$ a neighborhood of $S$, we  
study the microlocal Cauchy problem
\begin{equation}
\label{zeeeee}
\left\{ \begin{aligned}
&(P-z) u = 0 &&\text{ microlocally in } \Omega ,   \\
&u = u_0 &&\text{ microlocally in } U.
\end{aligned} \right.
\end{equation}
Here $u_0\in L^2(\R^d)$ the microlocal Cauchy data, and we have to suppose that $(P-z)u_0=0$ microlocally in  $S$. We assume that $\Omega$ is small enough, so that  $P$ is of principal type in $\Omega\setminus\{(0,0)\}$. In particular, we have the usual propagation of singularity results away from the critical point.

First, we address the uniqueness problem for (\ref{zeeeee}). If $u_0=0$, the solutions have to vanish on the incoming manifold $\Lambda_{-}$, and we ask the question if the corresponding solution  is identically $0$ in a neighborhood of $(0,0)$.
The first two theorems below state that this is true both in the analytic
category and in the $\CC^\infty$ category, for complex energies $z\in D(0,C_0  h)=\{z\in \C, \vert z\vert<C_0h\}$,  where $C_0>0$ is any positive constant, but  for $z$ in some
discrete set. The existence of this exceptional set should not be too
surprising, at least in the analytic case: It corresponds to that of
resonances generated by the barrier top, i.e. the existence of "purely
outgoing solutions". In the $\CC^\infty$ case also, one could have
conjectured such a result. Indeed, the principal symbol $p_{0}$ can be
written in suitable coordinates $(y,\eta)$ as
\begin{equation}\label{formenormale}
p_{0} (x,\xi)=B(y,\eta)y\cdot \eta,
\end{equation}
 where $B$ is a smooth map from a neighborhood of $(0,0)$ in $T^*\R^d$ to the space $\CM_{d}(\R)$ of $d\times d$ matrices. 
Therefore in the one-dimensional case, $p_{0}$ factorizes as $p_{0} =
q_{1} q_{2}$, with $q_1$ and $q_2$ of principal type, and using a
reduction to a normal form as in the work \cite{ha} by N. Hanges,
 concerning propagation of singularities for  operators
with multiple characteristics, this uniqueness  result can be shown to
hold for $z$ away from the set
\begin{equation}
\{-ih \lambda_1(\alpha+\frac 12)\; ;\; \alpha\in \N\}.
\label{gamma}
\end{equation}
In the present multidimensional  setting, we find it convenient to work with the form (\ref{formenormale}) for our operator, and using  $h$-pseudodifferential calculus in some suitable class of symbols as well as ideas from Sj\"{o}strand in \cite{sj01} and Burq and Zworski in \cite{buzw}, we show the following result.

\begin{theorem}\sl  \label{un1}
Let $\Omega$ be a small neighborhood of $(0,0) \in T^{*} \R^{d}$, and   $S = \Lambda_{-} \cap \{ (x, \xi ); \ \vert x \vert = \varepsilon\} \subset \Omega$ for some $\varepsilon>0$ small enough. 
Assume \eqref{hh1}--\eqref{hh2}. Let $N$, $C_0  >0$ be constants, and  $U \subset \Omega$ a neighborhood of $S$. There exists a  neighborhood $V$ of $(0,0)$ such that, for all $z \in D(0, C_0  h)\subset \C$, and $u \in L^2(\R^d)$, defined for $h$ small enough  with $\Vert u \Vert _{L^2}\leq 1$,  if 
\begin{equation}
\left\{ \begin{aligned}
&(P-z) u = 0 &&\text{ microlocally in } \Omega ,   \\
&u = 0 &&\text{ microlocally in } U ,
\end{aligned} \right.
\end{equation}
with $d (z, \Gamma(h))> h^{N}$, then  $u=0$ microlocally in $V$.

Here,  $\Gamma(h)$ is a discrete set, defined for any $h$ small enough, such that $\#\Gamma(h)\cap D(0,C_0  h)$ is bounded uniformly with respect to $h$, and $\Gamma(h)\subset \{\im z < -\delta_{0}h\}$ for some $\delta_{0}>0$.
\end{theorem}

% \begin{figure}
% \input{figure4.pstex_t}
% \caption{The geometry at the singular point.}
% \end{figure}

In the analytic category, we can be as precise about the  exceptional
set  as in the one-dimensional case, changing of course the notion of
$C^\infty$-microsupport  to that of analytic microsupport.
Indeed, if we denote by $\Gamma_0(h)$ the discrete subset of $\C$ defined by 
 \begin{equation}
\Gamma_0(h)=\Big\{-ih\sum_{j=1}^d \lambda_j(\alpha_j+\frac 12), \
\alpha=(\alpha_1,\dots, \alpha_d)\in \N^d \Big\},
\label{gammaa}
\end{equation}
we have the following theorem which is, in some sense, a semiclassical version of a part of the work of Sj\"{o}strand \cite{sjhokk}.

\begin{theorem}\sl Suppose that, in addition to assumptions
\eqref{hh1}--\eqref{hh2}, the  function $p (x, \xi ,h )$ extends
holomorphically in a complex neighborhood of $(0,0)$ in $\C^{2d}$.
Let $\nu$, $C_0  >0$ be constants, and  $U \subset \Omega$ a
neighborhood of $S$.

There exists a  neighborhood $V$ of $(0,0)$ such that, for all $z \in
D(0, C_0  h)\subset \C$, and $u \in L^2(\R^d)$, defined for $h$ small
enough  with $\Vert u \Vert _{L^2}\leq 1$,  if
\begin{equation}
\left\{ \begin{aligned}
&(P-z) u = 0 &&\text{ analytically microlocally in } \Omega ,   \\
&u = 0 &&\text{  analytically  microlocally in } U ,
\end{aligned} \right.
\end{equation}
with $d(z(h), \Gamma_0(h))>\nu h$, then  $u=0$ analytically microlocally in $V$.
\label{un2}
\end{theorem}

Notice that, as in \cite{sj87res}, and using the ideas there, the last
assumption in Theorem \ref{un2} about the distance to the exceptional
set can certainly be replaced by a weaker one as in Theorem
\ref{un1}, provided the set $\Gamma_{0} (h)$  is replaced by
$\widetilde{\Gamma}_{0} (h) = \{ \lambda_{\alpha} (h) ; \ \alpha \in
\N^{d} \}$, where the $\lambda_{\alpha} (h)$ have an
expansion in fractional powers of $h$ and satisfy $\lambda_{\alpha} (h) = - i h
\sum_{1\leq j\leq d} \lambda_{j} (\alpha_{j} +1/2) + o (h)$.

\begin{remark}\sl In the $\CC^\infty$ category, and when the $\lambda_j$ are $\N$-independent, one can perform WKB construction of purely outgoing solutions for energies $z\in D(0,C_0h)$ such that $d(z,\Gamma_0(h))>\nu h$ (see e.g. \cite{hsmw1}). Therefore, in that particular case at least, we have $\Gamma_0(h) \subset \Gamma(h)$.
\end{remark}

\begin{remark}\sl The two previous theorems  can  be proved under  slightly more general  assumptions. Indeed for Theorem \ref{un1},  it  is sufficient to suppose that $P=\op_{h} (p)$, where 
\begin{itemize}
\item $p (x, \xi ;h) = p_{0} (x, \xi ) + h p_{1} (x, \xi) + h^{1+ \varepsilon} p_{2} (x, \xi ;h)$ for some $\varepsilon >0$.
\item $p_{0} (x, \xi )$ is a real valued $C^\infty$ function which can be written, up to a symplectic change of variables,
\begin{align*}
p_{0} (x,\xi)=\xi^2 - \frac14\sum_{j=1}^d \lambda_{j}^2x_{j}^2 + \CO((x, \xi)^3) , 
\end{align*}
\item $p_{1}$, $p_{2} \in \CS_{h}^{0} (1)$ (see Appendix \ref{pdo} for the definition of $\CS_{h}^{0} (1)$).
\end{itemize}
In that case,  the statement $\Gamma(h) \subset  \{\im z < -\delta_{0}h\}$ in Theorem \ref{un1} should be replaced by $\Gamma(h) \subset \{\im z < h (\im p_{1} (0,0) - \delta_{0}) \}$.

For the proof of Theorem \ref{un2}, we have to suppose in addition that $p$ extends as a holomorphic function  to a (fixed) neighborhood of $(0,0)$ in $\C^{2d}$.
\label{mga}
\end{remark}

Now, using ideas from B. Helffer and J. Sj\"ostrand in \cite{hsmw3},
we address the question of the existence  of solutions for the problem  (\ref{zeeeee}). As in that paper, to perform our construction  we have to suppose that the data is not microlocally supported on
some manifold of codimension 1 in $\Lambda_{-}$.
Indeed, we know from \cite{hsmw3} that there exist functions
$\gamma_{j}^{\pm} (t, x, \xi)$, polynomials with respect to $t$, such that,
in the precise sense of Definition \ref{expdeb} below,
\begin{equation}
\exp( \mp tH_{p}) (x,\xi) \sim \sum_{j\geq 1} \gamma_{j}^{\pm} (t, x,\xi)
e^{-\mu_{j}t}, \quad t\to +\infty,
\label{eqfuwsd}
\end{equation}
for all $(x,\xi)\in \Lambda_{\pm}$ respectively. Here $(\mu_{j})_{j\geq 0}$ is the
increasing sequence of linear combinations over $\N$ of the
$\lambda_{j}$'s. Moreover, the function $\gamma_{1}^{\pm}$ is a constant
vector with respect to $t$ in $\ker (d_{(0,0)}H_{p} \mp \lambda_{1})$.
We shall also consider $x$-space projections of the trajectories, and for $\rho\in \Lambda_\pm$ respectively, we shall denote
 \begin{equation}\label{g+g-}
g_1^\pm(\rho)=\Pi_x\gamma_1^{\pm}(\rho).
\end{equation}
Let $m$ be the number of $\lambda_{j}$'s equal to $\lambda_{1}$. We
denote by $\widetilde{\Lambda }_{\pm}$ the subset of $\Lambda_{\pm}$
which consists in points $(x,\xi)$ such that $\gamma_{1}^{\pm}
(x,\xi)=0$.
%Notice that,
%using Hartmann's Theorem \cite{hart}, one can see
%that $\widetilde{\Lambda_{\pm}}$ is a $\CC^1$ submanifold of $\Lambda_{\pm}$ of
%dimension $d-m$, which is stable under the Hamiltonian flow. But, in fact, it
%will be proved at the end of Section \ref{subsec:phase} that
Notice that, using the stable manifold master theorem \cite[Theorem 7.2.8]{abma}, one can see
that $\widetilde{\Lambda_{\pm}}$ is a $\CC^{\infty}$ submanifold of $\Lambda_{\pm}$ of
dimension $d-m$, which is stable under the Hamiltonian flow. As above, we
denote by $S \subset \Omega$ the lift in $\Lambda_{-}$
of the sphere $\{ x \in \R^{d} ; \ \vert x \vert = \varepsilon\}$, with
$\varepsilon>0$ small enough.

%
%
%
%\begin{remark}  \label{p1}
%The manifold $\widetilde{\Lambda_{\pm}}$ is a $\CC^{\infty}$ submanifold of $\Lambda_{\pm}$ of dimension
%$d-m$. 
%\end{remark}
%
%It will be proved in Section \ref{secexis} that the
%$\widetilde{\Lambda_{\pm}}$ are in fact $C^{\infty}$ submanifold. 

\begin{theorem} \sl Suppose assumptions \eqref{hh1}--\eqref{hh2} hold. Let ${{C_0}}$, ${{C_1}}$, $\nu >0$ be constants, $z \in [-{{C_0}}h , {{C_0}}h] + i [- {{C_1}} h , {{C_1}}h]$ with $d(z, \Gamma_0(h)) > \nu h$, and  $u_{0} \in L^{2} (\R^{d})$ be such that $\Vert u_{0} \Vert_{L^{2}} \leq 1$ with $u_{0} = 0$ microlocally in $S \cap
\widetilde{\Lambda_{-}}$ and $(P-z) u_{0} =0$ microlocally in $S$, then the problem
\begin{equation}  \label{dfj}
\left\{ \begin{aligned}
&(P-z) u =0 &&\text{ microlocally in } \Omega , \\
&u = u_{0} &&\text{ microlocally in } S ,
\end{aligned} \right.
\end{equation}
has a solution $u (x,z,h)$ such that 
\begin{equation}
\Vert u\Vert_{L^2}\lesssim h^{-\E\big(\frac{{{C_1}}}{\lambda_1}-\frac{\sum\lambda_j}{2\lambda_1}+\frac d2\big)-1},
\end{equation}
where $\E(r)$ is the integer part of $r\in \R$.

Moreover, if $u_{0}$ is analytic with respect to $z \in [-{{C_0}}h , {{C_0}}h] + i [- {{C_1}} h , {{C_1}}h] \setminus (\Gamma_{0}(h) + D (0 ,\nu h))$, then $u$ is also analytic.
\label{exis}
\end{theorem}

%\begin{remark}\sl
%Even in this $\CC^\infty$ setting, we need  that the distance between $z$ and the set $\Gamma(h)$ to be bounded from below by $\nu h$ in order to obtain our representation formula  for the solution (see \ref{}).
%\end{remark}

We denote by ${\mathcal J} (z) u_{0}$ the solution of the problem (\ref{dfj}), which is
unique thanks to Theorem \ref{un1}. Using a microlocal partition of
unity, we can assume that the initial data $u_{0}$ is microlocally
supported only in a vicinity of a point $\rho_{-} = (x^{-}, \xi^{-})
\in S \setminus \widetilde{\Lambda_{-}}$. As for Theorem \ref{exis},
we are unable to calculate the solution ${\mathcal J} (z) u_{0}$ near
every point of $\Lambda_{+}$ and we must avoid some particular set of
points $\widetilde{\Lambda_{+}}(\rho_{-})$: Let $\varphi_{1}$ be the solution of the Cauchy problem
\begin{equation}\label{transpphi1} 
\left\{ \begin{aligned}
& (  \nabla_\xi p_0(x,\nabla \varphi_{+})\cdot \nabla - \lambda_{1} ) \varphi_{1}
  = 0 , \\
&\nabla \varphi_{1} (0) = -\lambda_{1} g_1^-(\rho_-) .
\end{aligned} \right.
\end{equation}
We set  $\widetilde{\Lambda_{+}}(\rho_{-}) = \{ (x, \xi) \in
\Lambda_{+} ; \ \varphi_{1}(x) = 0 \}$.
Then, $\widetilde{\Lambda_{+}}(\rho_{-})$ is a $\CC^{\infty}$
submanifold of $\Lambda_{+}$, of codimension $1$, which is stable under
the Hamiltonian flow, and we can compute ${\mathcal J} (z) u_{0}$
near any point $\rho_{+} = (x^{+}, \xi^{+}) \in \Lambda_{+} \setminus
\widetilde{\Lambda_{+}}(\rho_{-})$. As the operator $P$ is of principal type in a neighborhood of $\rho_-$, and since $u_0$ is in the kernel of $P-z$, $u_0$ is completely determined by its trace on any hypersurface transversal to the flow. Up to a change of variables, we can assume that $x_1=x_1(\rho_-)=\varepsilon$ is such an hypersurface (taking the first coordinate function to be collinear to $g_1^-(\rho_-)$), and we state the following result in that setting. Eventually, because of  \eqref{symbolpschro}, and for $x' = o (x_{1})$, $\xi ' = o (x_{1})$, the equation $p_{0} (x, \xi_{1} , \xi ') = 0$ has two solutions
\begin{equation}
\xi_{1} = f_{\pm} (x, \xi ') = \pm \frac{\lambda_{1}}{2} x_{1} + o (x_{1}).
\end{equation}
In the Schr\"odinger case where $p(x,\xi)=\xi^2+V(x)$, we would have $ f_{\pm} (x, \xi ')=\pm\sqrt{-\xi'^2-V(x)}$. Then, with these notations, we have the following description for ${\mathcal J} (z) u_{0}$ near $\Lambda_+$.

\begin{theorem}\sl \label{explicit} We suppose that the assumptions of Theorem \ref{exis} hold, and that $u_0$ is microlocally supported only in a vicinity of  $\rho_-\in\CS\setminus\widetilde{\Lambda_{-}}$. We set
\begin{equation*}
S(z/h)=\sum_{k=1}^d \frac{\lambda_k}2-i\frac{z}{h},
\end{equation*}
and we denote by $(\widehat{\mu}_j)_{j\geq 0}$ the increasing sequence of the linear combinations over $\N$ of the $(\mu_{k}-\mu_1)$.
Then, there exists a symbol $d(x,y',z,h) \sim \sum_{j \geq 0} d_{j} (x,y',z, \ln h) h^{\widehat{\mu}_{j} / \lambda_{1}} \in \CS_{h}^{0} (1)$, with $d_{j} (x,y',z, \ln h)$ polynomial with respect to  $\ln h$, such that
\begin{equation}\label{zfinal7}
\CJ(z) u_0(x,h)=\frac{h^{S(z/h)/\lambda_{1}}}{(2\pi h)^{d/2}}\int_{\R^{d-1}} d(x,y',z,h) e^{i(\varphi_+(x)-\varphi_-(\varepsilon,y'))/h} u_0(\varepsilon,y') dy' ,
\end{equation}
microlocally near $\rho_+ \in\Lambda_{+} \setminus \widetilde{\Lambda_{+}}(\rho_{-})$. The symbols $d$ and $d_{j}$ are analytic for $z \in [-{{C_0}}h , {{C_0}}h] + i [- {{C_1}} h , {{C_1}}h] \setminus (\Gamma_{0}(h) + D (0 ,\nu h))$. Moreover the principal symbol $d_0$ of $d$ is independent of $\ln h$, and can be written as
\begin{equation}
\begin{aligned}
d_0(x,y',z) =& \sqrt{\lambda_{1}} e^{-id\pi/4} \Gamma \left( S(z/h) /\lambda_{1} \right) \big( i\lambda_{1} \big\< g^-_{1} (\rho_{( \varepsilon,y' )}^{-}) \vert g^+_{1}(\rho_{x}^{+}) \big\> \big)^{-S(z/h)/\lambda_{1}}   \\
&\big\vert g^-_{1} (\rho_{( \varepsilon,y' )}^{-}) \big\vert
\big\vert \det \nabla^2_{y'y'}\varphi_-(\varepsilon,y') \big\vert^{1/2}
\sqrt{\nabla_{\xi_1}p_0(\varepsilon,y'  ,f_-(\varepsilon, y',\eta'),\eta')} \\
&e^{\int_{0}^{- \infty} ( \Delta \varphi_{+} (x(t)) - \sum \lambda_{k}/2 ) dt}
\lim_{t\to + \infty} \frac{e^{(\sum \lambda_{k} /2 -\lambda_{1})t}}{{\sqrt{\det \frac{\partial y(t,w',y')}{\partial (t,w')}}}}  ,
\end{aligned}
\end{equation}
where $\rho^{\pm}_{x} = (x, \nabla_{x} \varphi_{\pm} (x))$,  and $x(t)$ (resp. $y(t,w',y')$) denotes the $x$--space coordinate of the hamiltonian curve $\exp (t H_{p}) (\rho^{+}_{x})$ (resp.  $\exp (t H_{p}) (\varepsilon,  , \nabla_{y} \varphi_{-} (\varepsilon, y'))$). 
\end{theorem}

\begin{remark}\sl Notice that, since every quantity in the previous theorem depends smoothly on $\varepsilon$, one can also consider the operator $\CJ$ as an operator on $L^2(\R^d)$. Indeed one has
\begin{equation}\label{extra}
\CJ(z) u_0(x,h)=\frac{h^{S(z/h)/\lambda_{1}}}{(2\pi h)^{d/2}}\int_{\R^{d}} \widetilde d(x,y,z,h)e^{i(\varphi_+(x)-\varphi_-(y))/h} u_0(y) d y.
\end{equation}
where $\widetilde d(x,y,z,h)=\chi(y_{1}) d(x,y_{1},y',z,h)$ for any function $\chi \in \CC^\infty_0(]0,\varepsilon_0[)$, with $\varepsilon_0>0$ small enough, such that $\int\chi(y_{1}) d y_{1} =1$. Here $d(x,y_{1},y',z,h)$ is the symbol given by Theorem \ref{explicit} with $\varepsilon = y_{1}$.
\end{remark}

In order to make even  clearer the fact that the microlocal transition operator $\CJ$ does not really depend on the choice of $\varepsilon$, we shall use the terminology of \cite{sjzw}: For $z\in [-{{C_0}}h , {{C_0}}h] + i [- {{C_1}} h , {{C_1}}h]$, we denote by $\CK_{\rho_{\pm}}(z)$ the
set  of distributions $u$  microlocally defined near $\rho_{\pm}$,  such
that  $(P-z)u =0$ microlocally near $\rho_{\pm}$. Notice  that,
since $P$ is of principal type away from $(0,0)$, there exist $U_{\pm}$,
$V_{\pm}$ two neighborhoods of $\rho_{\pm}$ and $(0,0)$
respectively and an elliptic microlocal $h$-Fourier integral operator
(an $h$-FIO from now on), ${\mathcal U}_{\pm} (z)$ with canonical transformation
$\kappa_{\pm}: U_{\pm} \to V_{\pm}$  such that
$\kappa_{\pm} (\rho_{\pm} ) =(0,0)$ and
\begin{equation} \label{eee}
{\mathcal U}_{\pm} (P-z) = h D_{x_{1}} {\mathcal U}_{\pm} \quad \text{microlocally in } U_{\pm}.
\end{equation}
Moreover, we have $\kappa_{\pm}^{*} \xi_{1} := \xi_{1} \circ \kappa_{\pm} (x, \xi
) = p (x,\xi)$. (see e.g. \cite[Proposition 3.5]{sjzw} in this semiclassical setting). Then $\CK_{\rho_{\pm}}(z)$ can
be identified with $\CD'(\R^{d-1})$ using ${\mathcal U}_{\pm}$.

Let $v_{-} \in \CD'(\R^{d-1})$ be microlocally supported in a compact
subset of $V_{-}$. If $u_{-}$ is the corresponding element in
$\CK_{\rho_{-}}(z)$ and $u$ the solution of (\ref{dfj}) with
initial data $u_{-}$, we denote by ${\mathcal I}(z) v_{-}$ the
element of $\CD'(\R^{d-1})$ corresponding to $u$ near $\rho_+$. In other words, we have set
\begin{equation}
{\mathcal I} (z) = \imath^{*} \; {\mathcal U}_{+} \; {\mathcal J} (z)
\; {\mathcal U}_{-}^{-1} \; \pi^{*},
\label{zeI}
\end{equation}
where $\imath : x' \mapsto (0, x' )$ and $\pi : (x_{1} , x') \mapsto
 x'$.

\begin{theorem} \label{cho}  \sl
Assume $z\in [-{{C_0}}h , {{C_0}}h] + i [- {{C_1}} h , {{C_1}}h]$ and $d (z, \Gamma_{0} (h)) >  \nu h$. Then the operator
$\CI (z)$  is a $h$-Fourier integral operator of order $h^{\re S(z/h) - \frac{1}{2}}$ on $L^{2} (\R^{d-1})$,
microlocally defined near $(0,0)$, analytic with respect to $z$, associated to the canonical relation
\begin{equation}
\CC_{\CI} = \Pi \circ \kappa_{+} ( \Lambda_{+} )  \times \Pi \circ
\kappa_{-} ( \Lambda_{-} ),
\end{equation}
where $\Pi : (x_{1}, x' ,  \xi_{1} , \xi ') \mapsto (x', \xi ')$.
\end{theorem}

\begin{remark} \sl
The canonical relation does not depend on the choice of $\kappa_\pm$ in the following
sense. Suppose that $\widetilde{U}_{\pm}$ are others FIO's, with
canonical relation $\widetilde{\kappa}_{\pm}$, as in the discussion
before Theorem \ref{cho}. The operators $\widehat{\mathcal U}_{\pm}
= \widetilde{{\mathcal U}}_{\pm} {\mathcal U}^{-1}_{\pm}$ are FIO's with
canonical relation $\widehat{\kappa}_{\pm} = \widetilde{\kappa}_{\pm} \circ
\kappa_{\pm}^{-1}$. Then, we see that $\widehat{\kappa}_{\pm}$ must be
of the form
\begin{equation}
\widehat{\kappa}_{\pm} (x, \xi ) = \big( f_{1}^{\pm} (x, \xi ) , g_{x'}^{\pm} (x'  , \xi
') + \xi_{1} f_{2}^{\pm} (x, \xi ) , \xi_{1} , g_{\xi '}^{\pm} (x', \xi ') +
\xi_{1} f_{3}^{\pm} (x, \xi ) \big) ,
\end{equation}
where $(x, \xi ) = (x^{1},x',\xi^{1},\xi')$. Then, Lemma 3.4 of \cite{sjzw} implies that $\imath^{*} \;\widehat{\mathcal U}_{\pm}\;  \pi^{*}$ is an FIO on
$L^2(\R^{d-1})$ with canonical transformation
\begin{equation}
g_{\pm} : (x ', \xi ') \mapsto ( g_{x'}^{\pm} (x', \xi ') , g_{\xi '}^{\pm} (x', \xi ') ).
\end{equation}
Therefore, if we denote by $\widetilde{\CI}(z)$ the same operator as $\CI(z)$ but defined through $\widetilde{\mathcal U}_{\pm}$ instead of ${\mathcal U}_{\pm}$, we have
\begin{equation}
\CC_{\widetilde{\CI}}(z) = (g_{+}\times g_{-}) (\CC_{\CI}(z)).
\end{equation}
\end{remark}

%----------------------
%
%----------------------
\section{Uniqueness in the analytic case}
\label{unan}
 
We prove Theorem  \ref{un2}. Since this uniqueness statement  is essentially equivalent to the fact that there is no purely outgoing solution, it should not be surprising that our discussion is strongly related to the study of the resonances generated by a maximum of $V(x)$, and we use the same strategy as J. Sj\"{o}strand  in \cite{sj87res} (see also \cite{kk}), as well of some lemmas from that paper or from \cite{gesj}.

In this section, as for example in Figure 2,  we use the same notations for subsets of $T^*\R^d$ and their image in $\C^d$  by  $(x,\xi)\mapsto x-i\xi$. We recall that, using also this convention, we shall say that $u$ is microlocally 0 in $\Omega$ if  $\MS(u)\cap \Omega=\emptyset$,

We work under the assumptions \eqref{hh1}--\eqref{hh2}: We set  $P=\op_{h}(p)$, where $p$ is a holomorphic function, depending on $h\in ]0,1]$ say,  in a (fixed) complex neighborhood of $(0,0)$ in $\C^{2d}$.
We also assume that, up to a linear change of variables, $p_{0}$ can be written as
\begin{equation}
p_{0}(x,\xi ) = \sum_{j=1}^d\frac{ \lambda_j}{2} (\xi_j^2 - x_j^2) + \CO((x, \xi )^3)
\end{equation}
for some real and positive $\lambda_{j}$'s. We start with this expression  for $p_{0}$.

\begin{figure}% [!ht]
\begin{center}
\begin{picture}(0,0)%
\includegraphics{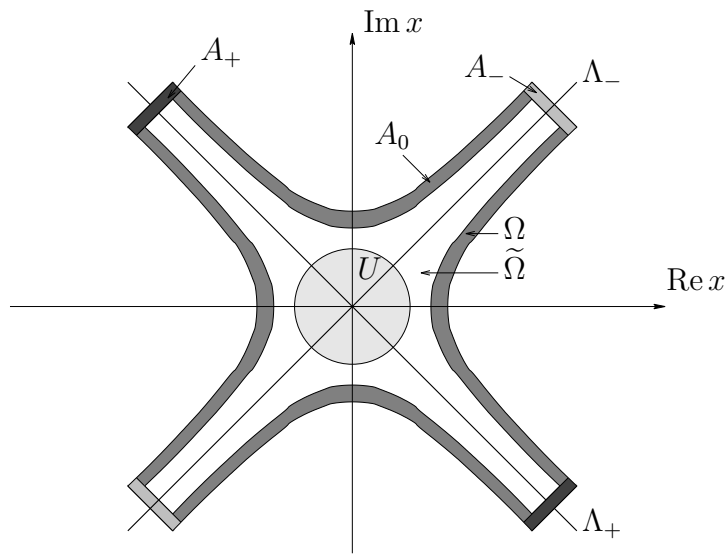}%
\end{picture}%
\setlength{\unitlength}{1855sp}%
\begingroup\makeatletter\ifx\SetFigFont\undefined%
\gdef\SetFigFont#1#2#3#4#5{%
  \reset@font\fontsize{#1}{#2pt}%
  \fontfamily{#3}\fontseries{#4}\fontshape{#5}%
  \selectfont}%
\fi\endgroup%
\begin{picture}(8799,7119)(1414,-8473)
\put(8026,-4786){\makebox(0,0)[lb]{\smash{{\SetFigFont{12}{14.4}{\rmdefault}{\mddefault}{\updefault}$\widetilde{\Omega}$}}}}
\put(6076,-4786){\makebox(0,0)[lb]{\smash{{\SetFigFont{12}{14.4}{\rmdefault}{\mddefault}{\updefault}$U$}}}}
\put(7501,-2086){\makebox(0,0)[lb]{\smash{{\SetFigFont{12}{14.4}{\rmdefault}{\mddefault}{\updefault}$A_{-}$}}}}
\put(6301,-2986){\makebox(0,0)[lb]{\smash{{\SetFigFont{12}{14.4}{\rmdefault}{\mddefault}{\updefault}$A_{0}$}}}}
\put(3976,-1861){\makebox(0,0)[lb]{\smash{{\SetFigFont{12}{14.4}{\rmdefault}{\mddefault}{\updefault}$A_{+}$}}}}
\put(9076,-8161){\makebox(0,0)[lb]{\smash{{\SetFigFont{12}{14.4}{\rmdefault}{\mddefault}{\updefault}$\Lambda_{+}$}}}}
\put(9076,-2161){\makebox(0,0)[lb]{\smash{{\SetFigFont{12}{14.4}{\rmdefault}{\mddefault}{\updefault}$\Lambda_{-}$}}}}
\put(8026,-4261){\makebox(0,0)[lb]{\smash{{\SetFigFont{12}{14.4}{\rmdefault}{\mddefault}{\updefault}$\Omega$}}}}
\put(6151,-1486){\makebox(0,0)[lb]{\smash{{\SetFigFont{12}{14.4}{\rmdefault}{\mddefault}{\updefault}$\im x$}}}}
\put(10201,-4936){\makebox(0,0)[lb]{\smash{{\SetFigFont{12}{14.4}{\rmdefault}{\mddefault}{\updefault}$\re x$}}}}
\end{picture}%
\end{center}
\caption{The domains.}
\label{figdomain}
\end{figure}

As in the discussion of Section \ref{geoset}, we work in some neighborhood $\Omega$ of the fixed point $(0,0)$, and we choose $ \Omega_{1} \Subset \Omega_{0}=\Omega$ as in Figure \ref{figdomain}. We write 
\begin{equation}
\Omega_{0} \setminus \Omega_{1}=A_+\cup A_0\cup A_-,
\label{omegas}
\end{equation}
where $A_\pm$ is close to $\Lambda_\pm$. We assume that $A_0$ is geometrically controlled by $A_-$, that is any point $(x,\xi)\in A_0$ can be written as $\exp tH_p (x_-,\xi_-)$ for some $(x_-,\xi_-)\in A_-$ and some $t>0$.
It is clear that one can find such a configuration when  $H_{p} = F_{p}$, and  Hartmann's  Theorem (see e.g. \cite{perko}) ensures  that we can do so in the general case as well.

We consider  the operator on $H_{\Phi}(\Omega)$ defined by
\begin{equation}
\widetilde{P} =\CT P\CT^*,
\label{Q1}
\end{equation}
where $\CT$ is the FBI transform given in (\ref{bargmann}), and $H_{\Phi}(\Omega)$ is defined in (\ref{sjostrandspace}). Then $\widetilde P$  is a pseudodifferential operator in the complex domain (see J. Sj\"{o}strand \cite{sjsam}). Its principal symbol is
\begin{equation}
\widetilde{p}_{0} (x,\xi)=p_{0}\circ\kappa_{\CT}^{-1}( x , \xi ) = \sum_{j=1}^d \frac{\lambda_j}{2} ( 2 \xi_j^2 -2i x_j \xi_{j} - x_{j}^2)+ \CO ((x,\xi)^3).
\label{q}
\end{equation}

First,  $u=0$ microlocally in $\Lambda_-\setminus \{(0,0)\}$, so that  we can assume that $u=0$ microlocally in $A_-$
provided $\Omega_0$ is small enough. Since $A$ is geometrically controlled by $A_-$,  we get from standard results on propagation of singularities, that, for some $\delta>0$,
\begin{equation}
\Vert\CT u\Vert_{H_{\Phi}(A_-\cup A_0)}=\CO (e^{-\delta/h}).
\label{mlk}
\end{equation}

Now, we choose $U$ an elliptic FIO with complex phase given by
\begin{equation}
\varphi_{U} (x,y)  = \frac{x^{2}}{2} - x y+ \frac{(1-i)}{4} y^{2} .
\end{equation}
This operator is associated to the complex canonical transform
\begin{equation}
\kappa_U : (x,\xi) \mapsto \left( \frac{2 \xi + (1-i) x}{2} , \frac{2 \xi - (1+i) x}{2} \right).
\label{kappaU}
\end{equation}
Notice that the operator $U$ cannot be realized on $H_{\Phi}$ since the function $y \mapsto -\im (\varphi_{U} (x,y)) + \Phi (y)$ has no saddle point.  However, if we set
\begin{equation}  \label{az12}
G(x) = - \re x \im x ,
\end{equation}
then, for $t>0$  fixed, $U$ is well-defined as an operator  from $H_{\Phi + t G}(\Omega_2)$ to $H_{\Psi_{t}}(\kappa_U(\Omega_3))$, where  $\Psi_{t}$ is some plurisubharmonic function. 
Here  $\Omega_3\subset \Omega_2$ are suitable neighborhoods of $(0,0)$ depending on $t$, since the saddle  point of $y \mapsto -\im (\varphi_{U} (x,y)) + \Phi (y) +tG(y)$ does. From \cite{sjsam}, we can invert $U$ by an FIO 
$V$ from  $H_{\Psi_{t}} $ to   $H_{\Phi + t G}$ up to exponentially small errors, taking care of domains. Now we set, after shrinking $\Omega_2$ and $\Omega_3$,
\begin{equation}
 Q=U \widetilde{P} V : H_{\Psi_{t}} (\kappa_U(\Omega_{2})) \to H_{\Psi_{t}} (\kappa_U(\Omega_{3})),
\end{equation}
which is a pseudodifferential operator with principal symbol
\begin{equation}
q_{0} (x,\xi) = \widetilde{p}_{0}  \circ \kappa_{U}^{-1} (x,\xi)  = \sum_{j=1}^d \lambda_j x_{j}  \xi_{j} +  \CO ((x,\xi)^3).
\end{equation}

Let us  recall Proposition 4.4 from \cite{gesj}:

\begin{proposition}[\cite{gesj}, Proposition 4.4]\sl \label{hhh}
Let $\chi\in \CC_0^\infty(\kappa_U(\Omega_{3}))$. There exists a classical symbol $\widetilde{q} (x, h)$ of order 0 such that
\begin{equation}
\bra \chi Q u, v\ket_{H_{\Psi_{t}}(\kappa_U(\Omega_{3}))} = \bra \widetilde{q} u,v \ket_{H_{\Psi_{t}} (\kappa_U(\Omega_{3}))} + r(u,v),
\label{puv}
\end{equation}
where
\begin{equation}
r(u,v)= \CO (h^\infty)\Vert u\Vert_{H_{\Psi_{t}} (\kappa_U(\Omega_{2}))}   \Vert v\Vert_{H_{\Psi_{t}} (\kappa_U(\Omega_{2}))}   ,
\label{ruv}
\end{equation}
and
\begin{equation}
\widetilde{q}_0 (x) = \chi(x) q_{0} \Big( x , \frac 2i \partial_x \Psi_{t} (x) \Big).
\label{q0}
\end{equation}
\end{proposition}

We  use  Proposition \ref{hhh}  with  $\chi =1$ near $\kappa_U(\Omega_{4})$, for some $\Omega_4\Subset\Omega_3$.  First, for $x\in \kappa_U(\Omega_3)$, we have
\begin{align}
q_{0} \Big( x , \frac 2i \partial_x \Psi_{t} (x) \Big) &= \widetilde{p}_{0} \Big(  y, \frac{2}{i} \partial_{y} (\Phi + tG) \Big)_{|_ {y = \kappa_{U}^{-1} (x)}}   \nonumber \\
&= p_{0} \Big( a +2 t \partial_{z} G (a-ib) , b-2i t \partial_{z} G (a-ib)  \Big)_{|_ {z = a  -ib = \kappa_{U}^{-1} (x)}},
\end{align}
with $\partial_{z} = (\partial_{a} + i \partial_{b}) /2$. In particular, we have
\begin{align}
- \im q_{0} (x,\frac 2i \partial_x \Psi_{t} (x) )\geq & \sum_{j=1}^{d} \lambda_{j} t (a_{j}^{2}  + b_{j}^{2} ) + \CO ( t^{2}  (a,b)^{2} + t (a,b)^{3} )_{|_ {a  -ib = \kappa_{U}^{-1} (x)}}    \nonumber  \\
\geq & \frac{t}{C_{0}} \vert x \vert^{2},
\end{align}
for $0<t$ and $x$ small enough. Therefore, since $z(h) \in D(0, C_{0}h)$, we obtain
\begin{align}
- \im \bra \chi & (Q - z(h)) u, u \ket_{H_{\Psi_{t}} (\kappa_U(\Omega_{3}))}  \nonumber  \\
&\geq \bra \frac{t}{C} \vert x \vert^{2}  u,u \ket_{H_{\Psi_{t}} (\kappa_U(\Omega_{4} ))} + \CO (h) \Vert u \Vert^{2}_{H_{\Psi_{t}} (\kappa_U(\Omega_{2}))}  \nonumber \\
&\geq h \Vert u \Vert^{2}_{H_{\Psi_{t}} (\kappa_U(\Omega_{4} ))} + \CO (h) \Vert u \Vert^{2}_{H_{\Psi_{t}} (\kappa_U(\Omega_{2}\setminus\Omega_{4}))}  + \CO (h) \Vert u \Vert^{2}_{H_{\Psi_{t}} ( \{ \vert x \vert < C_{1} \sqrt{h}\})}.    \label{qsd}
\end{align}

For  $n \in \N$, we denote by $\tau_n: H_{\Psi_{t}} (\kappa_U(\Omega_{2})) \to H_{\Psi_{t}} (\kappa_U(\Omega_{2}))$ the operator defined as
\begin{equation}
\tau_n(v)=\sum_{\vert\alpha\vert<n}\frac{1}{\alpha !}\partial_x^{\alpha}v(0) x^\alpha,
\label{taun}
\end{equation}
and we recall  the following

\begin{lemma}[\cite{gesj}, Lemme 4.5]\sl
Let $0<C_1<C_2$ be fixed constants. There exists  a sequence $(c_n)_{n}$ of real positive numbers such that $c_n\to 0$ as $n\to +\infty$, and, for any $n\in \N$, for any $v\in H_{\Psi_{t}} (\kappa_U(\Omega_{2})) \cap \ker \tau_{n}$ it holds that,
\begin{equation}
\Vert v\Vert_{H_{\Psi_{t}} ( \{ \vert x \vert < C_{1} \sqrt{h}\})}  \leq {c_n}\Vert v\Vert_{H_{\Psi_{t}} ( \{ \vert x \vert < C_{2} \sqrt{h}\})}.
\end{equation}
\label{aprioriholn}
\end{lemma}
Writing  (\ref{qsd})  for $u \in \ker \tau_{n}$ and $n$ large enough, we obtain
\begin{equation}
- \im \bra \chi (Q - z) u, u \ket_{H_{\Psi_{t}} (\kappa_U(\Omega_{3}))} \geq h \Vert u \Vert^{2}_{H_{\Psi_{t}} (\kappa_U(\Omega_{4}))} + \CO (h) \Vert u \Vert^{2}_{H_{\Psi_{t}} (\kappa_U(\Omega_{2}\setminus \Omega_{4}))} ,
\end{equation}
so that
\begin{equation}   \label{vbn}
h \Vert u \Vert_{H_{\Psi_{t}} (\kappa_U(\Omega_{4}))} \leq \Vert (Q-z) u \Vert_{H_{\Psi_{t}} (\kappa_U(\Omega_{3}))} + \CO (h) \Vert u \Vert_{H_{\Psi_{t}} (\kappa_U(\Omega_{2}\setminus \Omega_{4}))} .
\end{equation}

Now, we come back to the initial problem and we suppose that the assumptions of Theorem \ref{un2} hold. Since $(P-z(h)) u = 0$ analytically microlocally in $\Omega_0$, we have, for $0<t$ small enough,
\begin{equation}  \label{az7}
(Q-z) U \CT u = \CO ( e^{- \varepsilon / h}),
\end{equation}
in $H_{\Psi_{t}} (\kappa_U(\Omega_{2}))$ for some $\varepsilon >0$. Notice that the constant $\varepsilon $ may change from line to line in what follows,  and depends on  $t$. Applying $1- \tau_{n}$, we obtain
\begin{align}
(Q-z) (1-\tau_{n}) U \CT u &= [\tau_{n},Q] U \CT u + \CO ( e^{- \varepsilon  / h})  \nonumber  \\
&= \CO (h^{3/2}) \Vert U \CT u \Vert_{H_{\Psi_{t}} (\kappa_U(\Omega_{2}))} + \CO ( e^{- \varepsilon / h}) ,
\end{align}
in $H_{\Psi_{t}} (\kappa_U(\Omega_{3}))$. Here, we have used the fact that $\tau_{n} = \CO (1)$ and $[\tau_{n},Q] = \CO (h^{3/2})$ thanks to Proposition 4.3 in \cite{gesj} and Proposition 3.3 in \cite{sj87res}. Then, we get from (\ref{vbn}) the estimate
\begin{equation}  \label{az1}
\begin{aligned}
\Vert (1-\tau_{n}) U \CT u \Vert_{H_{\Psi_{t}} (\kappa_U(\Omega_{4}))} \leq \CO ( e^{- \varepsilon / h}) + \CO (h^{1/2}) \Vert & U \CT u \Vert_{H_{\Psi_{t}} (\kappa_U(\Omega_{2}))}   \\
&+ \CO (1) \Vert U \CT u \Vert_{H_{\Psi_{t}} (\kappa_U(\Omega_{2}\setminus\Omega_{4}))} .
\end{aligned}
\end{equation}

On the other hand, applying $\tau_{n}$ to (\ref{az7}), we get also
\begin{equation}  \label{az5}
\tau_{n} (Q-z) \tau_{n} U \CT u + \tau_{n} Q (1-\tau_{n}) U \CT u = \CO ( e^{- \varepsilon / h}).
\end{equation}
Now, if we set
\begin{equation}
\widetilde{Q} = \sum_{j=1}^{d} \lambda_{j} x_{j} h D_{x_{j}} - \frac{ih}{2} \sum_{j=1}^{d} \lambda_{j},
\end{equation}
 we see,  with Proposition 3.3 of \cite{sj87res}, that $Q- \widetilde{Q} = \CO (h^{3/2})$ and $\tau_{n} Q (1-\tau_{n}) = \CO (h^{3/2})$ as operators from $H_{\Psi_{t}} (\kappa_U(\Omega_{2}))$ to  $H_{\Psi_{t}} (\kappa_U(\Omega_{4}))$. Therefore (\ref{az5}) gives
\begin{equation} \label{az6}
\tau_{n} ( \widetilde{Q} - z ) \tau_{n} U \CT u = \CO (e^{- \varepsilon/ h}) + \CO (h^{3/2}) \Vert U \CT u \Vert_{H_{\Psi_{t}} (\kappa_U(\Omega_{2}))} ,
\end{equation}
in $H_{\Psi_{t}} ( \kappa_U(\Omega_{4}) )$. On $\ran \tau_{n}$, in the basis $(x^{\alpha})_{\vert \alpha \vert <n}$, the operator $\tau_{n} ( \widetilde{Q} - z ) \tau_{n}$ reduces to the diagonal matrix with entries $(-hi \sum (\alpha_{j} + 1/2) \lambda_{j} - z )_{\vert \alpha \vert <n}$. So, if $d (z,\Gamma (h)) > \nu h$, $\tau_{n} ( \widetilde{Q} - z ) \tau_{n}$ is invertible on $\ran \tau_{n}$, and its inverse is $\CO  (h^{-1})$. Then (\ref{az6}) gives
\begin{equation} \label{az9}
\tau_{n} U \CT u = \CO (e^{- \varepsilon/ h}) + \CO (h^{1/2}) \Vert U \CT u \Vert_{H_{\Psi_{t}} ( \kappa_U (\Omega_{2}) )} ,
\end{equation}
in $H_{\Psi_{t}} (\kappa_U(\Omega_{4}))$. Adding (\ref{az1}) and (\ref{az9}), we obtain, for $h$ small enough,
\begin{equation}  \label{az3}
\Vert U \CT u \Vert_{H_{\Psi_{t}} (\kappa_U(\Omega_{4}))} \leq \CO (e^{- \varepsilon / h}) + \CO (1) \Vert U \CT u \Vert_{H_{\Psi_{t}} (\kappa_U(\Omega_{2}\setminus \Omega_{4}))} .
\end{equation}
Then, we have, after shrinking $\Omega_4$,
\begin{equation}  \label{aq1}
\Vert \CT u \Vert_{H_{\Phi + tG} ( \Omega_{2} )} \leq \CO (e^{- \varepsilon / h}) + \CO (1) \Vert \CT u \Vert_{H_{\Phi + t G} (\Omega_{2}\setminus \Omega_{4} )} .
\end{equation}

Using the same kind of estimates as in Proposition \ref{hhh}, one can see that
\begin{equation}
\Vert \CT u \Vert_{H_{\Phi + tG} ( \Omega_{1} \setminus \Omega_{4} )} \leq \CO (e^{- \varepsilon / h}) + \CO (h^{1/2}) \Vert \CT u \Vert_{H_{\Phi + t G} (\Omega_{0}  )} .
\label{aq2}
\end{equation}
Adding (\ref{aq1}) and (\ref{aq2}), we obtain 
\begin{equation}
\Vert \CT u \Vert_{H_{\Phi + tG} ( \Omega_{1})} \leq \CO (e^{- \varepsilon / h}) + \CO (h^{1/2}) \Vert \CT u \Vert_{H_{\Phi + t G} (\Omega_{0}\setminus\Omega_1 )}.
\label{kgq}
\end{equation}

On the other hand, from the definition of $G$ (see (\ref{az12})), one can see that there exist $C>0$ and $\varepsilon_{1}>0$, 
 such that
 \begin{equation}  \label{az13}
e^{-tG/h} =
\left\{ \begin{aligned}
{}^{}&\CO (e^{C t /h}) &&\text{ on } \Omega_{0},  \\
{}^{}&\CO (e^{-\varepsilon_{1} t /h}) &&\text{ on }  A_{+}.
\end{aligned} \right.
\end{equation}
Moreover, for each $\varepsilon_{2} > 0$ there exists $\omega \subset \Omega_{1}$, a small enough neighborhood of $0$ such that, in $\omega$, we have
\begin{equation}
e^{-tG/h}\geq  e^{-\varepsilon_{2} t /h}.
 \label{kgb}
\end{equation}
Then (\ref{kgq}) gives
\begin{align}
e^{-\varepsilon_{2} t /2h} \Vert \CT u \Vert_{H_{\Phi} ( \omega )} &\leq \Vert \CT u \Vert_{H_{\Phi + tG} ( \Omega_{1} )}  \nonumber \\
&\leq \CO (e^{- \varepsilon / h}) + \CO (h^{1/2}) \Vert \CT u \Vert_{H_{\Phi + t G} (\Omega_{0} \setminus \Omega_{1} )}    \nonumber  \\
&\leq \CO (e^{- \varepsilon / h}) + \CO (e^{-\varepsilon_{1} t /h}) \Vert \CT u \Vert_{H_{\Phi} ( A_{+} )} + \CO (e^{ t C /h}) \Vert \CT u \Vert_{H_{\Phi} ( A_{-} \cup A_{0} )}   \nonumber  \\
&\leq \CO (e^{- \varepsilon / h}) + \CO (e^{-\varepsilon_{1} t /h})  + \CO (e^{ t C /h} e^{-\delta /h}) ,
\end{align}
since $\Vert \CT u \Vert_{H_{\Phi} (\C^{n})} = \Vert u \Vert_{L^{2} (\R^{n})} \leq 1$. Choosing first $t>0$ small enough and then $\varepsilon_{2}$ small enough, we get
\begin{equation}
\Vert \CT u \Vert_{H_{\Phi} ( \omega )} = \CO (e^{-\widetilde\delta/h}),
\end{equation}
for some $\widetilde\delta >0$, and Theorem \ref{un2} follows.

%----------------------
%
%----------------------
\section{Uniqueness in the $\CC^\infty$ case}
\label{uncinf}

This section is devoted to the proof of Theorem \ref{un1}. Let us recall briefly the assumptions  \eqref{hh1}--\eqref{hh2}: We suppose that $P=\op_{h}(p)$, where $p$ is a  real valued $\CC^\infty$ function, depending on the parameter $h\in ]0,1]$ say, in a  fixed neighborhood of $(0,0)$ in $T^*\R^d$. 
We also assume  that $p$ has an asymptotic expansion with respect to $h$:
\begin{equation}
p(x,\xi,h)\sim \sum_{k\geq 0} p_{k}(x,\xi)h^k,
\label{pas}
\end{equation}
and  that
\begin{equation}
p_{0}(x,\xi)= \sum_{j=1}^d\frac{ \lambda_j}{2} (\xi_j^2 - x_j^2) + \CO((x, \xi )^3),
\end{equation}
where the $\lambda_{j}$'s are real and positive numbers. Finally, we assume that $z \in D(0, C_{0}h)$ for some $C_{0}>0$. 

Recalling the discussion in Section \ref{geoset}, and since $\Lambda_{+}$ and $\Lambda_{-}$ are Lagrangian manifolds,  one can choose local symplectic coordinates $(y,\eta )$ such that
\begin{equation}  \label{aa1}
p_{0} (x, \xi) = B(y, \eta ) y\cdot \eta ,
\end{equation}
where $(y,\eta)\mapsto B(y,\eta)$ is a $\CC^\infty$ mapping from a neighborhood of (0,0) in $T^*\R^d$ to the space $\CM_{d}(\R)$ of $d\times d$ matrices with real entries such that, using the notations of Section \ref{secamr}, 
\begin{equation}
B (0,0) = \left( \begin{array}{ccc}
\lambda_{1} /2 & & \\
& \ddots & \\
& & \lambda_{d} /2
\end{array} \right).
\end{equation}
Now if  $U$ is a  unitary Fourier Integral Operator (FIO) microlocally defined in a neighborhood of $(0,0)$, whose canonical transformation is the map $(x, \xi )\mapsto(y, \eta )$, we denote
\begin{equation}
\widehat{P} = U P U^{-1}.
\end{equation}
Then $\widehat{P}$ is a pseudodifferential operator, with a real (modulo $\CO (h^{\infty})$) symbol $\widehat{p} (y, \eta )= \sum_{j} \widehat{p}_{j} (y, \eta ) h^{j}$, and such that
\begin{equation}
\widehat{p}_{0} = B(y , \eta )y\cdot \eta .
\label{hatp0}
\end{equation}

In order to turn our microlocal problem into a global one, we extend our symbol $p$ as a smooth function on the whole $T^*\R^d$.  Notice that this idea cannot be used in the analytic category. The way we perform this extension is reminiscent of the so-called Complex Absorption Potential Method, used by quantum chemists, and mathematically studied in a paper by P.~Stefanov \cite{st}.

 In the following, the notation $f \prec g$ means that $g=1$ near the support of $f$. Let $\chi_{5}$, $\chi_{8} \in C^{\infty}_{0} ({\rm T}^{*} (\R^{d}))$ be such that the support of $\chi_{8}$ is a small enough neighborhood of $0$ and ${\bf 1}_{\{ 0 \}} \prec \chi_{5} \prec \chi_{8}$. We define
\begin{equation}
\widetilde{p} (y, \eta ) = \widehat{p} (y, \eta ) \chi_{8} (y, \eta ) - i \sqrt{h} ( 1 - \chi_{5} (y, \eta )),
\end{equation}
and we also denote $\tilde P=\op_h(\tilde p)$.
Let us mention that, as one can see following the proof,  one could have taken $h^{\varepsilon}$ with $0<\varepsilon <1$ instead of $\sqrt{h}$ in front of the $1-\chi_{5}$ term.

Now we choose $\chi_{7} \in C_{0}^{\infty} ({\rm T}^{*} \R^{d})$ with $\chi_{5} \prec \chi_{7} \prec \chi_{8}$, and we set
\begin{equation}
g_{1} (y, \eta ) = (y^{2} - \eta^{2}) \chi_{7} (y, \eta ) \ln (1/h).
\end{equation}
Notice that $H_{\hat p_{0}}g_{1}(y,\eta)>0$ for any $(y,\eta)\neq (0,0)$.
Following  the appendix of  \cite{buzw}, we also define 
\begin{equation}
g_{2} (y, \eta ) = \Big( \ln \Big\< \frac{y}{\sqrt{h M}} \Big\> - \ln \Big< \frac{\eta}{\sqrt{h M}} \Big\> \Big) \chi_{3} (y, \eta ) .
\end{equation}
where $M >0$ will be fixed later and $\chi_{3} \in C_{0}^{\infty} ( {\rm T}^{*} (\R^{d}) )$ is such that ${\bf 1}_{\{ 0 \}} \prec \chi_{3} \prec \chi_{5}$. For $t_{1}$, $t_{2} >0$  we set 
\begin{equation}
G_{\pm 1} = \op_{h} (e^{\pm t_{1} g_{1} (y, \eta )} ) \mbox{ and } G_{\pm 2} = \op_{h} (e^{\pm t_{2} g_{2} (y, \eta )}), 
\label{G1G2}
\end{equation}
and we see that 
these $h$-pseudodifferential operators
satisfy $G_{\pm 1} \in \Psi_{h}^{0} (h^{- C t_{1}})$ as well as $G_{\pm 2}\in \Psi_{h}^{1/2} (h^{- C t_{2}})$, for some $C>0$. 

Now, as in  N. Dencker, J. Sj\"{o}strand and M. Zworski \cite[Section 4]{DeSjZw}, or in the very recent paper \cite{chris}, we set
\begin{eqnarray}
\nonumber
&&
Q_{z}= G_{-2} G_{-1} ( \widetilde{P} -z ) G_{1} G_{2}
\\
&&
\quad 
=G_{-2} G_{-1} \left (\op_h (\widehat{p} \chi_{8})+ i\sqrt{h}(1-\chi_{5})+z\right )G_{1} G_{2},
 \label{aa5}
\end{eqnarray}
and we consider each term of the above sum separately.

$\bullet$ First of all, we consider  the operator $G_{-1} \op_{h} \big( i \sqrt{h} (1- \chi_{5}  ) \big) G_{1}$.
By symbolic calculus in the class $\Psi_{h}^0(1)$ (see Proposition \ref{a1a2}), writing $F=\op_{h} \big( (1- \chi_{5}  ) \big)$, we have $F G_{1}=\op_{h}(\varphi_{1})$ with, for any $N_{1}\in \N$,
\begin{align}
\varphi_{1}(x,\xi) =& \sum_{k=0}^{N_{1}} \frac 1{k!}\Big (
\big(\frac{ih}2 \sigma(D_{x},D_{\xi};D_{y},D_{\eta})\big)^k
(1-\chi_{5})(x,\xi)e^{t_{1}g_{1}(y,\eta)}
\Big) \Big|_{y=x,\eta=\xi}   \nonumber  \\
&+ h^{N_{1}-Ct_{1}}\CS^0_h(1).  \label{psi1}
\end{align}
Then again, $G_{-1}\op_{h}(\varphi_{1})=\op_{h}(\varphi_{0})$ with
\begin{align}
\varphi_{0}(x,\xi) =& \sum_{k=0}^{N_{0}} \frac 1{k!}\Big (
\big(\frac{ih}2 \sigma(D_{x},D_{\xi};D_{y},D_{\eta})\big)^k
e^{-t_{1}g_{1}(x,\xi)} \varphi_{1}(y,\eta) \Big) \Big|_{y=x,\eta=\xi}  \nonumber  \\
&+ h^{N_{0}-2Ct_{1}}\CS^0_h(1).  \label{psi0}
\end{align}
But it is easy to see that the $k$-th term in the sum (\ref{psi0}) is $\CO(h^k)$, and choosing $N_{0}\in \N$ such that $N_{0}-2Ct_{1}\geq 0$,  we get that $G_{1}FG_{-1}\in \Psi^0_{h}(1)$.
We also see on (\ref{psi0}) that  the symbol of $G_{1}FG_{-1}$ is supported inside the support of $1-\chi_{5}$ modulo $\CO (h^{\infty})$.

Now since $\chi_{3} \prec \chi_{5}$, we also have, using the same kind of arguments, but in the class $\Psi_{h}^{1/2}(1)$, that 
\begin{equation*}
G_{-2} G_{-1} \op_{h} \big( i \sqrt{h} (1- \chi_{5} (y, \eta ) ) \big) G_{1} G_{2} =G_{-1} \op_{h} \big( i \sqrt{h} (1- \chi_{5} (y, \eta ) ) \big) G_{1} + \CO (h^{\infty}).
\end{equation*}
Notice that without explicit notification, any error term in equalities  between pseudodifferential operator has to be understood in the sense of bounded operators on $L^2$. Finally, keeping only the first term in the expansion (\ref{psi0}), we get
\begin{equation}
G_{-2} G_{-1} \op_{h} \big( i \sqrt{h} (1- \chi_{5} (y, \eta ) ) \big) G_{1} G_{2}= \op_{h} \big( i \sqrt{h} (1- \chi_{5} (y, \eta )) \big) + O \big( h^{3/2} \big).
 \label{aa6}
\end{equation}

$\bullet$ We consider now the second term in (\ref{aa5}), and we set 
\begin{equation}
\widehat{Q} = G_{-1} \op_{h} \big( \widehat{p} (y, \eta ) \chi_{8} (y, \eta ) \big) G_{1}.
\end{equation}
We obtain again by symbolic calculus  in the class $\Psi_{h}^0(1)$ that
\begin{equation}
\widehat{Q} = \op_{h} (\widehat{q}) +  \CO (h^{\infty}),
\end{equation}
where $\widehat{q} (y, \eta ) \in \CS^{0}_{h} (1)$ is supported inside the support of $\chi_{8}$ and satisfies
\begin{equation}
\widehat{q} = \widehat{p}_{0} \chi_{8} + h \widehat{p}_{1} \chi_{8} + ih t_{1} \{ g_{1} , \widehat{p}_{0} \chi_{8} \} +  h^{2} \ln^{2} (1/h)\CS_{h}^{0}(1).
\label{aa3}
\end{equation}

As in \cite{buzw}, since  $G_{\pm 2}$ is in some $\Psi_{h}^{1/2}$, we need to rescale the variables in order to compute the symbol of $G_{-2}\widehat QG_{2}$: We define a unitary transformation $V$ on $L^2(\R^d)$ by
\begin{equation}
V f(y) = \lambda^{-d/2} f \big( \lambda^{-1} y \big), \quad \lambda=\sqrt{hM}, 
\label{scaling1}
\end{equation}
and,  if $a(y, \eta,h)$ is a family of distributions in ${\mathcal S}' ({\rm T}^{*} (\R^{d}))$, we have
\begin{equation}
V^{-1} \op_{h} \big( a (y, \eta ,h) \big) V = \op_{\frac1M} \big( a \big(   \lambda Y ,  \lambda H, \frac{\lambda^2}M\big) \big).
\label{scaling2}
\end{equation}
Notice that here and in what follows,  we always assume that $\lambda \ll 1$. 

Then we set 
$\widetilde{Q} = V^{-1} G_{-2} \widehat{Q} G_{2} V $
and we notice that
\begin{equation}
\widetilde{Q} =
\op_{\frac1M} \big( e^{-t_{2} \widetilde{g}_{2}(Y, H ) } \big) \op_{\frac1M} \big( \widehat{q} \left( \lambda Y , \lambda H \right) \big) \op_{\frac1M} \big( e^{t_{2} \widetilde{g}_{2}(Y, H ) } \big) +  \CO (h^{\infty}),
\end{equation}
where 
\begin{equation} 
\widetilde{g}_{2} (Y,H) = \left( \ln \left\< Y \right\> - \ln \left< H \right\> \right) \chi_{3} \big( \lambda (Y,H) \big).
 \label{aa2}
\end{equation}
We notice that, for any $\alpha,\beta\in \N^d$, and for some constants $C_{\alpha, \beta}$ and $C$ that are independent of  $\lambda$,
\begin{equation} 
\left\vert \partial_{Y}^{\alpha} \partial_{H}^{\beta} e^{\pm t_{2} \widetilde{g}_{2} (Y,H)} \right\vert 
\leq C_{\alpha, \beta} \left\< \frac{\< Y \> }{ \ln \< Y \> }\right\>^{-\vert \alpha \vert}
 \left\< \frac{\< H \> }{\ln \< H \> }\right\>^{-\vert \beta \vert} \< (Y,H) \>^{C t_{2}} ,
 \label{aa4bis}
\end{equation}
Using (\ref{hatp0}), and  since $\lambda\<(Y,H)\>$ can be considered as $\CO(1)$ for $\hat p_{0}\chi_{8}$ is compactly supported, we see  also that, for any $\alpha,\beta\in \N^d$,
\begin{equation} 
\left\vert \partial_{Y}^{\alpha} \partial_{H}^{\beta} (\widehat{p}_{0} \chi_{8}) (\lambda Y, \lambda H ) \right\vert \leq C_{\alpha, \beta} \lambda^{2} \< (Y,H) \>^{2-\vert \alpha \vert - \vert \beta \vert} ,
\label{aa4}
\end{equation}

%%%%%%%%26/11/2004

At this point, it is convenient to  introduce a new class of symbols: We shall write that $f(Y,H,\frac1M)$  belongs to $\widetilde\CS_{\frac1M}(m)$ if it is a smooth function of $(Y,H)$ such that, for any $\alpha,\beta\in \N^d$, there exists a constant $C_{\alpha,\beta}>0$ such that
\begin{equation}
\vert \partial_{Y}^{\alpha} \partial_{H}^{\beta} f (Y,H,\frac 1M) \vert \leq C_{\alpha,\beta}
 \< Y \> ^{-\vert \alpha \vert/2}  \< H \> ^{-\vert \beta \vert/2} 
 m (Y,H).
\end{equation}
Here  the function $m$ is any order function  in the sense of \cite{disj}, Chapter 7 (see also Appendix \ref{pdo}).
With these notations, we have $ e^{-t_{2} \widetilde{g}_{2}}\in \widetilde\CS_{\frac1M}(\bra (Y,H)\ket^{Ct_{2}})$ and $(Y,H)\mapsto\widehat{p}_{0} \chi_{8}(\lambda Y,\lambda H) \in \widetilde\CS_{\frac1M}(\lambda^2\bra (Y,H)\ket^{2})$, uniformly with respect to $\lambda$. 
Notice also that  if $a(y,\eta,h)\in \CS^0_h(1)$ for example, then $ a \big(   \lambda Y ,  \lambda H, \frac{\lambda^2}M\big)\in\widetilde\CS_{\frac1M}(1)$.

Now we compute the symbol of  $\widetilde{Q}$, and we shall again consider each term in (\ref{aa3}) separately.

From $M^{-1}$-pseudodifferential calculus for symbols in $\widetilde\CS_{\frac 1M}$, we get that
\begin{equation}
\op_{\frac1M} \big( e^{-t_{2} \widetilde{g}_{2}} \big)  \op_{\frac1M} (\widehat{p}_{0} \chi_{8}(\lambda Y, \lambda H) ) = \op_{\frac1M}(\widetilde \ell_0),
\end{equation}
where
\begin{align}
\widetilde \ell_0 =&
 e^{-t_{2} \widetilde{g}_{2}} \Big( \widehat{p}_{0} \chi_{8} + i\frac{t_{2}}{2M} \{ \widetilde{g}_{2} , \widehat{p}_{0} \chi_{8} \} \Big) 
 - \frac{1}{8M^2} \Big (
\partial_{Y}^{2} e^{-t_{2} \widetilde{g}_{2}} \partial_{H}^{2} (\widehat{p}_{0} \chi_{8})
 \nonumber  \\
&- 2 \partial^{2}_{Y,H} e^{-t_{2} \widetilde{g}_{2}} \partial^{2}_{Y,H} (\widehat{p}_{0} \chi_{8})
+ \partial_{H}^{2} e^{-t_{2} \widetilde{g}_{2}} \partial_{Y}^{2} (\widehat{p}_{0} \chi_{8}) 
\Big )
\nonumber  \\
&+  e^{-t_{2} \widetilde{g}_{2}} \widetilde\CS_{\frac1M} (M^{-3} \lambda^{2}) + \widetilde\CS_{\frac1M} \big( \lambda^{2} M^{-\infty} \< ( Y,H) \>^{- \infty} \big).
\end{align} 
Notice that we have used (\ref{aa4}) for the first error term above. Using the particular form of $p_{0}$ in (\ref{aa1}) and that of $\widetilde{g}_{2}$ in (\ref{aa2}), and the fact that $\lambda\<(Y,H)\>=\CO(1)$ since $\hat p_0\chi_8$ is compactly supported, we obtain, for some $\varepsilon>0$,
\begin{align}
\tilde \ell_0=&
e^{-t_{2} \widetilde{g}_{2}} \Big( \widehat{p}_{0} \chi_{8} + i\frac{t_{2}}{2M} \{ \widetilde{g}_{2} , \widehat{p}_{0} \chi_{8} \} \Big) 
\nonumber  \\
&+
e^{-t_{2} \widetilde{g}_{2}} \big (\widetilde\CS_{\frac1M} ( t_{2}^{2} \lambda^{2} M^{-2}) +\CO_{M} ( h^{1+\varepsilon})\widetilde\CS_{\frac1M} (1) + \widetilde\CS_{\frac1M} (M^{-3} \lambda^{2}) \big)  
\nonumber  \\
&+ \widetilde\CS_{\frac1M}\big( \lambda^{2} M^{-\infty} \< ( Y,H) \>^{- \infty} \big) .  
\end{align}
Here, the notation $\CO_{M} (m)$ means that the function is bounded by $m$ with bound depending on $M$. Now we compute the symbol $\tilde q_0$ defined by
\begin{equation}
\op_{\frac1M} (\widetilde \ell_0)\op_{\frac1M} \big( e^{t_{2} \widetilde{g}_{2}} \big) =\op_{\frac1M}(\tilde q_0).
\end{equation}
We have
\begin{align}
\widetilde{q}_0 =& \widehat{p}_{0} \chi_{8} + i\frac{t_{2}}{2M} \{ \widetilde{g}_{2} , \widehat{p}_{0} \chi_{8} \} - \frac{it_{2} }{2M} \Big\{ \widehat{p}_{0} \chi_{8} + i\frac{t_{2}}{2M} \{ \widetilde{g}_{2} , \widehat{p}_{0} \chi_{8} \}, \widetilde{g}_{2} \Big\}    \nonumber  \\
&+ \widetilde\CS_{\frac1M} (t_{2}^{2} \lambda^{2} M^{-2}) +\CO_{M}(h^{1+\varepsilon}) \widetilde\CS_{\frac1M}  (1)+ \widetilde\CS_{\frac1M}  (M^{-3} \lambda^{2}) + \widetilde\CS_{\frac1M} \big( \lambda^{2} M^{-\infty} \big)   \nonumber  \\
=& \Big( \widehat{p}_{0} \chi_{8} + i \frac{t_{2}}{M} \{ \widetilde{g}_{2} , \widehat{p}_{0} \chi_{8} \} \Big) + \widetilde\CS_{\frac1M}  ( t_{2}^{2} \lambda^{2} M^{-2}) + O_{M} (h^{1+\varepsilon})\widetilde\CS_{\frac1M}(1) + \widetilde\CS_{\frac1M} (M^{-3} \lambda^{2}). \label{aa7}
\end{align}
Notice that  we have used the following explicit expression:
\begin{equation}   \label{ab3}
\begin{aligned}
\{ \widetilde{g}_{2} , \widehat{p}_{0} \chi_{8} \} =& - \lambda^{2} \chi_{3} \big( \lambda (Y,H) \big) \Big( \frac{H}{\< H \>^{2}}\cdot \big( B(\lambda Y, \lambda H) H  + \partial_{y}B(\lambda Y, \lambda H) \lambda Y\cdot H \big) \\
&+ \frac{Y}{\< Y \>^{2}}\cdot\big( B(\lambda Y, \lambda H) Y + \partial_{\eta}B(\lambda Y, \lambda H) \lambda Y.H \big) \Big)     \\
&+ \lambda^{2} \left( \ln \left\< Y \right\> - \ln \left< H \right\> \right) \left\{ \chi_{3}, \widehat{p}_{0} \chi_{8} \right\} \big( \lambda (Y,H) \big),
\end{aligned}
\end{equation}
so that, in particular,  we have written  $\{ \widetilde{g}_{2} , \{ \widetilde{g}_{2} , \widehat{p}_{0} \chi_{8} \} \} = \CO (\lambda^{2})$  in (\ref{aa7}).

Now we compute the contribution of the second term in (\ref{aa3}). Let us  define the symbol $\tilde \ell_1$ by 
\begin{equation}
\op_{\frac1M} \big( e^{-t_{2} \widetilde{g}_{2}} \big) \op_{\frac1M} \big( (h \widehat{p}_{1} \chi_{8} ) (\lambda Y, \lambda H) \big) = \op_{\frac1M} (\tilde \ell_1).
\end{equation}
We have first
\begin{equation}
\op_h(e^{-t_{2} {g}_{2}})\op_h(h \widehat{p}_{1} \chi_{8} )=\op_h(\ell_1),
\end{equation}
where 
\begin{equation}
\ell_1=h e^{-t_{2} {g}_{2}} \widehat{p}_{1} \chi_{8} +  \frac{h^{2}}{\lambda}e^{-t_{2}{g}_{2}} \CS^{1/2}_h(1) + h^{\infty}\CS^{1/2}_h(1).
\end{equation}
Restoring the $(Y,H)$ variables,  we obtain, also  since $h^2/\lambda\leq \lambda h$,
\begin{equation}
\tilde \ell_1(Y,H)=h e^{-t_{2} \widetilde{g}_{2}} \widehat{p}_{1} \chi_{8}(\lambda Y,\lambda H) +  e^{-t_{2} \widetilde{g}_{2}} \widetilde\CS_{\frac 1M}(\lambda h) + \widetilde\CS_{\frac 1M} ( \lambda^{\infty} ).
\end{equation}
Then, using the symbolic calculus  in the class $\widetilde\CS_{\frac 1M}$, we get
\begin{align}
\op_{\frac1M} (\tilde q_1)  : =& \op_{\frac1M} (\tilde \ell_1)\op_{\frac1M} \big( e^{t_{2} \widetilde{g}_{2}} \big)  \nonumber  \\
=& \op_{\frac1M} ( h  \widehat{p}_{1} \chi_{8}) + \CO (h M^{-2}) + \CO_{M} (h^{1+\varepsilon}).  
\label{aa8}
\end{align}

As for the third term in (\ref{aa3}), we write
\begin{equation}
\op_{\frac1M} \big( e^{-t_{2} \widetilde{g}_{2}} \big) \op_{\frac1M}  \big( (ih t_{1} \{ g_{1} , \widehat{p}_{0} \chi_{8} \}) (\lambda Y, \lambda H) \big)   =\op_{\frac1M}(\tilde \ell_2),
\end{equation}
with
\begin{equation}
\tilde\ell_2= ih t_{1} e^{-t_{2} \widetilde{g}_{2}} \{ g_{1} , \widehat{p}_{0} \chi_{8} \} +  e^{-t_{2} \widetilde{g}_{2}} \lambda h \ln (1/h)\widetilde\CS_{\frac 1M}(1)  + \CO ( h^{\infty} ).
\end{equation}

Then we remark that, for any $\alpha,\beta\in \N^d$, we have, for some $C_{\alpha,\beta}>0$,
\begin{equation}
\big\vert \partial_{Y}^{\alpha} \partial_{H}^{\beta} \big( ih t_{1} \{ g_{1} , \widehat{p}_{0} \chi_{8} \} \big) \big\vert \leq C_{\alpha, \beta} \lambda h \ln (1/h) \< (Y,H) \>^{1-\vert \alpha \vert - \vert \beta \vert},
\end{equation}
so that the function $(Y,H)\mapsto ih t_{1} \{ g_{1} , \widehat{p}_{0}
\chi_{8} \} (\lambda Y,\lambda H)$ belongs to $\widetilde\CS_{\frac 1M}(\lambda h \ln (1/h)\<(Y,H)\>)$. Therefore,  using (\ref{aa4}) and the symbolic calculus in $\widetilde\CS_{\frac 1M}$, we have
\begin{align}
\op_{\frac1M} (\tilde q_2)  : =& \op_{\frac1M}(\tilde \ell_2)\op_{\frac1M} \big( e^{t_{2} \widetilde{g}_{2}} \big) \nonumber  \\
=& \op_{\frac1M} \big( ih t_{1} \{ g_{1} , \widehat{p}_{0} \chi_{8} \}  \big) + \CO(\lambda h \ln (1/h) ).   
\label{aa9}
\end{align}

Finally, let $r(y, \eta)$ be the remainder term in (\ref{aa3}). We see that  in $r\in \CS^{0}_h (h^{3/2})$, and that $r$ has  compact support inside the support of $\chi_{8}$. In the variables $(Y,H)$, we have
\begin{equation}
\big\vert \partial_{Y}^{\alpha} \partial_{H}^{\beta} r (\lambda Y, \lambda H) \big\vert \leq C_{\alpha, \beta} \ h^{2} \ln^2(1/h) \< (Y,H) \>^{-\vert \alpha \vert - \vert \beta \vert}.
\end{equation}
Therefore, working again  in the class $\widetilde\CS_{\frac 1M}$, we obtain
\begin{equation}
\op_{\frac1M} \big( e^{-t_{2} \widetilde{g}_{2}} \big) \op_{\frac1M} \left( r (\lambda Y, \lambda H) \right) \op_{\frac1M} \big( e^{t_{2} \widetilde{g}_{2}} \big) = \CO (h^{2} \ln^2(1/h)).   \label{aa10}
\end{equation}

$\bullet$ It remains to study $G_{-2} G_{-1} z G_{1} G_{2}$.  First of all,  since $\{ e^{-t_{1} \widetilde{g}_{1}} , e^{t_{1} \widetilde{g}_{1}} \} = 0$, we have
\begin{equation}
G_{-1} z G_{1}= \op_{h}(z(1+\CS_{h}^{0}(h^2\ln^2(1/h))).
\end{equation}
Then working in $\widetilde\CS_{\frac 1M}$, we obtain
\begin{equation}
G_{-2} G_{-1} z G_{1} G_{2} = z + \CO (zM^{-2}).
 \label{aa11}
\end{equation}

Finally,  collecting  (\ref{aa5}), (\ref{aa6}), (\ref{aa3}), (\ref{aa7}), (\ref{aa8}), (\ref{aa9}), (\ref{aa10}), we have obtained that
\begin{equation}  \label{ab1}
\begin{aligned}
Q_{z} =& \op_{h} \left( \widehat{p}_{0} \chi_{8} + h \widehat{p}_{1} \chi_{8} + ih t_{1} \{ g_{1} , \widehat{p}_{0} \chi_{8} \} \right)   \\
&+ \op_{h} \left( i t_{2} M^{-1} \{ \widetilde{g}_{2} , \widehat{p}_{0} \chi_{8} \} - i \sqrt{h} (1- \chi_{5} )\right) - z   \\
&+ \CO ( t_{2}^{2} h M^{-1}) + \CO _{M} (h^{1+\varepsilon}) + \CO (h M^{-2} ),  
\end{aligned}
\end{equation}
and we are able to prove the following

\begin{proposition}  \label{ooo}\sl 
Let $\delta$, $C_0  >0$, $t_{1} \gg 1$ and $t_{2} \gg 1$ be fixed. For $M^{-1}$ fixed and $h$ both small enough, we have:
\begin{enumerate}
\item For $z \in D(0, C_0  h)$ and $\im z > \delta h$, the operator $Q_{z} : L^{2} (\R^{d}) \to L^{2} (\R^{d}) $ is invertible and
\begin{equation}  \label{er1}
\Vert Q_{z}^{-1} \Vert = \CO (h^{-1}).
\end{equation}

\item There exists an operator $K = K(h)$ with $\rank K = \CO (1)$ and $K = \CO (1)$ such that $Q_{z} + hK  : L^{2} (\R^{d}) \to L^{2} (\R^{d})$ is invertible for $z \in D(0, C_0  h)$ and
\begin{equation} \label{er2}
\Vert (Q_{z} + h K )^{-1} \Vert = \CO (h^{-1}).
\end{equation}
\end{enumerate}
\end{proposition}

\begin{proof}
For  $u \in \CS (\R^{d})$, we have using (\ref{ab1}),
\begin{align}
- \im \left\< Q_{z} u , u \right\>_{L^{2} (\R^{d})} \geq & - \left\< \left( \op_{h} \Big( h t_{1} \{ g_{1} , \widehat{p}_{0} \chi_{8} \} + \frac{t_{2}}{M} \{ \widetilde{g}_{2} , \widehat{p}_{0} \chi_{8} \} - \sqrt{h} (1- \chi_{5} ) \Big) -\! \im z \right) u,u \right\>     \nonumber  \\
&+ \left(  \CO ( t_{2}^{2} h M^{-1}) +  \CO_{M} (h^{1+\varepsilon}) +  \CO ( h M^{-2} )  \right) \Vert u \Vert^{2}_{L^{2} (\R^{d})}.   \label{ab10}
\end{align}
Let $\chi_{2}$, $\varphi_{1} \in C_{0}^{\infty} ({\rm T}^{*} (\R^{d}); [0,1])$ be such that ${\bf 1}_{\{ 0 \}} \prec \chi_{2} \prec \chi_{3}$,  $\varphi_{1}=0$ near $(0,0)$,  and $(1-\chi_{2}) \chi_{3} \prec \varphi_{1} \prec \chi_{5}$. From (\ref{ab3}),  and since $\chi_{2}$ vanishes on the support of $\{\chi_{3},\widehat p_{0}\chi_{8}\}$, we have, for some $\varepsilon>0$, 
\begin{equation}   \label{ab5}
- \frac1M \{ \widetilde{g}_{2} , \widehat{p}_{0} \chi_{8} \} \chi_{2}^{2} (\lambda Y, \lambda H)
\left\{ \begin{aligned}
&\in \CS_{\frac1M}^{0} (h)  \\
&\geq \varepsilon h \Big( \frac{Y^{2}}{\< Y \>^{2}} + \frac{H^{2}}{\< H \>^{2}} \Big) \chi_{2}^{2} (\lambda Y, \lambda H) ,   \\
\end{aligned} \right.
\end{equation}
if the support of $\chi_{3}$ is small enough. On the other hand, using again the fact that  $\lambda (Y,H)=\CO(1)$, we notice that
\begin{equation}   \label{ab6}
- \frac1M \{ \widetilde{g}_{2} , \widehat{p}_{0} \chi_{8} \} (1- \chi_{2}^{2} ) (\lambda Y, \lambda H) \in \CS_{\frac1M}^{0} \big( h  \ln (1/h)  \big).
\end{equation}
Working in the variables $(Y,H)$, using (\ref{ab5}) and  G{\aa}rding's inequality, we get
\begin{equation}   \label{ab7}
\Big\< \op_{h} \Big(- \frac1M \{ \widetilde{g}_{2} , \widehat{p}_{0} \chi_{8} \} \chi_{2}^{2} \Big) u, u \Big\>  \geq - \frac{C h}{M} \Vert u \Vert^{2}.
\end{equation}
Now, since $1-\varphi_{1}$ and $(1-\chi_{2}^2)\chi_{3}$ have disjoint supports, we have
\begin{equation}  \label{ab7bis}
\begin{aligned}
\op_{h} \Big(- \frac1M \{ \widetilde{g}_{2} , \widehat{p}_{0} & \chi_{8} \} (1- \chi_{2}^{2} ) \Big)   \\
&= \op_{h} (\varphi_{1}) \op_{h} \Big(- \frac1M \{ \widetilde{g}_{2} , \widehat{p}_{0} \chi_{8} \} (1- \chi_{2}^{2}) \Big) \op_{h} (\varphi_{1}) + \CO (h^{\infty}),
\end{aligned}
\end{equation}
and we get from (\ref{ab6})  and Calder\`on--Vaillancourt's theorem, that
\begin{equation}   \label{ab8}
\Big\< \op_{h} \Big(- \frac 1M\{ \widetilde{g}_{2} , \widehat{p}_{0} \chi_{8} \} (1-\chi_{2}^{2}) \Big) u, u \Big\> \geq - C h \ln (1/h) \Vert \op_{h} (\varphi_{1}) u \Vert^{2} + \CO_{M} (h^{\infty}) \Vert u \Vert^{2}.
\end{equation}
Here $C>0$ is uniform with respect to $M$ and $h$.

Let $\chi_{1}$, $\chi_{6} \in C_{0}^{\infty} ({\rm T}^{*} (\R^{d}))$ with ${\bf 1}_{\{0\}} \prec \chi_{1} \prec \chi_{2} \prec \chi_{5} \prec \chi_{6} \prec \chi_{7}$ and $\varphi_{1} \prec \chi_{6} - \chi_{1}$. Then
\begin{align}
-h t_{1} \{ g_{1} , \widehat{p}_{0} \chi_{8} \} &= - t_{1} h \ln (1/h) \left \{ (y^{2} - \eta^{2}) \chi_{7} (y, \eta) , \chi_{8} (y, \eta ) B(y, \eta ) y\cdot\eta  \right \} \nonumber \\
&
\left\{ \begin{aligned}
&\in \CS_{h}^{0} (t_{1} h \ln (1/h) )   \\
&\geq \varepsilon t_{1} h \ln (1/h) ( y^{2} + \eta^{2}) \text{ near the support of } \chi_{6} .  \\
\end{aligned} \right.
\end{align}
Let $\varphi_{2} \in C_{0}^{\infty} ({\rm T}^{*} (\R^{d}); [0,1])$ with $(1- \chi_{6}) \chi_{7} \prec \varphi_{2} \prec (1- \chi_{5} ) \chi_{8}$. Using  G{\aa}rding's inequality for symbols in  $\CS_{h}^{0} ( t_{1} h \ln (1/h))$, we obtain 
\begin{equation}   \label{ab9}
\begin{aligned}
\left\< \op_{h} \left( -h t_{1} \{ g_{1} , \widehat{p}_{0} \chi_{8} \} \chi_{6}^{2}\right)  u, u \right\> \geq \varepsilon t_{1} h \ln (1/h) \Vert \op_{h} ( \chi_{6} & - \chi_{1} ) u \Vert^{2}   \\
&+ \CO (t_1h^{2}\ln(1/h)) \Vert u \Vert^{2} .
\end{aligned}
\end{equation}
We also have, as in (\ref{ab7bis}--\ref{ab8}), 
\begin{equation}   \label{ab11}
\left\< \op_{h} \left( -h t_{1} \{ g_{1} , \widehat{p}_{0} \chi_{8} \} (1- \chi_{6}^{2}) \right)  u, u \right\> \geq - C t_{1} h \ln (1/h) \Vert \op_{h} \left( \varphi_{2} \right) u \Vert^{2} + \CO (h^{\infty}) \Vert u \Vert^{2} .
\end{equation}
Then, collecting (\ref{ab7}), (\ref{ab8}), (\ref{ab9}) and (\ref{ab11}), the inequality (\ref{ab10}) becomes
\begin{equation}
\begin{aligned}  \label{ac2bis}
- \im \left\< Q_{z} u , u \right\> \geq & \varepsilon t_{1} h \ln (1/h) \Vert \op_{h} \left( \chi_{6} - \chi_{1} \right) u \Vert^{2} + \sqrt{h} \left\< \op_{h} (1-\chi_{5}) u,u \right\> + \im z \Vert u \Vert^{2}     \\
&- C h  \ln (1/h)  \Vert \op_{h} (\varphi_{1}) u \Vert^{2} - C t_{1} h \ln (1/h) \Vert \op_{h} \left( \varphi_{2} \right) u \Vert^{2}  \\
&+ \CO ( h M^{-1}) \Vert u \Vert^{2} + \CO_{M} (h^{1+\varepsilon}) \Vert u \Vert^{2},
\end{aligned}
\end{equation}
where $C$ and $\varepsilon$ are uniform with respect to $h$ and $M$. Now, since $\chi_{6}-\chi_{1}=1$ on $\supp \varphi_{1}$, and $1-\chi_{5}=1$ on $\supp\varphi_{2}$, G{\aa}rding's inequality in $\CS_{h}^{0} (\sqrt h)$ gives us,  for any chosen $t_{1}$  large enough 
\begin{equation}  \label{ac3}
- \im \left\< Q_{z} u , u \right\>_{L^{2} (\R^{d})} \geq \im z \Vert u \Vert^{2} + \CO ( h M^{-1}) \Vert u \Vert^{2} + \CO_{M} (h^{1+\varepsilon}) \Vert u \Vert^{2}.
\end{equation}
Then we have
\begin{equation}
- \im \left\< Q_{z} u , u \right\>_{L^{2} (\R^{d})} \geq \frac{\delta h }2 \Vert u \Vert^{2}.
\end{equation}
provided $\im z \geq \delta h$ and $M$ is fixed large enough (and $h$ is small enough).
Since $\vert \im \left\< Q_{z} u , u \right\> \vert \leq \Vert Q_{z} u \Vert \Vert u \Vert$, we get
\begin{equation}
\Vert Q_{z} u \Vert \geq \delta h /2 \Vert u \Vert.
\end{equation}
We can obtain the same way the same estimate for $Q_{z}^{*}$, and this finishes the proof of the first point of the proposition.

Now we consider the second point. Let $\varphi \in C_{0}^{\infty} ({\rm T}^{*} (\R^{d}); [0,1])$ be such that $\varphi = 1$ near $0$. We denote
\begin{equation}
 \label{ac1}
\widetilde{K} = C_{1} \op_{h} \Big( \varphi \Big( \frac{y}{\sqrt{Mh}}, \frac{\eta}{\sqrt{Mh}} \Big) \Big),
\end{equation}
where $C_{1}>0$ is a large constant. Since its symbol is real, $\widetilde{K}$ is self-adjoint and $\widetilde{K} = \CO (C_{1} )$. Recalling (\ref{scaling1}--\ref{scaling2}), we have
\begin{equation}
 \label{ac12}
V^{-1} \widetilde{K} V = C_{1} \op_{\frac1M} \left( \varphi (Y,H) \right),
\end{equation}
and, therefore, $\Vert \widetilde{K} \Vert_{tr} = \CO (C_{1} M^{d})$ (see  \cite[Theorem 9.4]{disj}). Now, using (\ref{ab8}), (\ref{ab9}) and (\ref{ab11}), we get
\begin{equation}   \label{ac2}
\begin{aligned}
- \im \big\< (Q_{z} -ih \widetilde{K}) u , u \big\> \geq & \big\< \op_{\frac1M} \left( - t_{2} M^{-1} \{ \widetilde{g}_{2} , \widehat{p}_{0} \chi_{8} \} \chi_{2}^{2} + C_{1} h \varphi \right) u,u \big\> + \im z \Vert u \Vert^{2}     \\
&+\varepsilon t_{1} h \ln (1/h) \Vert \op_{h} \left( \chi_{6} - \chi_{1} \right) u \Vert^{2} + \sqrt{h} \left\< \op_{h} (1-\chi_{5}) u,u \right\>         \\
&- C h  \ln (1/h) \Vert \op_{h} (\varphi_{1}) u \Vert^{2} - C t_{1} h \ln (1/h) \Vert \op_{h} \left( \varphi_{2} \right) u \Vert^{2}     \\
&+ \CO ( h M^{-1}) \Vert u \Vert^{2} + \CO_{M} (h^{1+\varepsilon}) \Vert u \Vert^{2}.
\end{aligned}
\end{equation}
Here, instead of using (\ref{ab7}),  we notice that, recalling (\ref{ab5}), the term $- \frac{t_{2}}{M} \{ \widetilde{g}_{2} , \widehat{p}_{0} \chi_{8} \} \chi_{2}^{2} + C_{1} h \varphi$ belongs to $ \CS^{0}_{\frac1M} (h)$ and satisfies
\begin{equation}
- t_{2} M^{-1} \{ \widetilde{g}_{2} , \widehat{p}_{0} \chi_{8} \} \chi_{2}^{2} + C_{1} h \varphi \geq \varepsilon \min (t_{2} , C_{1}) h \chi_{2}^{2}.
\end{equation}
Thus, from G{\aa}rding's inequality in $\CS_{\frac1M}^{0} (h)$, we obtain
\begin{equation}
\begin{aligned}
\big\< \op_{\frac1M} \big( - t_{2} M^{-1} \{ \widetilde{g}_{2} , \widehat{p}_{0} \chi_{8} \} & \chi_{2}^{2}+ C_{1} h \varphi \big) u,u \big\>   \\
&\geq \varepsilon \min (t_{2} , C_{1}) h \Vert \op_{h} (\chi_{2}) u \Vert^{2} + \CO ( hM^{-1} ) \Vert u \Vert^{2}.
\end{aligned}
\end{equation}
Now, as in (\ref{ac3}), the inequality (\ref{ac2}) becomes
\begin{align}
- \im \big\< (Q_{z} -ih \widetilde{K}) u , u \big\> \geq & \varepsilon \min (t_{2} , C_{1}) h \Vert \op_{h} (\chi_{2}) u \Vert^{2} + \im z \Vert u \Vert^{2}   \nonumber  \\
&+\varepsilon t_{1} h \ln (1/h) /2 \Vert \op_{h} \left( \chi_{6} - \chi_{1} \right) u \Vert^{2} + \sqrt{h} /2 \left\< \op_{h} (1-\chi_{5}) u,u \right\>       \\
&+ \CO ( h M^{-1}) \Vert u \Vert^{2} + \CO_{M} (h^{1+\varepsilon}) \Vert u \Vert^{2},  \nonumber
\end{align}
for $t_{1}$ large enough. Now, if $t_{1}$, $t_{2}$ and $C_{1}$ are large enough, we get as in (\ref{ac2bis}--\ref{ac3}),
\begin{align}
- \im \big\< (Q_{z} -ih \widetilde{K}) u , u \big\> &\geq 2 C_0  h \Vert u \Vert^{2} + \im z \Vert u \Vert^{2} + \CO ( h M^{-1}) \Vert u \Vert^{2} + \CO_{M} (h^{1+\varepsilon}) \Vert u \Vert^{2}   \nonumber \\
&\geq C_0  h \Vert u \Vert^{2} + \CO ( h M^{-1}) \Vert u \Vert^{2} +\CO_{M} (h^{1+\varepsilon}) \Vert u \Vert^{2},
\end{align}
for $z \in D (0, C_0  h)$. And this implies 
\begin{equation}
\Vert (Q_{z} -i h \widetilde{K} )^{-1} \Vert = \CO (h^{-1}).
\end{equation}
for $M$ large enough. Since $\widetilde{K} = \CO (1)$ is self-adjoint and $\Vert \widetilde{K} \Vert_{tr} = \CO (1)$, one can find a bounded operator $K$ such that   $\rank K$ is finite and independent of $h$, and  $-i \widetilde{K} - K$ is as small as needed (uniformly with respect to $h$, $M$, \ldots), and the proposition is proved.
\end{proof}

Now we  can estimate $(Q-z)^{-1}$ for $z$ away from some discrete set $\Gamma (h)$. We follow J. Sj\"{o}strand \cite{sj01} and S. H. Tang and M. Zworski \cite{tazw}.

\begin{proposition}\sl   \label{ras} Suppose that  the assumptions of Proposition \ref{ooo} hold.
Then there is a discrete set $\Gamma (h)$ independent of $t_{1}$, $t_{2}$ and $M$, with $\# \Gamma (h) \cap D(0, C_0  h)= \CO(1)$, such that $Q_{z} : L^{2} (\R^{d}) \to L^{2} (\R^{d})$ is invertible for $z \notin \Gamma (h) \cap D(0,C_0  h)$. 

Moreover,  if $d ( z, \Gamma (h)) > \nu h^{N}$, for some $\nu >0$ and $N \geq 1$, we have
\begin{equation} \label{er3}
\Vert Q_{z}^{-1} \Vert =\CO (h^{- C}),
\end{equation}
where $C$ depends only on $N$ and $C_0 $.
% where $C = 1+ (N-1) \# \big( \Gamma (h) \cap D(0, C_0  h) \big)$.
\end{proposition}

\begin{proof}
We begin the proof by showing that $Q_{z}$ is invertible outside a finite set $\Gamma (h)$. Here again, we use ideas developed for the study of resonances. Let
\begin{equation}
F(z) = \det \left( Q_{z} (Q_{z} + hK)^{-1} \right) = \det \left( 1 - hK (Q_{z} + hK)^{-1} \right) .
\end{equation}
Since $K$ is trace class, $F(z)$ is well-defined and holomorphic in $D (0, 2 C_0  h)$. From (\ref{er2}), we get
\begin{equation}  \label{ad3}
 F (z)  = \CO (1),
\end{equation}
for $z \in D(0, 2 C_0  h)$. On the other hand, for $\im z > \delta h$, we see from (\ref{er1}) that $(F (z))^{-1} = \det \left( (Q_{z} + hK) Q_{z}^{-1} \right) = \det \left( 1 + hK Q_{z}^{-1} \right)$, so that
\begin{equation}  \label{ad4}
\vert F (z) \vert > \varepsilon ,
\end{equation}
still for $\im z > \delta h$. The estimates (\ref{ad3}), (\ref{ad4}) and Jensen's formula imply that the number of zeros of $F(z)$ in $D(0, C_0  h)$ is bounded. Using the properties of the determinant of an analytic family of operators, we get that $Q_{z}$ is invertible outside a bounded set $\Gamma (h)$ and that the algebraic multiplicity of the poles of $Q_{z}^{-1}$ is bounded. At this point, $\Gamma (h)$ depends on $t_{1}$, $t_{2}$ \ldots, but if $h$ and $M^{-1}$ are small enough, we have
\begin{alignat}{2}
&G_{1} G_{-1} = 1 + \CO (h^{2}\ln^2\frac 1h),          && G_{-1} G_{1} = 1 + \CO (h^{2}\ln^2\frac 1h), \\
&G_{2} G_{-2} = 1 + \CO (M^{-2}), \qquad   && G_{-2} G_{2} = 1 + \CO (M^{-2}),
\end{alignat}
so that $G_{1}$, $G_{-1}$, $G_{2}$ and $G_{-2}$ are invertible. Thus, $\Gamma (h)$ is nothing but the set of eigenvalues of $\widetilde{P}$, which are independent of $t_{1}$, $t_{2}$ and $M$. These eigenvalues have finite multiplicity.

In order to estimate $\Vert Q_{z}^{-1}\Vert$ for $z$ away from $\Gamma (h)$, we use the same strategy as in \cite{sj01}. Let $e_{1},$ $\ldots, e_{N}$ be an orthonormal basis of $\im  K^{*} = (\ker K)^{\bot}$ and $(e_{j})_{j\geq N+1}$ an orthonormal basis of $\ker K$. We denote $R_{+} : L^{2} (\R^{d})\to \C^{N}$ and $R_{-} : \C^{N} \to L^{2} (\R^{d})$ the operators given by
\begin{equation}
R_{+} (u)= (\< u, e_{j} \>)_{j=1\dots N}, \qquad R_{-} {u_{-}}_{j} = \sum_{j=1}^{N} {u_{-}}_{j} (Q_{z} + hK)e_{j}.
\end{equation}
We study the following operator on $L^{2}(\R^{d}) \times \C^{N}$
\begin{equation}
{\mathcal P}_{z} = \left( \begin{array}{cc}
Q_{z} & R_{-} \\
R_{+} & 0
\end{array} \right) ,
\end{equation}
which is associated to the Grushin problem
\begin{equation}  \label{ad1}
\left\{ \begin{aligned}
&Q_{z} u + R_{-} u_{-} &&= v \\
&R_{+} u  &&= v_{+}
\end{aligned} \right.  ,
\end{equation}
with $u$, $v \in L^{2}(\R^{d})$ and $u_{-}$, $v_{+} \in \C^{N}$. Since $Q_{z} = Q_{z} + h K -hK$ with $Q_{z} + hK$ invertible and $K$ compact, $Q_{z}$ and then ${\mathcal P}_{z}$ are holomorphic families of Fredholm operators of index $0$. It is therefore enough  to show that $\CP$ is injective to show that it is invertible. Assume that
\begin{equation}
{\mathcal P} \left( \begin{aligned}
&u \\
&u_{-}
\end{aligned} \right) =0 ,
\end{equation}
with $u = \sum_{j=1}^{\infty} u_{j} e_{j}$. Then, since $R_{+} u =0$, we get $u_{j} = 0$ for $1 \leq j \leq N$, and, since $K e_{j} = 0$ for $N < j $, the equation $Q_{z} u + R_{-} u_{-} = 0$ becomes
\begin{equation}
(Q_{z} +hK) \Big( \sum_{j=N+1}^{\infty} u_{j} e_{j} + \sum_{j=1}^{N} {u_{-}}_{j} e_{j} \Big) =0.
\end{equation}
Then, from (\ref{er2}), we get $u=0$, $u_{-} =0$, and ${\mathcal P}$ is invertible. We denote its inverse by
\begin{equation}
{\mathcal P}^{-1} = \left( \begin{array}{cc}
E (z) & E_{+} (z)  \\
E_{-} (z) & E_{-+} (z)
\end{array} \right) ,
\end{equation}
and we look for estimates on the entries of $\CP^{-1}$. Assume (\ref{ad1}) and write $u = (u',u'')$ with $u' \in {\rm vect} \{ e_{1}, \ldots,e_{N} \}$ and $u'' \in {\rm vect} \{ e_{N+1}, \ldots \}$. Since $Q_{z} u + R_{-} u_{-} = v$, we have
\begin{equation}
(Q_{z} +hK) \bigg( u'' + \sum_{j=1}^{N} {u_{-}}_{j} \bigg) = v - Q_{z} u' ,
\end{equation}
and 
\begin{equation}
u'' + \sum_{j=1}^{N} {u_{-}}_{j} = (Q_{z} +hK)^{-1} v - \big( 1- (Q_{z} +hK)^{-1} h K \big) u' .
\end{equation}
Therefore, since $Q_{z} = \CO (1)$ and using (\ref{er2}), we obtain
\begin{equation}
\Vert u'' \Vert_{L^{2}} + \Vert u_{-} \Vert_{\C^{N}} \leq C \big( h^{-1} \Vert v \Vert_{L^{2}} + \Vert u' \Vert_{L^{2}} \big),
\end{equation}
and then, since  we have
$
\Vert u'  \Vert_{L^{2}} = \Vert v_{+} \Vert_{\C^{N}}
$ because $R_{+} u  = v_{+}$, we get
\begin{equation} \label{ad2}
\Vert u \Vert_{L^{2}} + \Vert u_{-} \Vert_{\C^{N}} \leq C \big( h^{-1} \Vert v \Vert_{L^{2}} + \Vert v_{+} \Vert_{\C^{N}} \big).
\end{equation}
Thus we have $E = \CO (h^{-1})$, $E_{-} =  \CO (h^{-1})$, $E_{+} =  \CO(1)$ and $E_{-+} =  \CO (1)$.

Now we follow S. H. Tang and M. Zworski \cite{tazw}. For $z \in D(0, 2 C_0  h)$, $Q_{z}$ is invertible if and only if  $E_{-+} (z)$ is invertible, and in that case,
\begin{equation}
Q_{z}^{-1} = E_{+} (z) E_{-+ }^{-1} (z) E_{-} (z) - E(z),
\end{equation}
which implies
\begin{equation}
Q_{z}^{-1} =  \CO (h^{-1}) \big( 1+ \Vert E_{-+}^{-1} (z) \Vert \big).
\label{schurnorm}
\end{equation}
Since $E_{-+} (z) =  \CO (1)$ as an operator on $\C^{N}$, we also have 
$
\Vert E_{-+}^{-1} (z) \Vert =  \CO (\vert D(z) \vert^{-1})
$, where $D(z)=\det  E_{-+ } (z)$. Now we set
\begin{equation}
D_{w} (z) = \prod_{z_{j} \in \Gamma (h) \cap D(0, C_0  h)} \frac{z-z_{j}}{h},
\end{equation}
and  we know that $D(z) = D_{w} (z) \times G (z)$, where $G(z)$ is holomorphic. Here, we use the fact that the order of the zeros of $D (z;h)$ coincides with the multiplicity of the eigenvalues of $\widetilde{P}$. Since $\# \big (\Gamma (h)\cap D(0, C_0  h)\big )$ is uniformly bounded, we have for $z \in D(0, C_0  h)$,
\begin{equation} \label{ae2}
D_{w} (z) =  \CO(1).
\end{equation}
 On the other hand, one can find $r (h) \in ]C_0 , 2 C_0 [$ such that, for $z$ on the circle $\partial D(0,r(h))$, we have
\begin{equation}  \label{ae3}
\vert D_{w} (z) \vert \geq \varepsilon .
\end{equation}
For $z \in D(0, 2 C_0  h)$ we also have
\begin{equation}  \label{ae4}
D (z) = \CO (1),
\end{equation}
since $E_{-+} =  \CO (1)$ as an operator on $\C^{N}$. Finally if $\im z > \delta h$, we have $E_{-+}^{-1} (z) = -R_{+} Q_{z}^{-1} R_{-}$, thus
\begin{equation}  \label{ae1}
\vert D (z) \vert > \varepsilon.
\end{equation}
 Using (\ref{ae3}) and (\ref{ae4}), we obtain
\begin{equation}
G(z) =  \CO(1)
\end{equation}
for $z \in D(0, r(h))$. Now  (\ref{ae2}) and (\ref{ae1}) imply that
\begin{equation}
\vert G(z) \vert > \varepsilon,
\end{equation}
for $\im z > \delta h$. Then  Harnack's inequality for the function $C - \ln \vert G(z) \vert$, where $C$ is chosen so that this function is non-negative,  implies
\begin{equation}
G(z)^{-1} =  \CO (1)
\end{equation}
for $z \in D(0, C_0  h)$. Therefore, if $d(z,\Gamma (h)) > \nu h^{N}$, one has  $\det (E_{-+} (z) )^{-1} =  \CO(h^{-C})$,  $(E_{-+} (z) )^{-1} =  \CO(h^{-C})$ and the proposition follows from  (\ref{schurnorm}).
\end{proof}

Finally, we extend the domain of validity of the estimate (\ref{er1}) on $Q_{z}^{-1}$ as far as possible into the lower half complex plane.

\begin{proposition}\sl   \label{super}
Assume that $t_{2}$ and $M^{-1}$ are small enough. There exists $\delta_{0} >0$ such that, for all $t_{1}$ large enough, $Q_{z}$ is invertible on $L^{2} (\R^{d})$ for $z \in D (0, C_0  h)$ with $\im z > - \delta_{0} h$, and, for such $z$,
\begin{equation}
\Vert Q_{z}^{-1} \Vert = \CO (h^{-1}).
\end{equation}
\end{proposition}

\begin{proof}
From (\ref{ab5}) and  Fefferman-Phong's inequality, we have
\begin{align*}
\left\< \op_{h} \left(- M^{-1} \{ \widetilde{g}_{2} , \widehat{p}_{0} \chi_{8} \} \chi_{2}^{2} \right) u, u \right\>  \geq  \varepsilon h \op_{\frac1M} \Big( \Big( \frac{Y^{2}}{\< Y \>^{2}} + \frac{H^{2}}{\< H \>^{2}} \Big) \chi_{2}^{2} (\lambda Y, \lambda H) \Big) + \CO (h M^{-2}) \Vert u \Vert^{2}.
\end{align*}
But using the Appendix of \cite{buzw}, we have
\begin{equation}
\op_{\frac1M} \Big( \Big( \frac{Y^{2}}{\< Y \>^{2}} + \frac{H^{2}}{\< H \>^{2}} \Big) \chi_{2}^{2} (\lambda Y, \lambda H) \Big) \geq \varepsilon M^{-1} \chi_{1}^{2} (\lambda Y, \lambda H) + \CO (M^{-2}),
\end{equation}
and (\ref{ab10}) becomes, as for (\ref{ac2bis}),
\begin{align}
\im \left\< Q_{z} u , u \right\>_{L^{2} (\R^{d})} \geq& \varepsilon t_{2} h M^{-1} \Vert \op_{h} (\chi_{1}) u \Vert^{2} + \varepsilon t_{1} h \ln (1/h) \Vert \op_{h} \left( \chi_{6} - \chi_{1} \right) u \Vert^{2}  \nonumber  \\
&+ \sqrt{h} \left\< \op_{h} (1-\chi_{5}) u,u \right\> + \im z \Vert u \Vert^{2}    \nonumber   \\
&- C h (\ln (1/h) + \ln (M)) \Vert \op_{h} (\varphi_{1}) u \Vert^{2} - C t_{1} h \ln (1/h) \Vert \op_{h} \left( \varphi_{2} \right) u \Vert^{2} \nonumber \\
&+ \left( \CO ( t_{2}^{2} h M^{-1}) + \CO_{M} (h^{1+\varepsilon}) + \CO ( h M^{-2} ) + \CO (zM^{-2}) \right) \Vert u \Vert^{2} .
\end{align}
If $t_{2}$ is fixed small enough and $t_{1}$ large enough, we obtain
\begin{align}
- \im \left\< Q_{z} u , u \right\>_{L^{2} (\R^{d})} \geq& 
\big( \varepsilon t_{2} h M^{-1} + \im z \big) \Vert u \Vert^{2}  
  \nonumber  \\ 
&+ \left( \CO ( t_{2}^{2} h M^{-1}) + \CO_{M} (h^{1+\varepsilon}) + \CO ( h M^{-2} ) + \CO (hM^{-2}) \right) \Vert u \Vert^{2} .
\end{align}
This gives the Proposition, provided  $M$ is chosen large enough.
\end{proof}

Now we finish the proof of Theorem \ref{un1}. Proposition \ref{ras} implies that, for $u \in L^{2} (\R^{d})$ and if $d ( z, \Gamma (h)) > \nu h^{N}$, we have
\begin{equation}  \label{ae7}
\Vert G_{2}^{-1} G_{1}^{-1} u \Vert = \CO (h^{-C}) \Vert G_{-2} G_{-1} (\widetilde{P} -z) u \Vert .
\end{equation}
Since $G_{-2} \in \Psi_{h}^{1/2} (h^{-C t_{2}})$, we also have  $\Vert G_{-2} \Vert = \CO (h^{-C t_{2}})$. Working in $\widetilde\CS_{\frac 1M}$, we get $\op_{h} \big( e^{-t_{2} g_{2}} \big) \op_{h} \big( e^{t_{2} g_{2}} \big) = 1+ \CO (M^{-2})$, so that
\begin{equation}
G_{2}^{-1} = (1+ \CO (M^{-2})) G_{-2} = \CO (h^{-C t_{2}}).
\end{equation}
Then (\ref{ae7}) becomes
\begin{equation}
\Vert G_{1}^{-1} u \Vert = \CO (h^{-C-2C t_{2}}) \Vert G_{-1} (\widetilde{P} -z) u \Vert,
\end{equation}
and since, from pseudodifferential calculus in $\CS_{h}^{0}$, we have
% $\op_{h} \big( e^{-t_{1} g_{1}} \big) \op_{h} \big( e^{t_{1} g_{1}} \big) = 1+ O (h^{2})$ and 
$G_{1}^{-1} = (1+\CO (h^{2})) G_{-1}$, this gives
\begin{equation}  \label{ae8}
\Vert G_{-1} u \Vert = \CO (h^{-C-2C t_{2}}) \Vert G_{-1} (\widetilde{P} -z) u \Vert .
\end{equation}
Now we use a Cordoba--Fefferman type estimate, as given by A. Martinez (see \cite[Corollary 3.5.3]{mabk}) and in a more precise form in \cite[Theorem 3]{bomi}.

\begin{lemma} \sl \label{mama}
Let $f(y,\eta)$, $a(y,\eta,h) \in \CS^0_{h}(1)$. There exist a symbol $b(y,\eta ,h)\sim\sum_{j\geq 0}h^jb_{j}(y,\eta)\in \CS^{0}_{h}(1)$, and an operator $R(h)=  \CO(h^\infty)$ such that, for all 
$u,v\in \CC^\infty_{0}(\R^d)$, one has
\begin{equation}
\<f\CT'\op_{h}(a)u,  \CT' v\>_{L^{2} (\R^{d})} = \big\< (b(y,\eta,h)+R(h)) \CT' u, \CT' u \big\>_{L^{2} ({\rm T}^{*} (\R^{d}))} ,
\end{equation}
where $\supp b_{j}\subset\supp f$ for all $j$, and $b_j$ is given in terms of derivatives of $a$ and $f$ of order at most $2j$. In particular
\begin{equation*}
b_{0}(y,\eta)=f(y,\eta)a_{0}(y,\eta).
\end{equation*}
\end{lemma}

From this result, one can obtain the following

\begin{corollary}\sl 
There exists a function $b (y, \eta,h ) = 1+ \CO_{t_{1}} (h \ln^{2} (1/h))$ such that
\begin{equation}
\Vert G_{-1} u \Vert^{2}_{L^{2} (\R^{d})} = \big\< b e^{- t_{1} g_{1} (y, \eta )} \CT' u, e^{- t_{1} g_{1} (y, \eta )} \CT' u \big\>_{L^{2} ({\rm T}^{*} (\R^{d}))} + \CO (h^{\infty}) \Vert u \Vert^{2}.
\end{equation}
\end{corollary}

Then, in view of this estimate,  (\ref{ae8}) gives
\begin{equation}  \label{ae9}
\Vert e^{- t_{1} g_{1}} \CT' u \Vert_{L^{2} ({\rm T}^{*} (\R^{d}))} = \CO (h^{-C-2C t_{2}}) \Vert e^{- t_{1} g_{1}} \CT' (\widetilde{P} -z) u \Vert_{L^{2} ({\rm T}^{*} (\R^{d}))} + \CO (h^{\infty}) \Vert u \Vert .
\end{equation}

Now assume that $u$ and $z$ satisfy the assumptions of Theorem \ref{un1}. Let $\chi_{4} \in C_{0}^{\infty} ({\rm T}^{*} (\R^{d}))$ with $\chi_{3} \prec \chi_{4} \prec \chi_{5}$ and suppose that $\FS((\widetilde P-z) u)$ does not intersect a neighborhood of the support of $\chi_{4}$. Then (\ref{ae9}) gives
\begin{align} 
\Vert e^{- t_{1} g_{1}} \CT' \op_{h}(\chi_{4}) u \Vert =& \CO (h^{-C-2C t_{2}}) \Vert e^{- t_{1} g_{1}} \CT' \op_{h}(\chi_{4}) (\widetilde{P} -z) u \Vert  \nonumber  \\
&+ \CO (h^{-C-2C t_{2}}) \Vert e^{- t_{1} g_{1}}\CT' [ \widetilde{P} , \op_{h}(\chi_{4}) ] u \Vert + \CO (h^{\infty}) \Vert u \Vert .
\end{align}
Let $\varphi_{3} \in C_{0}^{\infty} (T^{*} (\R^{d}))$ with $\chi_{4} ' \prec \varphi_{3} \prec (1- \chi_{3} ) \chi_{5}$. Since the operator $[ \widetilde{P} , \op_{h}(\chi_{4}) ] \in \Psi_{h}^{0} (h)$ has its symbol supported inside the support of $\chi_{4} '$, we get, using again Lemma \ref{mama},
\begin{equation}
\Vert e^{- t_{1} g_{1}} \chi_{4} \CT'u \Vert_{L^{2} ({\rm T}^{*} (\R^{d}))} = \CO (h^{\infty}) + \CO (h^{-C-2C t_{2}}) \Vert e^{- t_{1} g_{1}} \varphi_{3} \CT' u \Vert  .
\end{equation}
Here the constant $C$ is uniform with respect to $t_{1}$. Now, we can assume that the support of $\chi_{4}$ and $\varphi_{3}$ satisfies the  properties of the domains $\Omega_1$ and $\Omega_0 \setminus {\Omega_1}$ as in (\ref{omegas})  and Figure \ref{figdomain},  and we get  the main part of Theorem \ref{un1}.
The remaining statement concerning the fact that  the exceptional set $\Gamma(h)$ can be chosen so that $\Gamma(h)\subset \{\im z\leq - \delta_{0} h\}$ for some $\delta_{0}>0$, follows from the above discussion,  using Proposition \ref{super} instead of Proposition \ref{ras}.

%------------------------------------
%
%------------------------------------

\section{Existence}
\label{secexis}

This section is devoted to the proof of Theorem \ref{exis}. We use the ideas and the constructions of B.~Helffer and J.~Sj\"{o}strand in \cite{hsmw3}, concerning the study of the
tunnel effect between potential wells. At many places in the following sections, we shall use some terminology and some general results from  \cite{hsmw3} that we recall now, here in a  slightly different setting.

Let $(\mu_{j})_{j\geq 0}$ be the strictly growing  sequence of linear 
combinations over $\N$ of the $\lambda_{j}$'s. Let $u (t,x,\eta ')$ be a function defined on $\R^+\times U \times V$, $U \subset\R^d$, $V\subset \R^{m}$.

\begin{definition}\label{expdeb}\sl We say that $u:[0,+\infty[\times U \times V \to \R$, a smooth function, is expandible, if, for any $N\in\N$, $\varepsilon >0$, $\alpha , \beta , \gamma \in \N^{1 + d + m}$,
\begin{equation}
\partial^{\alpha}_{t} \partial_{x}^\beta \partial_{\eta '}^\gamma \Big( u(t,x,\eta ')-\sum_{j= 1}^N 
u_{j}(t,x,\eta ') e^{-\mu_{j}t} \Big)= \CO \big( e^{-(\mu_{N+1} - \varepsilon ) t} \big) ,
\label{defexpandible2}
\end{equation}
for a sequence of  $(u_{j})_{j}$ smooth functions, which are 
polynomials in $t$.
We shall write 
\begin{equation*}
u(t,x,\eta ' ) \sim \sum_{j\geq 1} u_{j}(t,x,\eta ') e^{-\mu_{j}t} ,
\end{equation*}
when (\ref{defexpandible2}) holds.
\label{def1}
\end{definition}

\noindent
As the following result shows, this symbol class is the suitable one for our geometric setting at $(0,0)$.

\begin{proposition}[\cite{hsmw3}, Section 3]  \label{p22}
Let $\nu(t,x , \eta ' )$ be a time-dependent vector field. Suppose that there exists a matrix-valued map $(x,\eta ')\mapsto A(x,\eta ')$ from $U \times V$ to 
$\CM_d(\R)$ such that
\begin{enumerate}
\item $A(0)=\diag(\lambda_1,\lambda_2,\dots, \lambda_d)$, with 
$0<\lambda_1\leq\lambda_2\leq \dots\leq\lambda_d$.
\item $(t,x,\eta')\mapsto  \nu(t,x,\eta')-A(x,\eta')x$ is a smooth real expandible 
matrix.
\end{enumerate}
Then, if $v(t,x,\eta ' )$ is expandible and vanishes at $x=0$,  and $u_{0}(x,\eta')$ is a smooth function,  the solution $u(t,x, \eta ' )$ to the Cauchy problem
\begin{equation}
\left\{
\begin{aligned}
&\partial_t u + \nu (t,x, \eta ' )u=v, \quad  t\geq 0, \ x \in U , \ \eta ' \in V, \\
&u_{\vert_{t=0}}=u_{0},
\end{aligned}
\right .
\label{evol}
\end{equation}
is expandible.
\end{proposition}

Notice that this result implies in particular that, as we have already mentioned in Section \ref{secamr},  $\gamma(t,x,\xi)=\exp(\pm tH_p)(x,\xi)$ is expandible when $(x,\xi)\in \Lambda_\mp$.

\begin{definition}\sl We say that $u (t,x, \eta ' ,h)$, a smooth function is of class $\CS^{A,B}$ if, for any $\varepsilon >0$, $(\alpha , \beta , \gamma) \in \N^{1 + d + m}$,
\begin{equation} \label{defexpandible3}
\partial^{\alpha}_{t} \partial_{x}^\beta \partial_{\eta '}^\gamma u(t,x,\eta ' , h) = \CO \big( h^{A} e^{-(B - \varepsilon ) t} \big) .
\end{equation}
Let $\CS^{\infty , B} = \bigcap_{A} \CS^{A,B}$. We say that $u (t,x, \eta ' ,h)$ is a classical expandible function of order $(A,B)$, if, for any $K \in \N$,
\begin{equation}
u(t,x,\eta ' ,h) - \sum_{k=A}^{K} u_{k} (t,x,\eta ') h^{k} \in \CS^{K+1,B},
\end{equation}
for a sequence of  $(u_{k})_{k}$ expandible functions. We shall write 
\begin{equation*}
u(t,x,\eta ' ,h) \sim \sum_{k\geq A} u_{k}(t,x,\eta ') h^{k} ,
\end{equation*}
in that case.
\end{definition}

We recall from Section \ref{secamr} that $\Omega$ is a small neighborhood of $(0,0) \in T^{*} \R^{d}$, that $\varepsilon >0$ is small enough such that $S = \Lambda_{-} \cap \{ (x, \xi ); \ \vert x \vert = \varepsilon\} \subset \Omega$, and  $U \subset \Omega$ a neighborhood of $S$. We look for a solution of the problem
\begin{equation}\label{zeproblem}
\left\{ \begin{aligned}
&(P-z) u = 0 &&\text{ microlocally in } \Omega ,   \\
&u = u_0 &&\text{ microlocally in } U .
\end{aligned} \right.
\end{equation}
Since this problem 
is linear with respect to the  initial data
$u_{0}$, we can  assume that 
$u_{0}$ vanishes  microlocally outside a small neighborhood of some point
${\rho}_{-} \in (\Lambda_{-} \cap S) \setminus
\widetilde{\Lambda_{-}}$. We recall that by assumption, $u_0$ vanishes on $\widetilde{\Lambda_{-}}$. Since $P$ is of principal type in $\Omega\setminus \{(0,0)\}$, 
$u_{0}$ can be extended as a microlocal solution of $(P-z) u_{0} = 0$
near each point of $\Lambda_{-} \setminus \{ (0,0) \}$. As
${\rho}_{-} \notin \widetilde{\Lambda_{-}}$, we know from
(\ref{eqfuwsd}), that 
\begin{equation}  \label{uio}
\gamma_{-} (t) = \exp (t H_{p}) ({\rho}_{-} ) \sim \sum_{j=1}^{+\infty}
e^{-\mu_{j} t} \gamma^-_{j} (t), \text{ as } t \to + \infty,
\end{equation}
where $\gamma^-_{1} \neq 0$ is an eigenvector of $F_{p}=d_{(0,0)}H_p$ associated to
the eigenvalue $-\lambda_{1}$. We recall that  $(\mu_{j})_{j\geq 0}$ is  the strictly growing  sequence of linear 
combinations over $\N$ of the $\lambda_{j}$'s.

Here and from now on, we shall write points in $T^*\R^d$ as
$(x,\xi)=(x_{1},x', \xi_{1},\xi')$ with $x_{1}$, $\xi_{1}$ in $\R$ and
$x'$, $\xi'$ in $\R^{d-1}$. We can always assume, up to a linear
change of variables, that $g_1^-(\rho_-)=\Pi_{x} \gamma^-_{1}$ is collinear to the direction
$x_{1}$. In these coordinates, we set $H_{-}: x_{1}=\varepsilon$. Of
course, the lift $H_-\times \R^d$ of $H_{-}$ in $T^*\R^d$  is
transverse to $\gamma_{-}$ for
$\varepsilon$ small enough, and we can suppose so.
Here and in the sequel we may have to change a certain finite number of times for a smaller $\varepsilon>0$, and  therefore to change (silently) for another $\rho_-$ on the curve $\gamma_-$.
In the rest of this section, we prove  Theorem \ref{exis}  under a more precise form. As in \cite{hsmw3}, the main idea is  to look for a solution to (\ref{zeproblem}) of the form
\begin{equation}  \label{solution}
u(x,h) = \frac{1}{(2 \pi h)^{d-\frac{1}{2}}} \iint_{T^*\R^{d-1}} \int_{-1}^{+\infty} 
e^{i( \varphi (t,x,\eta ' )- y' \eta ' )/h} a(t,x, \eta ' ,z, h) u_{0}(\varepsilon , y' )
dt dy' d\eta ' .
\end{equation}
Therefore we shall  look for a phase  function   $\varphi$ and  a symbol $a$ such that
\begin{equation}\label{mastereq}
(h D_{t} + P(x,hD) - z) a e^{i\varphi/h} = {\mathcal O} (h^{\infty}),
\end{equation}
in a sense that we will precise later on.
One of the differences with respect to \cite{hsmw3} is that we shall do so for  each $\eta'$ in  a neighborhood of $\xi^-{}'$,  so that we can also fulfill the initial condition in (\ref{zeproblem}).

However, as in \cite{hsmw3}, in general, the integral with respect to $t$ in (\ref{solution}) does not converge for the functions $a$ and $\varphi$ we build, and our representation of the solution is somewhat more complicated than (\ref{solution}). Recalling that we suppose $z\in [-{{C_0}}h , {{C_0}}h] + i [- {{C_1}} h , {{C_1}}h]$ for some ${{C_0}}$, ${{C_1}} >0$, we denote
\begin{equation}
S=S(z/h)=\sum_{j=1}^d \frac{\lambda_j}2-i\frac zh, \mbox{ and }
K_1={\E} \Big( \frac{{{C_1}}}{\lambda_{1}} - \frac{\sum \lambda_{j}}{2 \lambda_{1}} \Big) +1,
\end{equation}
where  ${\mathbb E} (r)$ denotes the integer part of $r\in \R$.

\begin{theorem}\label{repres}
\sl 
Assume that $u_{0}$ vanishes microlocally  in $\Lambda_{-} \cap (H_{-}
\times \R^{d})$
outside a small neighborhood of $\rho_{-}$. Then, there
exist a neighborhood $U$ (resp.  $W$) of $\gamma_{-} ([-1 , + \infty[) \cup \{ 0
\}$ (resp.  $\xi^{-} {}'$) in $\R^{d}$ (resp. $\R^{d-1}$), a phase
function $\varphi(t,x,\eta ')$, a symbol $A_+(t,x, \eta ' ,z,h)$ defined
on $[-1,+\infty[\times U \times W$, and  a symbol $A_-(x,\eta',z,h)$ defined
on $U \times W$ such that
\begin{enumerate}
\item  There exists a smooth function $\tilde\psi(\eta')$ such that  the function  $\varphi-\varphi_+(x)-\tilde\psi (\eta')$ is expandible: 
\begin{equation*}
\varphi (t,x,\eta')-(\varphi_+(x)+\tilde\psi(\eta'))\sim \sum_{j\geq 
1}e^{-\mu_j t}\varphi_j(t,x,\eta').
\end{equation*}
Moreover $\tilde\psi$ is a generating function for $\Lambda_-$, in the 
sense that, the projection of $\Lambda_-$ onto $T^*H_-$ can be written as the set of 
$(\nabla\tilde\psi(\eta'), \eta')$'s, with $\eta'\in W$. \medskip

\item The symbol $A_+$ is classically expandible: $A_+\in \CS^{-K_1,-\delta}$ for some $\delta>0$, and it is an analytic function with respect to $z\in [-{{C_0}}h , {{C_0}}h] + i [- {{C_1}} h , {{C_1}}h]$. \medskip

\item The function  $A_-$ is a semiclassical symbol of order $-K_1$, and it is an analytic function with respect to $z\in [-{{C_0}}h , {{C_0}}h] + i [- {{C_1}} h , {{C_1}}h]$. \medskip

\item For any cut-off function $\chi \in C^{\infty} (]-1, + \infty[)$ equal to $1$ near $[0, + \infty[$, the function
\begin{align}
u(x, & z,h)= \frac{1}{(2 \pi h)^{d-\frac{1}{2}}} \iint_{T^*\R^{d-1}}  e^{i( \varphi_+(x)+\tilde\psi(\eta ' )- y' \eta ' )/h}A_{-} (x, \eta ' ,z ,h) u_0(\varepsilon,y')dy' d\eta ' \nonumber \\
& +\frac{1}{(2 \pi h)^{d-\frac{1}{2}}} \iint_{T^*\R^{d-1}} \int_{-1}^{+ \infty} e^{i( \varphi (t,x,\eta ' )- y' \eta ' )/h}\chi (t) A_{+} (t,x, \eta ' ,z ,h) u_0(\varepsilon,y')dt dy' d\eta' , \label{p12}
\end{align}
is a solution to (\ref{zeproblem}) for any $z\in [-{{C_0}}h , {{C_0}}h] + i [- {{C_1}} h , {{C_1}}h]$.
\end{enumerate}
\end{theorem}

Precise definitions  for $A_+$ and $A_-$ are given in Section \ref{subsec:symbol} below. Notice that different choices for the cut-off function $\chi$ in (\ref{p12}) would lead to  the same microlocal solution in $\Omega$.

\Subsection{The phase function}
\label{subsec:phase}

We start with the construction of the phase function $\varphi$. From \eqref{symbolpschro}, for $x' = o (x_{1})$ and $\xi ' = o (x_{1})$, the equation $p_{0} (x, \xi_{1} , \xi ') = 0$ has two solutions
\begin{equation}
\xi_{1} = f_{\pm} (x, \xi ') = \pm \frac{\lambda_{1}}{2} x_{1} + o (x_{1}).
\end{equation}
Since $\gamma_{-}$ is a simple characteristic for the operator $P$, 
by usual Hamilton-Jacobi theory we have first the

\begin{lemma}\sl There exists a neighborhood $U_{-}$ of $x^{-}$, which depends on $\varepsilon$,  such that,
for all $\eta ' \in \R^{d-1}$ close enough to ${\xi^{-}}{}'$, there is a unique smooth function $\psi_{\eta '} : \R^d \to \R$, defined in $U_{-}$, verifying
\begin{equation}
\left\{ \begin{aligned}
& p_{0} (x,\nabla \psi_{\eta '}(x))=0, \\
& \psi_{\eta '}(x)=x' \cdot \eta' \text{ for } x \in  H_{-}\cap U_{-},\\
& \partial_{x_{1}} \psi_{\eta '}(x^-) = f_{-} (x^{-} , \eta ') .
\end{aligned}
\right.
\label{eikonale}
\end{equation}
\end{lemma}

If we denote by $\Lambda_{\psi_{\eta '}}$ the corresponding Lagrangian manifold
\begin{equation}
\Lambda_{\psi_{\eta '}} = \{ (x,\xi) \in T^*\R^d ; \ x \in U_{-},\ \xi=\nabla\psi_{\eta '}(x)\}, 
\label{Lambdapsi}
\end{equation}
 we have the following 
\begin{lemma}\sl
The Lagrangian manifolds $\Lambda_{-}$ and $\Lambda_{\psi_{\eta '}}$ intersect along an integral curve $\gamma_{\eta '}$ for $H_p$, and they intersect transversally. This curve is $\gamma_{-}$ when $\eta ' = \xi^{-}{}'$.
\label{intersectionpropre}
\end{lemma}

\begin{proof} First we study $\Lambda_{\psi_{\eta '}} \cap (H_{-} \times \R^d)$:  a point $(x_{1},x', \xi_{1}, \xi')$ belongs to this intersection if and only if $x_{1}= x^{-}_{1} = \varepsilon$ and $(\xi_{1},\xi')=\nabla_{x} \psi_{\eta '}(x_{1}^-,x')$. But  we have
 \begin{equation}
\nabla_{x'} \psi_{\eta '}(x_{1}^-,x')=\eta ',
\label{nabla'psi}
\end{equation}
and, moreover, $\psi_{\eta '}$ satisfies the eikonal equation.
Thus, using also the third equation of \eqref{eikonale}, we get by continuity
\begin{equation}
\partial_{x_{1}}\psi(x_{1}^-,x') = f_{-} (x_{1}^{-} , x' , \eta ' ),
\label{d1psi}
\end{equation}
and 
\begin{equation}
\Lambda_{\psi_{\eta '}} \cap (H_{-} \times \R^{d}) = \big\{ (x_{1}^-,x', f_{-} (x_{1}^{-} , x' , \eta ' ) ,\eta' ),\ x'\in\R^{d-1} \big\}.
\label{interH0}
\end{equation}
Then the intersection of $\Lambda_{\psi_{\eta '}} \cap (H_{-}\times \R^{d})$ with $\Lambda_{-}$ is given by
the equation
\begin{equation}
(x_1^-,x', f_{-} (x_{1}^{-} , x' , \eta ' ) ,\eta ' )=
(x_1^-,x',\nabla_{x} \varphi_{-}(x_1^-,x')),
\label{intersection}
\end{equation}
where $\varphi_{-}$ is a generating function for $\Lambda_{-}$ in $\Omega$ as in (\ref{phi+-en0}).

Let $g: \R^{d-1} \to \R^{d-1}$ be the function defined by $g(x') = \nabla_{x'} \varphi_{-} (x_{1}^{-} , x')$. In view of (\ref{phi+-}), we have $g(x^{-} {}') = \xi^{-}{}'$ and $\nabla_{x'} g(x^{-}{}') = \nabla_{x', x'}^2 \varphi_{-}(x^{-} {}' ) = -L'/2 +o(1)$ as $\varepsilon \to 0$.  Here $L'$ is the $(d-1)\times (d-1)$ matrix given by $L'=\diag(\lambda_2,\dots,\lambda_d)$ (see (\ref{linearise})). Thus the inverse function theorem implies that $g(x') = \eta '$ has a unique solution $x' = x' (\eta ')$ in a neighborhood of $x^{-} {}'$, for $\eta '$ in a neighborhood of $\xi^{-} {}'$. Notice also that,
\begin{equation} \label{x'(eta)}
x'(\eta ')= x^{-} {} ' + \CO (\vert \eta ' - \xi^- {} ' \vert) ,
\end{equation}
uniformly as $\varepsilon\to 0$ in a neighborhood of $\xi^- {} '$ which depends on $\varepsilon$. Since $\Lambda_{-}\subset p_{0}^{-1}(0)$, we  have 
\begin{equation}
\partial_{x_{1}} \varphi_{-}(x_1^-,x'(\eta '))= f_{-} (x_{1}^{-} , x'(\eta ') , \eta ' ) ,
\label{automatic}
\end{equation}
so that finally the equation (\ref{intersection}) has a unique solution $x'(\eta ')$ in a neighborhood of $x^{-} {}'$ for $\eta '$ close enough to $\xi^{-} {}'$.
 
Let us denote by 
\begin{equation}
\rho_{\eta '}=(x(\eta'),\xi(\eta'))=\big(x_1^- , x'(\eta '), f_{-} (x_{1}^{-} , x'(\eta ') , \eta ' ) ,\eta '\big)
\label{rhoeta} 
\end{equation}
the corresponding point. We show now that 
the tangent spaces at $\rho_{\eta '}$ to $\Lambda_{\psi_{\eta '}}$ and $\Lambda_{-}$ intersect along a one-dimensional space.

First it is clear that $H_{p}$ belongs to both $T_{\rho_{\eta '}} \Lambda_{\psi_{\eta '}}$ and $T_{\rho_{\eta '}}\Lambda_{-}$, since $\Lambda_{-}$ as well as $\Lambda_{\psi_{\eta '}}$ are invariant under the $H_{p}$ flow, or otherwise stated, because these Lagrangian manifolds are generated by solutions of the eikonal equation for $p$.

On the other hand, a vector $(\delta_{x},\delta_{\xi})$ belongs to $T_{\rho_{\eta '}}\Lambda_{\psi_{\eta '}} \cap T_{\rho_{\eta '}} \Lambda_{-}$ if and only if 
\begin{equation}
\left\{
\begin{array}{l}
  \delta_{\xi}=(\nabla^2_{x,x} \psi_{\eta} ) (x(\eta')) \delta_{x}, \\[6pt]
  \delta_{\xi}=(\nabla^2_{x,x} \varphi_{-} ) (x(\eta')) \delta_{x},
\end{array}
\right .
\end{equation}
or $\delta_{x}\in \ker \big( (\nabla^2_{x,x} \psi_{\eta '}) (x(\eta'))- (\nabla^2_{x,x} \varphi_{-} ) (x(\eta ')) \big)$. But we have seen that $\nabla^2_{x,x} \varphi_{-}( \rho_{-})=-L/2 +o(1)$ as $\varepsilon\to 0$, and that $(\nabla_{x',x'}^2 \psi_{\eta '}) (\rho_x(\eta'))=0$, so that the matrix $(\nabla^2_{x,x} \psi_{\eta '}) (x(\eta')) - (\nabla^2_{x,x} \varphi_{-}) (x(\eta'))$ has a $(d-1)\times (d-1)$ non-vanishing minor. Thus its rank is larger than $d-1$, and finally $H_{p}$ generates $T_{\rho_{\eta '}}\Lambda_{\psi_{\eta '}} \cap T_{\rho_{\eta '}}\Lambda_{-}$.
\end{proof}

Let $\gamma_{\eta '}$ be the hamiltonian curve with initial data $\rho (\eta ')$. We denote by $\Gamma^{\eta '}_0$ the set of level $\psi_{\eta '}(x (\eta '))$ for $\psi_{\eta '}$:
\begin{equation}
\Gamma^{\eta '}_0=\{(x,\xi) \in \Lambda_{\psi_{\eta '}} ; \ \psi_{\eta '}(x) = \psi_{\eta '}(x(\eta '))\},
\label{gammaOeta}
\end{equation}
and, possibly after shrinking $U_-$,  we have the

\begin{lemma}\sl
For $\varepsilon$ small enough, there exists a neighborhood $V_-$ of $\xi^- {}'$ such that,  for any $\eta '\in V_-$,  one can find a Lagrangian manifold $\Lambda_{0}^{\eta '}$ defined above $U_-$ such that
\begin{equation}
\Lambda_{0}^{\eta '}\cap \Lambda_{\psi_{\eta '}}= \Gamma_{0}^{\eta '},
\end{equation}
where the intersection is clean. Moreover $\Lambda_{0}^{\eta '}$ depends smoothly on $\eta '$, and $\Pi_x:\Lambda_{0}^{\eta '}\to U_-$ is a diffeomorphism.
\label{Lambda0}
\end{lemma}

\begin{proof}
It is sufficient  to prove the Lemma for $\eta ' = \xi^{-} {}'$. Indeed, every object that appears below evaluated at $\xi^-{}'$ is a smooth functions of $\eta'\in V_-$.
In particular  the estimates  below hold uniformly with respect to $\eta'$.

The vector $(\delta_{x},\delta_{\xi})$ belongs to $T_{\rho_-}\Gamma_{0}^{\xi^-{} '}$ if and only if 
\begin{equation}
\left\{
\begin{aligned}
&\delta_{\xi}=( \nabla^2_{x,x} \psi_{\xi^-{} '}) (x^-) \delta_{x}, \\
&( \nabla_{x} \psi_{\xi^-{} '} ) (x^-) \delta_{x}=0.
\end{aligned}
\right .
\label{tangentGamma0}
\end{equation}
Indeed $\Gamma^{\xi^-{} '}_0\subset \Lambda_{\psi}^{\xi^-{} '}$, and
$\Gamma_{0}^{\xi^-{} '}$ is a level curve for $\psi$. Thanks to
(\ref{nabla'psi}) and (\ref{d1psi}), the second equation becomes
\begin{equation}
\delta_{x}^1= - \frac{\xi^-{} ' \cdot \delta_{x}'}{f_{-} (x^{-}, \xi^-{} ' )}
\label{orthoniveau}
\end{equation}
and we see in particular that $T_{\rho_-}\Gamma_{0}^{\xi^-{} '}$ is parametrized by $\delta_{x}'$.

Let us compute the entries of the matrix $M_{\varepsilon}= (
\nabla^2_{x,x} \psi_{\xi^-{} '} ) (x^-)$.
We have already seen (see (\ref{nabla'psi})) that, for $i,j\geq 2$, $m_{ij}=0$.
We also know  (see (\ref{d1psi})) that
\begin{equation}
( \nabla_{x '} \partial_{x_{1}} \psi ) (x^-) = \nabla_{x'} f_{-} (x^{-} , \xi^{-}{}' ) =
\frac{L'{}^{2} x^{-}{}' + \CO ( \vert x_{1}^{-} \vert^{2})}{4 f(x^{-} , \xi^{-}{}' ) },
\label{dx'd1psi}
\end{equation}
and we are left with the computation of $\partial^2_{x_{1} , x_{1}}
\psi(x^-)$. But we have seen that $H_{p}(\rho_-)=(\nabla_{\xi} p_{0} (\rho_{-}), - \nabla_{x} p_{0} (\rho_{-}) )$ belongs to
$T_{\rho_-}\Lambda^{\xi^-{} '}_{\psi}$, that is satisfies the first
equation in (\ref{tangentGamma0}), so that
\begin{equation}
\frac{1}{2} L^{2} x^{-} + \CO ( \vert x_{1}^{-} \vert^{2}) = M_{\varepsilon} \big( 2 \xi^{-} + \CO ( \vert x_{1}^{-} \vert^{2}) \big)
\end{equation}
which gives in particular
\begin{equation}
\frac{\lambda_{1}^{2}}{2} x^{-}_{1} + \CO ( \vert x_{1}^{-} \vert^{2}) = m_{11} ( 2 \xi_{1}^{-} + \CO ( \vert x_{1}^{-} \vert^{2})) + \frac{L'{}^{2} x^{-} {}' \cdot \xi^{-} {}' + \CO ( \vert x_{1}^{-} \vert^{3})}{2 \xi_{1}^{-} } ,
\end{equation}
so that 
\begin{equation}
m_{11} = - \frac{\lambda_{1}}{2} + o (1).
\label{m11}
\end{equation}
as $\varepsilon\to 0$. Here we recall that, as $x_{1}^{-} = \varepsilon$ goes to $0$, we have $\xi_{1}^{-} = -\frac{\lambda_{1}}{2} \varepsilon + o (\varepsilon)$, $x^{-} {}' = o (\varepsilon)$ and $\xi^{-} {}' = o (\varepsilon)$.

Thus 
\begin{equation}
M_{\varepsilon} = \left( \begin{array}{ccc}
- \frac{\lambda_{1}}{2} & 0 & \cdots \\
0 & 0 \\
\vdots & & \ddots \\
\end{array} \right) + o(1),
\end{equation}
and the equation (\ref{orthoniveau}) becomes
\begin{equation}
\delta_{x}^1 = o(\delta_{x}').
\end{equation}

Summing up, we see, using (\ref{tangentGamma0}), that  the vectors of $T_{\rho_-}\Gamma_{0}^{\xi^-{} '}$ can be written, when $\varepsilon\to 0$, as
\begin{equation}
\big( (0,\delta_{x}'),(0,0) \big) + o(\delta_{x}'),\  \delta_{x}'\in \R^{d-1}.
\label{tangentgamma0delta}
\end{equation}
Let us denote by $\CE_{0}$ the ``limit space" for $T_{\rho_-}\Gamma_{0}^{\xi^-{} '}$, that is the linear subspace of $\R^{2d}$ generated by the $e_{j}$'s for $j=2, \dots, d$, where $e_{j}=(\delta_{i,j})_{i=1, \dots , 2d}$. It is clear that $e_{1} \R \oplus \CE_{0}$ is a Lagrangian subspace of $\R^{2d}$. Then, using the Gram--Schmidt orthonormalization principle, one can find a unitary vector $v (\varepsilon) \in \R^{2d}$ such that
\begin{equation}
\sigma \left( u, v (\varepsilon) \right) = 0, \ \text{for all } u \in T_{\rho_-}\Gamma_{0}^{\xi^-{} '},
\end{equation}
and
\begin{equation}  \label{vepsilon}
v (\varepsilon) = e_{1} + o (1),
\end{equation}
as $\varepsilon \to 0$. Then, $v (\varepsilon) \R \oplus T_{\rho_-} \Gamma_{0}^{\xi^-{} '}$ is a Lagrangian vector space at $\rho_-$, and, extending $\Gamma_{0}^{\xi^-{} '}$ along a suitably chosen Hamilton field,   one can find locally close to $\rho_-$, a Lagrangian manifold $\Lambda_{0}^{\xi^-{} '}$ such that $\Gamma_{0}^{\xi^-{} '} \subset \Lambda_{0}^{\xi^-{} '}$ and
\begin{equation}
T_{\rho_-} \Lambda_{0}^{\xi^-{} '} = v (\varepsilon) \R \oplus T_{\rho_-} \Gamma_{0}^{\xi^-{} '}.
\end{equation}
 Moreover, if $(\delta_{x}, \delta_{\xi}) \in T_{\rho_-} \Lambda_{0}^{\xi^-{} '}$, we get, from (\ref{tangentgamma0delta}) and (\ref{vepsilon}),
\begin{equation}  \label{ffd}
\delta_{\xi} = o ( \delta_{x} ).
\end{equation}

To show that the intersection $\Lambda_{0}^{\xi^-{} '}\cap \Lambda_{\psi}^{\xi^-{} '}$ is clean, is is enough to show that $H_{p}(\rho_{-} ) \in T_{\rho_-} \Lambda_{\psi}^{\xi^-{} '}$ is not in $T_{\rho_-} \Lambda_{0}^{\xi^-{} '}$. As $\varepsilon\to 0$, we have
\begin{equation}
H_{p}(\rho_-)=\left(
\begin{array}{c}
2 \xi^-_1  \\
2 \xi^-{}'   \\
\lambda_{1} x_{1}^{-} /2   \\
L' x^{-} {}' /2
\end{array}
\right ) + \CO (\varepsilon^{2})
=
\left(
\begin{array}{c}
-\lambda_{1}\varepsilon /2 \\
0\\
\lambda_{1} \varepsilon /2 \\
0
\end{array}
\right )
+o(\varepsilon).
%\label{}
\end{equation}
Then (\ref{ffd}) implies that $H_{p}(\rho_-) \notin T_{\rho_-} \Lambda_{0}^{\xi^-{} '}$ and the dimension of the intersection
$T_{\rho_-}\Lambda_{0}^{\xi^-{} '} \cap T_{\rho_-}\Lambda_{\psi}^{\xi^-{} '}$
is exactly $d-1$.

Finally, the Lagrangian manifold $\Lambda^{\xi^-{} '}_0$ projects nicely on the $x$-space: Indeed if $(\delta_{x},\delta_{\xi})\in T_{\rho_-} \Lambda_{0}^{\xi^-{} '}$, we know by (\ref{ffd}), that as $\varepsilon\to 0$, $\delta_{\xi} = o (\delta_{x})$, so that $\delta_{x}\neq 0$ for any $\varepsilon$  small enough.
\end{proof}

\begin{figure}
\begin{center}
\begin{picture}(0,0)%
\includegraphics{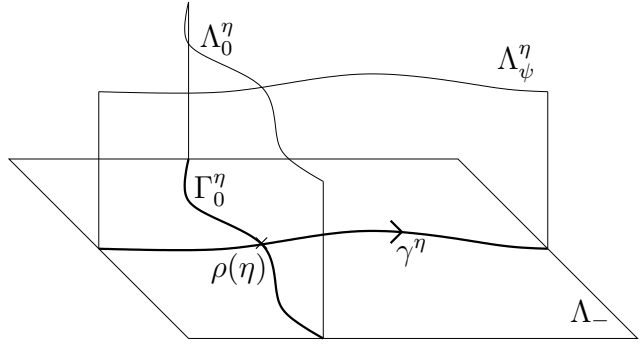}%
\end{picture}%
\setlength{\unitlength}{1855sp}%
\begingroup\makeatletter\ifx\SetFigFont\undefined%
\gdef\SetFigFont#1#2#3#4#5{%
  \reset@font\fontsize{#1}{#2pt}%
  \fontfamily{#3}\fontseries{#4}\fontshape{#5}%
  \selectfont}%
\fi\endgroup%
\begin{picture}(8424,4545)(2389,-6394)
\put(4951,-2461){\makebox(0,0)[lb]{\smash{{\SetFigFont{12}{14.4}{\rmdefault}{\mddefault}{\updefault}$\Lambda_{0}^{\eta}$}}}}
\put(9901,-6136){\makebox(0,0)[lb]{\smash{{\SetFigFont{12}{14.4}{\rmdefault}{\mddefault}{\updefault}$\Lambda_{-}$}}}}
\put(8926,-2761){\makebox(0,0)[lb]{\smash{{\SetFigFont{12}{14.4}{\rmdefault}{\mddefault}{\updefault}$\Lambda_{\psi}^{\eta}$}}}}
\put(4876,-4486){\makebox(0,0)[lb]{\smash{{\SetFigFont{12}{14.4}{\rmdefault}{\mddefault}{\updefault}$\Gamma_{0}^ {\eta}$}}}}
\put(7576,-5311){\makebox(0,0)[lb]{\smash{{\SetFigFont{12}{14.4}{\rmdefault}{\mddefault}{\updefault}$\gamma^{\eta}$}}}}
\put(5101,-5536){\makebox(0,0)[lb]{\smash{{\SetFigFont{12}{14.4}{\rmdefault}{\mddefault}{\updefault}$\rho(\eta)$}}}}
\end{picture}%
\end{center}
\caption{The Lagrangian manifolds.}
\label{figlag}
\end{figure}

Now we consider the associated Lagrangian manifold
\begin{equation}   \label{lambdat}
\Lambda_{t}^{\eta '} = \exp(tH_{p}) (\Lambda_{0}^{\eta '}).
\end{equation}
The manifold $\Lambda_{t}^{\eta '}$ projects nicely on $\R^{d}_{x}$. In fact, possibly after shrinking $V_-$, we have the

\begin{lemma}\sl
There exists $T_{0} >0$ such that for any $\varepsilon >0$ small enough, there exist $\delta >0$ and $V_-$ a neighborhood of $\xi^- {}'$ such that for all $\eta ' \in V_-$, the manifold $\Lambda_{t}^{\eta '}$ projects nicely onto $U_{t} = B (x_{-} (t), \delta)$ for $t \in [-1 , T_{0} ]$ and onto $U_{\infty}$ for $t> T_{0}$. Here $U_{\infty}$ is a neighborhood of $0\in \R^d$, such that $B (x_{-} (t), \delta) \subset U_{\infty}$, for $t> T_{0} $.
\label{bonneprojection}
\end{lemma}

\begin{proof}
Let $(\delta_{x} , \delta_{\xi})$ be in the tangent space
$T_{\rho_{\eta '} (t)} \Lambda_{t}^{\eta '}$. The proof of Lemma 2.1
of \cite{hsmw3} implies that
\begin{equation}
\delta_{\xi} - L/2 \delta_{x} = B_{t} \left( \delta_{\xi} + L/2
\delta_{x} \right),
\end{equation}
with $B_{t} = {\mathcal O} \left( e^{-\lambda_{1} t} \right)$,
uniformly with respect to $\varepsilon$ and $\eta '$. Then, for each
$\widetilde{\varepsilon}>0$, there is a $T_{0} >0$, such that
\begin{equation}   \label{hjk1}
|\delta_{\xi} - L/2 \delta_{x}| \leq \widetilde{\varepsilon} |\delta_{x}|,
\end{equation}
for $(\delta_{x} , \delta_{\xi} ) \in T_{\rho_{\eta '} (t)}
\Lambda_{t}^{\eta '}$, $t> T_{0}$, uniformly with respect to
$\varepsilon$ and $\eta '$. This inequality, together  with \cite[Lemma 2.2]{hsmw3}, gives the Proposition for $t > T_{0} $.

For $t \in [-1, T_{0} ]$, it is enough to prove the Lemma for $\eta ' = \xi^{-} {}'$, as in the proof of Lemma \ref{Lambda0}.
We shall use the fact that, on $[-1,T_{0} ]$, the evolution of a tangent vector is
closed to the evolution for the reference operator $p_{0} = \xi^{2} - \sum \lambda_{j}^{2} x_{j}^{2} /4$, provided $\varepsilon$ is small enough. If
$(\delta_{x} (t) , \delta_{\xi}(t)) \in T_{\rho_{-} (t)}
\Lambda_{t}^{\xi^{-} {}'}$ is the evolution of a tangent vector $(\delta_{x}, \delta_{\xi})$
along the integral curve $\gamma_{-}$, we have
\begin{align}
\delta_{x}^j (t) =&\frac12 (e^{\lambda_{j}t}+e^{-\lambda_{j}t})
\delta_{x}^j + \frac1{\lambda_{j}}(e^{\lambda_{j}t}-e^{-\lambda_{j}t})
\delta_{\xi}^j + o(\delta_{x})   \nonumber  \\
=& \frac12 (e^{\lambda_{j}t}+e^{-\lambda_{j}t}) \delta_{x}^j +
o(\delta_{x}),   \label{hjk2}  \\
\delta_{\xi}^j (t) =& \frac{\lambda_{j}}4
(e^{\lambda_{j}t}-e^{-\lambda_{j}t}) \delta^j_{x} +
\frac1{2}(e^{\lambda_{j}t}+ e^{-\lambda_{j}t})\delta^j_{\xi} +
o(\delta_{x})   \nonumber  \\ 
=& \frac{\lambda_{j}}4
(e^{\lambda_{j}t}-e^{-\lambda_{j}t})\delta^j_{x} + o(\delta_{x}),
\label{hjk3} 
\end{align}
since $\delta_{\xi} = o(\delta_{x})$ by (\ref{ffd}). From
(\ref{hjk2}) and (\ref{hjk3}), we see that $\delta_{\xi} (t)$ is a
function of $\delta_{x} (t)$, and that proves the Lemma.
\end{proof}

We set
\begin{equation}
\widetilde{U}_{t} = \left\{ \begin{aligned}
&U_{t} \quad &&\text{for } t \in [-1,T_{0} ],  \\
&U_{\infty} \quad &&\text{for } t \in ] T_{0}, + \infty[.
\end{aligned} \right.
\end{equation}
Thanks to Lemma \ref{bonneprojection}, there is a smooth function $\varphi (t,x,\eta ')$ defined on $]-1,+\infty[\times \widetilde{U}_{t} \times V_-$ such that the Lagrangian manifold $\Lambda_t^{\eta '}$ is given by
\begin{equation}
\xi=\nabla_x \varphi (t,x,\eta ') \quad \text{for } x\in \widetilde{U}_{t}.
\label{lambdatdef}
\end{equation}
It satisfies of course the eikonal equation
\begin{equation}  \label{eikonalet}
\partial_t\varphi(t,x,\eta ')+p_{0} (x,\nabla_x\varphi(t,x,\eta '))=0.
\end{equation}
Therefore, it follows from \cite[Theorem 3.12]{hsmw3}, that $\varphi(t,x,\eta')$ is expandible in the sense of Definition \ref{expdeb}: There exists a sequence $\varphi_{j}$ of smooth functions on $]-1,+\infty[\times \widetilde{U}_{t} \times V_-$ which are polynomials in $t$, such that for any $N, k\in\N$, $\alpha, \beta\in \N^d$,
\begin{equation}
\partial^k_{t}\partial_{x}^\alpha \partial_{\eta'}^\beta \Big( \varphi(t,x,\eta')-\sum_{j= 0}^N 
\varphi_{j}(t,x,\eta')e^{-\mu_{j}t} \Big)=\CO(e^{-\mu_{N}t}) .
\label{defexpandible}
\end{equation}

Now  we set
\begin{equation}
\Gamma_{t}^{\eta '} = \exp (t H_{p}) \Gamma_{0}^{\eta '},
\end{equation}
and we have, possibly after shrinking $\widetilde{U}_{t}$ and $V_-$, the

\begin{proposition}\sl
For each $\eta ' \in V_{-}$ and $x \in \bigcup_{t} \widetilde{U}_{t} \cap \{x; \  0< \vert x' \vert < \nu x_{1} \}$, for some $\nu >0$, there is a unique time $t=t(x,\eta ')$ such that
$x \in \Pi_{x} \Gamma_{t}^{\eta '}$. Moreover, it is the only critical
point for the function $t\mapsto \varphi(t,x,\eta ')$, and it is a non-degenerate critical point.
\label{critique}
\end{proposition}

\begin{proof}
If $x \in \Pi_{x} \Gamma_{t}^{\eta '}$, there is a $\xi \in \R^{n}$
such that $(x,\xi ) \in \Gamma_{t}^{\eta '}$. Then
\begin{equation}   \label{yyy}
\xi = \nabla_x \varphi(t,x, \eta ' ),
\end{equation}
and $p_{0} (x,\xi) =0$ since $\Gamma_{0}^{\eta '} \subset p_{0}^{-1} (0)$ and the
Hamiltonian flow preserves the energy. Together with (\ref{eikonalet}),
we get that $t$ is a critical point for the function $t\mapsto
\varphi(t,x,\eta ')$ if and only if $x \in \Pi_{x} \Gamma_{t}^{\eta '}$.

In the case $x \in U_{\infty}$, the proposition follows from \cite[Lemma 3.14]{hsmw3}. For $x \notin U_{\infty}$, it is enough to see that $\partial_{t}^{2} \varphi (t,x,\xi^{-} {}') >0$ for $x = x^{-} (t)$. The eikonal equation (\ref{eikonalet}) implies
\begin{align}
\nabla_{x} \partial_{t} \varphi =& - 2 {\rm Hess} ( \varphi ) \nabla_{x} \varphi + L^{2}x/2 + \CO \big( (x^{2} + \vert \nabla_{x} \varphi \vert^{2} ) ( \vert {\rm Hess} ( \varphi ) \nabla_{x} \varphi \vert +1) \big)  \label{triste1}  \\
\partial_{t}^{2} \varphi =& - 2 \nabla_{x} \partial_{t} \varphi \cdot \nabla_{x} \varphi + \CO \big( \vert \nabla_{x} \partial_{t} \varphi \vert (x^{2} + \vert \nabla_{x} \varphi \vert^{2} ) \big).   \label{triste2}
\end{align}
From (\ref{ffd}), (\ref{hjk1}), (\ref{hjk2}) and (\ref{hjk3}), we get that ${\rm Hess} ( \varphi )=\CO(1)$, and that
\begin{equation}
{\rm Hess} ( \varphi )>o(1) ,
\end{equation}
as $\varepsilon \to 0$ because $\delta_{\xi} = {\rm Hess} ( \varphi ) \delta_{x}$ for
$(\delta_{x},\delta_{\xi}) \in T_{\rho_{-} (t)} \Lambda^{\xi^{-} {}'}_t$. Since we assume
that $\Pi_{x} \gamma_{1}$ is collinear to $x_{1}$ (see the remark
after (\ref{uio})), we also have
\begin{align*}
x^{-} (t) &= (x^{-}_{1} (t),0,\ldots,0) + o(x^{-}_{1} (t)) \\
\xi^{-} (t) &= (-\lambda_{1} x^{-}_{1} (t) /2 ,0,\ldots,0) +
o(x^{-}_{1} (t)).
\end{align*}
and then \eqref{triste1} and \eqref{triste2} become
\begin{align}
\nabla_{x} \partial_{t} \varphi (t , x^{-} (t) , \xi^{-} {}') =& - 2 {\rm Hess} ( \varphi ) \xi^{-} (t) + L^{2} x^{-} (t) /2 + \CO \big( x^{-} (t)^{2} \big)   \\
\partial_{t}^{2} \varphi (t , x^{-} (t) , \xi^{-} {}') =& - 2 \nabla_{x} \partial_{t} \varphi \cdot \xi^{-}(t) + \CO \big( x^{-} (t)^{3} \big)    \nonumber  \\
=& 4 {}^{t}\xi^{-} (t) {\rm Hess} ( \varphi ) \xi^{-} (t) - L^{2} x^{-} (t) \cdot \xi^{-} (t) + \CO \big( x^{-} (t)^{3} \big)     \nonumber \\
\geq & - L^{2} x^{-} (t) \cdot \xi^{-} (t) + o\big( x^{-} (t)^{2} \big)  \nonumber \\
\geq & \lambda_{1}^{3} (x^{-} (t))^{2} + o\big( x^{-} (t)^{2} \big)  \nonumber \\
>& 0.
\end{align}
\end{proof}

As a consequence of Proposition \ref{critique}, we get in particular
that, in 
\begin{equation}
\widehat{U} = \bigcup_{t} \widetilde{U}_{t} \cap \{x; \  0< \vert x' \vert < \nu x_{1} \}
\end{equation}
where both these functions are defined, we have
\begin{equation}
\nabla_x \psi_{\eta '} (x) = \nabla_x \left( \varphi (t(x,\eta '),x, \eta ' ) \right) .
\label{psietphi}
\end{equation}
Therefore $x\mapsto \psi_{\eta '} (x)$ and $x\mapsto \varphi (t(x,\eta '),x,\eta ')$ differ from a constant.
Then, adding a constant (with respect to  $t$, $x$) to $\varphi (t ,x, \eta ' )$, we can assume that
\begin{equation}
\varphi ( t( x,\eta '), x,\eta ')=x' \cdot \eta ' ,
\label{phiinitiale}
\end{equation}
for any $x\in H_{-}\cap \widehat{U}$. 
Furthermore, we can compute the first term in the expansion (\ref{defexpandible}):

\begin{lemma}\sl
In the sense of expandible functions, we have
\begin{equation}  \label{pi1}
\varphi(t,x,\eta) \sim \varphi_+(x) + \widetilde{\psi} (\eta ') + \sum_{j\geq
  1} e^{-\mu_j t} \varphi_j (t,x,\eta '),
\end{equation}
where the $\varphi_{j} (t,x, \eta ')$ are polynomials in $t$ with
smooth coefficients in $x, \eta '$, and 
\begin{equation} \label{pi2}
\widetilde{\psi} (\eta ') = x' (\eta ') \cdot \eta ' - \varphi_- (x
(\eta ')).
\end{equation}
\end{lemma}

\begin{proof}
As we have already mentioned, the asymptotic (\ref{pi1}) follows from the proofs of sections 2 and 3 of
\cite{hsmw3}, and we are left with the proof of  (\ref{pi2}). Let us denote by
$(x (t),\xi (t))$ the points on the curve $\gamma_{\eta '}$
defined in Lemma \ref{intersectionpropre}, with
$(x (0),\xi (0)) = \rho_{\eta '} = (x (\eta ') , \xi (\eta ') ) \in
H_- \times \R^d$ given by (\ref{rhoeta}). We notice that, by
(\ref{pi1}),
\begin{equation}
\widetilde{\psi}(\eta ') = \lim_{t \to +\infty} \varphi (t,x (t),\eta ').
\end{equation}
On the other hand, by the eikonal equation (\ref{eikonalet}) and since $(x(t),\xi(t))\in\gamma_{\eta'}\subset p_{0}^{-1}(0)$, we have
\begin{align}
\partial_{t} ( \varphi (t,x (t),\eta ') ) =& (\partial_t \varphi) (t,x
(t) , \eta ') + (\partial_x \varphi )(t,x(t),\eta ' ) \cdot
 (\partial_{t} x ) (t)    \nonumber  \\ 
=& \xi (t) \cdot (\partial_{t} x) (t) = \nabla \varphi_- (x(t)) \cdot
( \partial_{t} x) (t)    \nonumber  \\
=& \partial_{t} ( \varphi_{-} (x (t)) ) 
\end{align}
where we use also the fact that $\gamma_{\eta '}
\subset \Lambda_{-}$. Therefore, we get, with (\ref{phiinitiale}),
\begin{align}
\widetilde{\psi}(\eta ') =& \lim_{t \to +\infty} \varphi (t,x (t),\eta
')- \varphi_{-} (x(t))    \nonumber   \\
=& \varphi (0,x (0),\eta ')- \varphi_{-} (x(0)) = x' (\eta ') \cdot
\eta '- \varphi_{-} (x (\eta ')),
\end{align}
which is (\ref{pi2}).
\end{proof}

% Before the study of the symbol, we prove Remark \ref{p1}.
% 
% \begin{proof}[Proof of Remark \ref{p1}]
% Since $m$ of the $\lambda_{j}$ are equal to $\lambda_{1}$, one can find $m$ independent vectors
%$v_{1}, \ldots, v_{m}$ in $\ker ( F_{p} - \lambda_{1})$. Using Hartmann's Theorem, one can find
%$m$ Hamiltonian curves $\gamma_{1} (t), \ldots , \gamma_{m} (t)$ with $g_{1}^{+} (\gamma_{j} (1) )
%= v_{j}$.
%
%
%\begin{equation}
%\left\{ \begin{aligned}
%& (  \nabla_\xi p_0(x,\nabla \varphi_{+})\cdot \nabla - \lambda_{1} ) \varphi_{1}
%  = 0 , \\
%&\nabla \varphi_{1} (0) = -\lambda_{1} g_1^-(\rho_-) .
%\end{aligned} \right.
%\end{equation}
%We set  $\widetilde{\Lambda_{+}}(\rho_{-}) = \{ (x, \xi) \in
%\Lambda_{+} ; \ \varphi_{1}(x) = 0 \}$.
%Then, $\widetilde{\Lambda_{+}}(\rho_{-})$ is a $\CC^{\infty}$
%submanifold of $\Lambda_{+}$, of codimension $1$, which is stable under
%the Hamiltonian flow,
%\end{proof}

\Subsection{The symbol}
\label{subsec:symbol}

Now we look for a symbol $a(t,x,\eta',z,h)=\sum_k a_k(t,x,\eta',z)h^k$ such that (\ref{mastereq}) holds. This leads to the usual transport equations for the $a_j$'s (see \cite[Theorem IV-19]{Ro}):
\begin{equation}\label{transp}
\left\{
\begin{aligned}
&\partial_t a_0 + \partial_{\xi} p_{0} (x, \partial_{x} \varphi) \partial_{x} a_{0} + \Big( \frac{1}{2} \tr \big( \partial^{2}_{\xi , \xi} p_{0} (x, \partial_{x} \varphi) \partial^{2}_{x,x} \varphi \big) -i\frac{z}h \Big) a_{0} =0, \\
&\partial_t a_k + \partial_{\xi} p_{0} (x, \partial_{x} \varphi) \partial_{x} a_{k} + \Big( \frac{1}{2} \tr \big( \partial^{2}_{\xi , \xi} p_{0} (x, \partial_{x} \varphi) \partial^{2}_{x,x} \varphi \big) -i\frac{z}h \Big) a_{k} = F_{k} , \qquad k \geq 1,
\end{aligned}
\right.
\end{equation}
where $F_{k} (a_{0}, \ldots , a_{k-1})$ is a differential operator on the $a_{0}, \ldots , a_{k-1}$ with smooth coefficients. In the Schr\"{o}dinger case ($p= \xi^{2} + V(x)$), these equations become
the more familiar
\begin{equation}\label{transp2}
\left\{
\begin{aligned}
&\partial_t a_0+2\nabla_x\varphi\cdot\nabla_x a_0 +(\Delta_x\varphi\  -i\frac{z}h)a_0=0,\\
&\partial_t a_k +2\nabla_x\varphi\cdot\nabla_x a_k +(\Delta_x\varphi\
-i\frac{z}h)a_k = i\Delta_x a_{k-1}, \qquad k \geq 1.
\end{aligned}
\right .
\end{equation}
Let us denote by $x_{\eta'} (t)$ the spacial projection of the curve $\gamma_{\eta '}$ defined in Lemma \ref{intersectionpropre}. As in \cite{hsmw3}, using the time-dependent change of coordinates $y=x-x_{\eta '} (t)$, the transport equations \eqref{transp} can be written as
\begin{eqnarray}\label{transp3}
\nonumber
&&
\partial_t a_k + \Big( \partial_{\xi} p_{0} \big( x_{\eta '} (t)+  y, \partial_{x} \varphi (t, x_{\eta '} (t)+y) \big)- \partial_{\xi} p_{0} \big(x_{\eta '} (t), \partial_{x} \varphi (t, x_{\eta '} (t)) \big) \Big) \partial_{y} a_{k} \\
&&
\qquad\qquad + \Big( \frac{1}{2} \tr \big( \partial^{2}_{\xi , \xi} p_{0} ( \cdot , \partial_{x} \varphi ( t, \cdot) ) \partial^{2}_{x,x} \varphi ( t, \cdot ) \big) -i\frac{z}h \Big) (x_{\eta '} (t)+y) a_{k} = F_{k}.
\end{eqnarray}

We also want that the function $u$ given by (\ref{solution}) satisfies
the initial condition $u=u_0$  microlocally in $U$.
Performing a formal stationary phase expansion with respect to $t$ in
(\ref{solution}), we get, for $x=(\varepsilon ,x') \in H_{-}$,
\begin{equation}  \label{erp1}
u (x,h) = \frac{1}{(2 \pi h)^{d-1}} \iint_{T^* \R^{d-1}} e^{i (x' \cdot \eta '- y' \cdot \eta ' )/h} \widetilde{a} (x' , \eta ',z,h) u_{0}(\varepsilon , y' ) dy' d\eta '  ,
\end{equation}
where $\widetilde{a} (x' , \eta ', z,h)$ is another classical symbol, whose principal part is given by
\begin{equation} \label{kjn1}
\widetilde{a}_0= e^{i\pi /4} \frac{a_0(t(x, \eta '), x, \eta
  ',z,h)}{\vert \partial_{tt}^{2} \varphi (t(x, \eta '), x, \eta ')
  \vert^{1/2}} \cdotp
\end{equation}
Since we want that $u (x,h)$ coincides with $u_{0} (x,h)$ on
$H_{-}$, we look for a symbol $a (t,x,\eta ', z,h)$ such that
\begin{equation}  \label{rtt}
\widetilde{a} (x' , \eta ', z,h) =1 + {\mathcal O} (h^{\infty}).
\end{equation}
From the structure of the stationary phase expansion,
there exists a unique formal classical symbol $a_{ini} (x', \eta ',z,h)$
which solves the problem (\ref{rtt}). And since the vector field
$(\partial_t, \nabla_x \varphi \cdot \nabla_x)$ is not tangent to
the hypersurface $\R \times H_-$ in $\R^{d+1}$, we can determine
uniquely solutions $a_j$ to the problem (\ref{transp}) which satisfy
\begin{equation}  \label{condinit}
a (t, x, \eta ',z,h) = a_{ini} (x',\eta ',z,h),
\end{equation}
for all $x = (\varepsilon, x') \in H_{-}$ and $t \in \R$. Notice that the $a_j$'s depend holomorphically on the parameter $z$.
 
Moreover, by \eqref{transp3} and Proposition \ref{p22}, the functions $a_k$ are expandible with respect to $(x, \eta')$ in the modified sense that the family of exponents is now $(S+\mu_j)_{j\in \N}$, where 
 \begin{equation}\label{defS}
S=S(z/h)=\Delta\varphi_+(0)-i\frac zh=\sum_{j=1}^d \frac{\lambda_j}2 -i\frac zh\cdotp
\end{equation}
We can also  find a realization $a$, holomorphic with respect to $z$, of the asymptotic sum $\sum a_k (t, x, \eta ',z ) h^k$ such that
\begin{equation}\label{realza}
a (t, x, \eta ',z,h) \in \CS^{0,\re S}
\end{equation}
and 
\begin{equation} \label{ram9}
r=e^{-i\varphi/h}(h D_{t} + P(x,hD) - z) a e^{i\varphi /h} \in \CS^{\infty,\re S}.
\end{equation}

Now we want to give a meaning to  the integral
\begin{equation}\label{}
\int_{-1}^{+\infty} e^{i\varphi(t,x,\eta')/h} \chi (t) a(t,x,\eta',z,h)dt ,
\end{equation}
where $\chi \in C^{\infty} (]-1, + \infty[)$ equal to $1$ near $[0, + \infty[$.
Notice that with respect to the situation in \cite[Section 4]{hsmw3}, here we have to deal with an oscillatory integral. As soon as $\re S>0$, this integral is absolutely convergent. But if  $\re S \leq 0$, there might exist $j$'s in $\N$ such that $\re S+\mu_j \leq 0$, and then the integral above has no obvious meaning. Nevertheless,  we explain now how to obtain a solution even in that case.

We set
\begin{equation}\label{K1gh}
K_{1} = {\mathbb E} \Big( \frac{{{C_1}}}{\lambda_{1}} - \frac{\sum \lambda_{j}}{2 \lambda_{1}} \Big) +1,
\end{equation}
and we also  denote $\varphi_{\infty}(x, \eta ') = \varphi_+(x) + \widetilde{\psi} (\eta ')$, $\varphi_{\star} (t, x, \eta ') = \varphi - \varphi_{\infty} = \CO (e^{-\lambda_{1}t})$. Then we can write
\begin{equation}
a e^{i\varphi /h} = a e^{i\varphi /h} - \sum_{k< K_{1}} \frac{a}{k !} \Big( \frac{i \varphi}{h} \Big)^{k} e^{i \varphi_{\infty} /h} + \sum_{k< K_{1}} \frac{a}{k !} \Big( \frac{i \varphi}{h} \Big)^{k} e^{i \varphi_{\infty} /h} .
\end{equation}
From our  choice for $K_{1}$, there exists $\delta >0$ such that for all $( \alpha , \beta , \gamma ) \in \N^{1+d+(d-1)}$ and $z \in [-{{C_0}}h , {{C_0}}h] + i [- {{C_1}} h , {{C_1}}h]$,
\begin{equation}  \label{ram1}
\partial_{t}^{\alpha} \partial_{x}^{\beta} \partial_{\eta '}^{\gamma} \Big( a e^{i \varphi /h} - \sum_{k< K_{1}} \frac{a}{k !} \Big( \frac{i \varphi_{\star}}{h} \Big)^{k} e^{i \varphi_{\infty} /h} \Big) \lesssim h^{-K_{1} - \vert \alpha \vert - \vert \beta \vert - \vert \gamma \vert} e^{- 3 \delta t},
\end{equation}
uniformly with respect to  $h$ and $t$.

On the other hand, 
\begin{equation}  \label{ram45}
b = \sum_{k< K_{1}} \frac{a}{k !} \Big( \frac{i \varphi_{\star}}{h} \Big)^{k} \sim \sum_{1-K_{1} \leq k} b_{k} (t, x, \eta ' ,z) h^{k}
\end{equation}
is expandible  for the family of exponents $(S + \mu_j)_{j}$:
\begin{equation}  \label{ram5}
b_{k} (t, x, \eta ' ,z) \sim \sum_{j} b_{k,j} (t, x, \eta ',z) e^{- (S + \mu_{j}) t},
\end{equation}
where $b_{k,j}$ is  polynomial with respect to $t$. Let $J_{1} \in \N$ be such that
\begin{equation}
\mu_{J_{1}} > 2 \delta - \sum_{j} \lambda_{j} /2 + {{C_1}}.
\label{voilaJ1}
\end{equation}
As in \cite{hsmw3}, for an expandible symbol satisfying \eqref{ram45} and \eqref{ram5}, we define
\begin{equation}
[ b_{k} ]_{-} = \sum_{j< J_{1}} b_{k,j} e^{- (S + \mu_{j}) t} \in
\CS^{0 , \re S},   \mbox{ and }
[ b_{k} ]_{+} = b_{k} - [ b_{k} ]_{-} \in \CS^{0 , 2 \delta}.     \label{ram7}
\end{equation}
Using Borel's lemma, we can find  $[ b ]_{+}$ and then $[ b ]_{-}$, holomorphic with respect to $z$ in $[-{{C_0}}h , {{C_0}}h] + i [- {{C_1}} h , {{C_1}}h]$, such that
\begin{equation}
[ b ]_{+} \sim \sum_{-K_{1} \leq k} [ b_{k} ]_{+} h^{k} \in
\CS^{1-K_{1} ,2 \delta}, \mbox{ and }
[ b ]_{-} = b - [ b ]_{+} \in \CS^{1-K_{1} , \re S} .  \label{ram3}
\end{equation}
Then the function
\begin{equation}  \label{ram4}
A_{+} (t,x, \eta ',z,h) = a e^{i\varphi /h} - \sum_{k< K_{1}}
\frac{a}{k !} \Big( \frac{i \varphi_{\star}}{h} \Big)^{k} e^{i
\varphi_{\infty} /h} + [ b ]_{+} e^{i
\varphi_{\infty} /h} ,
\end{equation}
satisfies an estimate like \eqref{ram1}, with $\delta$ instead of $3 \delta$. As in \cite[Lemma 4.1]{hsmw3}, $A_{+}$ satisfies

\begin{proposition} \label{conv}\sl 
For all $(\alpha , \beta , \gamma ) \in \N^{1+d+(d-1)}$ and $N>0$, we have, uniformly with respect to $z \in [-{{C_0}}h , {{C_0}}h] + i [- {{C_1}} h , {{C_1}}h]$,
\begin{equation}
\big\vert \partial_{t}^{\alpha} \partial_{x}^{\beta} \partial_{\eta
'}^{\gamma} (h D_{t} + P(x,hD) - z) A_{+} \big\vert \leq C_{\alpha ,
\beta , N} h^{N} e^{- \delta t}
\end{equation}
\end{proposition}

\begin{proof}
The main difference with \cite[Lemma 4.1]{hsmw3} is that, here, $P$ is a
pseudodifferential operator. Let $c (t , x, \eta ' , z, h)$ be an
expandible symbol like $b$ (see \eqref{ram5}). From the definition of $[ c_{k} ]_{+}$ given by \eqref{ram7}, we have
\begin{align}
\partial_{t} [ c_{k} ]_{+} = [ \partial_{t} c_{k} ]_{+} \quad \text{and} \quad \partial_{\eta '} [ c_{k} ]_{+} = [ \partial_{\eta '} c_{k} ]_{+},
\end{align}
so that
\begin{align}
\partial_{t} [c]_{+} - [ \partial_{t} c ]_{+} \in \CS^{\infty, 2
\delta}    \mbox{  and  }
\partial_{\eta '} [c]_{+} - [ \partial_{\eta '} c ]_{+} \in \CS^{\infty, 2
\delta} .
\end{align}

Let $Q$ be a pseudodifferential operator with classical symbol $q (x , \eta ',
\xi ,z,h) \in \CS_{h}^{0} (1)$ that doesn't
depend on $t$. Then, there exist $(Q_{\widetilde{k}})_{\widetilde{k}
\in \N}$, a family of differential operators
in $x$ with $S^{0} (1)$ coefficients, such that, for all $d (t,x,
\eta ',z,h) \in \CS^{A,B}$ with $A,B>0$,
\begin{equation}
Q d = \sum_{\widetilde{k} \geq 0} ( Q_{\widetilde{k}} d )
h^{\widetilde{k}} \quad \text{modulo } \CS^{\infty, B}.  \label{ram10}
\end{equation}
Moreover, if $d$ is a classical expandible symbol, $Q d$ is also a classical expandible symbol.

Using this property with the $c_{k}$'s, we get
\begin{align}
Q c \sim& \sum_{k \geq 0} Q  c_{k} h^{k}  \quad\text{modulo }
\CS^{\infty, \re S}   \nonumber \\
\sim& \sum_{l \geq 0} \Big( \sum_{k+\widetilde{k} =l} Q_{\widetilde{k}}
c_{k} \Big) h^{l}   \quad\text{modulo } \CS^{\infty, \re S}.   \label{ram11}
\end{align}
Since $Q$ doesn't depend on $t$, we have
\begin{equation} \label{ram12}
[ Q_{\widetilde{k}} c_{k}]_{+} = Q_{\widetilde{k}} [c_{k}]_{+}.
\end{equation}
Then, \eqref{ram11} and \eqref{ram12} imply that $Qc$ is a classical expandible symbol and
\begin{align}
[ Q c]_{+} \sim& \sum_{l \geq 0} \Big[ \Big( \sum_{k+\widetilde{k} =l} Q_{\widetilde{k}}
c_{k} \Big) \Big]_{+} h^{l}  \quad\text{modulo } \CS^{\infty, 2 \delta}   \nonumber   \\
\sim& \sum_{l \geq 0} \Big( \sum_{k+\widetilde{k} =l} [ Q_{\widetilde{k}}
c_{k} ]_{+}  \Big) h^{l}  \quad\text{modulo } \CS^{\infty, 2 \delta}   \nonumber    \\
\sim& \sum_{l \geq 0} \Big( \sum_{k+\widetilde{k} =l}
Q_{\widetilde{k}} [ c_{k} ]_{+}  \Big) h^{l}   \quad\text{modulo }
\CS^{\infty, 2 \delta}   \nonumber   \\
\sim& Q [c]_{+} \quad \text{modulo } \CS^{\infty, 2 \delta}.  \label{ram31}
\end{align}
It follows that $[ Q c]_{-} \sim Q [c]_{-}$ modulo $\CS^{\infty, 2 \delta}$.

Let $q (x, \eta ', \xi ,z,h) \in S^{0}_{h}(1)$ be the (time independent) symbol of the pseudodifferential operator
\begin{equation}
Q = e^{-i \varphi_{\infty} /h} P(x,hD) e^{i \varphi_{\infty} /h} .
\end{equation}
From \eqref{ram9}, we get, for all $\varepsilon , N>0$,
\begin{align*}
\big\vert \partial_{t}^{\alpha} \partial_{x}^{\beta} \partial_{\eta '}^{\gamma} (h D_{t} + Q - z) a e^{i\varphi_{\star} /h} \big\vert \lesssim h^{N} e^{-(\re S + \varepsilon )t}.
\end{align*}
This estimate, combined with \eqref{ram1}, gives
\begin{equation}
\big\vert \partial_{t}^{\alpha} \partial_{x}^{\beta} \partial_{\eta '}^{\gamma} (h D_{t} + Q - z) b \big\vert \lesssim h^{-1-K_{1} - \vert \alpha \vert - \vert \beta \vert - \vert \gamma \vert} e^{-\delta t} + h^{N} e^{-(\re S + \varepsilon )t}  .  \label{ram21}
\end{equation}
Since $b$ is a classical expandible symbol, $d =\partial_{t}^{\alpha} \partial_{x}^{\beta} \partial_{\eta '}^{\gamma} (h D_{t} + Q - z) b$ is also a classical expandible symbol. Then \eqref{ram21} implies the

\begin{lemma}\sl
We have
\begin{equation}
\big[ \partial_{t}^{\alpha} \partial_{x}^{\beta} \partial_{\eta '}^{\gamma} (h D_{t} + Q - z) b \big]_{-} = 0 \quad \text{modulo } \CS^{\infty, 2 \delta} .  \label{ram30}
\end{equation}
\end{lemma}

\begin{proof}
If there exists $k$ such that $[d_{k} ]_{-} \neq 0$, we set $\widehat{\jmath} < J_{1}$  the first index such that there exists $k$ with $d_{k,\widehat{\jmath}} \neq 0$. Then, let $\widehat{k}$ be the first index with $d_{\widehat{k}, \widehat{\jmath}} \neq 0$. Using \eqref{ram21}, we get that, for all $N>0$ and $\varepsilon >0$,
\begin{align}
\vert d_{\widehat{k}, \widehat{\jmath}} \vert \lesssim& h^{-C} e^{(\lambda_{\widehat{\jmath}} -3 \delta) t} + h^{N} e^{Ct} + \sum_{k< \widehat{k}} \vert d_{k} \vert h^{k-\widehat{k}} e^{\lambda_{\widehat{\jmath}} t} + h e^{\varepsilon t}  \nonumber \\
\lesssim& h^{-C} e^{(\lambda_{\widehat{\jmath}} -3 \delta) t} + h^{N} e^{Ct} + h^{-C} e^{(\lambda_{\widehat{\jmath}} - \lambda_{\widehat{\jmath}+1} + \varepsilon) t} + h e^{\varepsilon t},
\end{align}
where the constant $C$ doesn't depend on $N, \varepsilon , t ,h , x, \eta ' ,z$. Notice that $\lambda_{\widehat{\jmath}} -3 \delta <0$ and $\lambda_{\widehat{\jmath}} - \lambda_{\widehat{\jmath}+1} <0$. Taking $h = e^{-\mu t}$ with $\mu >0$ small enough, we get
\begin{align}
\vert d_{\widehat{k}, \widehat{\jmath}} \vert \lesssim e^{- \mu t /2},
\end{align}
for $\varepsilon$ small enough and $N$ large enough. This implies $d_{\widehat{k}, \widehat{\jmath}} =0$, and this is a contradiction.
\end{proof}

Now we finish the proof of Proposition \ref{conv}. Using \eqref{ram9}, \eqref{ram4}, \eqref{ram31} and \eqref{ram30}, we get
\begin{align}
\big\vert \partial_{t}^{\alpha} \partial_{x}^{\beta} \partial_{\eta
'}^{\gamma} (h D_{t} + & P(x,hD) - z) A_{+} \big\vert    \nonumber  \\
\lesssim& h^{N} e^{(\re S + \varepsilon)t} + \big\vert \partial_{t}^{\alpha} \partial_{x}^{\beta} \partial_{\eta '}^{\gamma} (h D_{t} + P(x,hD) - z) [b]_{-} e^{i \varphi_{\infty} /h} \big\vert    \nonumber   \\
=& h^{N} e^{(\re S + \varepsilon)t} + \big\vert \partial_{t}^{\alpha} \partial_{x}^{\beta} \partial_{\eta '}^{\gamma} (h D_{t} + Q(x,hD) - z) [b]_{-} \big\vert    \nonumber  \\
\lesssim& h^{N} e^{(\re S + \varepsilon)t} . \label{ram40}
\end{align}
The proposition follows, taking a geometric mean between the two estimates \eqref{ram1} and \eqref{ram40}.
\end{proof}

Recalling that the functions
\begin{equation}
b_{k,j} (t, x, \eta ' ,z) = \sum_{l} b_{k,j,l} (x, \eta ' ,z) t^{l},
\end{equation}
are polynomial with respect to $t$, we can find a function $A_-$, holomorphic with respect to $z \in [-{{C_0}}h , {{C_0}}h] + i [- {{C_1}} h , {{C_1}}h]$, such that
\begin{equation}  \label{ram61}
A_{-} (x, \eta ' ,z ,h) \sim  \sum_{k \geq 1-K_{1}} h^{k}
\sum_{j<J_{1}, l} \frac{l!}{(S + \mu_{j})^{l+1}} b_{k,j,l} (x, \eta '
,z).
%+ \int_{-1}^{0} \chi (t) [ b ]_{-} (t, x, \eta ' ,z ,h) dt ,
\end{equation}
Notice that, formally, $A_{-} = \int_{0}^{+ \infty} \chi (t) [ b ]_{-} (t, x, \eta ' ,z ,h) dt$. At last, we set
\begin{equation}
u (x, \eta ',z,h) = A_{-} (x, \eta ' ,z ,h) + \int_{-1}^{+ \infty} \chi (t) A_{+} (t,x, \eta ' ,z ,h) dt ,
\end{equation}
and we have

\begin{proposition}[see {\cite[Proposition 4.2]{hsmw3}}]\sl
The function $u (x, \eta ',z,h)$ is holomorphic with respect to $z \in [-{{C_0}}h , {{C_0}}h] + i [- {{C_1}} h , {{C_1}}h]$ and satisfies, for all $( \beta , \gamma ) \in \N^{d + (d-1)}$,
\begin{eqnarray}
&\partial_{x}^{\beta} \partial_{\eta '}^{\gamma} u =
\CO \big( h^{-K_{1} - \vert \beta \vert - \vert \gamma \vert} \big),  \label{ram60} \\
&\partial_{x}^{\beta} \partial_{\eta '}^{\gamma} (
P(x,hD) - z) u = \CO ( h^{\infty} ),   \label{ram52}
\end{eqnarray}
for $x \in \bigcup_{t>-1/2} \widetilde{U}_{t}$ and $\eta ' \in
V_{-}$. Moreover, for $x \in H_{-}$, we have
\begin{equation}  \label{ram53}
u = ( 1 + r (x, \eta ' ,z ,h) )e^{i x' \cdot \eta ' /h} ,
\end{equation}
where $r \in \CS^{\infty} (1)$.
\end{proposition}

\begin{proof}
The estimate \eqref{ram60} follows from \eqref{ram4} and \eqref{ram61}.
Now from \eqref{conv}, we get
\begin{align}
\partial_{x}^{\beta} \partial_{\eta '}^{\gamma} (P(x,hD) - z) \int_{-1}^{+ \infty} \chi A_{+} dt   =& \int_{-1}^{+ \infty} \partial_{x}^{\beta} \partial_{\eta '}^{\gamma} (hD_{t} + P(x,hD) - z)  \chi A_{+} dt  \nonumber \\
=& \CO (h^{\infty}) + \int_{-1}^{+ \infty} \partial_{x}^{\beta} \partial_{\eta '}^{\gamma} (hD_{t}  \chi) A_{+} dt
 \label{ram51}
\end{align}
so that the L.H.S.  is microlocally $0$ in $\Omega$
since $A_{+} =0$ microlocally in that set for $t \in \supp
(\partial_{t} \chi) \subset ]-1,-1/2[$.

On the other hand, from \eqref{ram30}, we have
\begin{align}
\sum_{1-K_{1} \leq \widetilde{k} \leq k} Q_{\widetilde{k}} b_{k -
\widetilde{k}, j ,l} - \frac{z}{h} b_{k-1,j,l} = i (l+1) b_{k-1,
j,l+1} - i (S+\mu_{j}) b_{k-1,j,l}.
\end{align}
Then
\begin{align}
(P-z) A_{-} e^{i \varphi_{\infty} /h} =& e^{i \varphi_{\infty} /h}
(Q-z) A_{-}    \nonumber  \\
\sim& e^{i \varphi_{\infty} /h} \sum_{k} h^{k} \sum_{j < J_{1},l}
\frac{l !}{(S+\mu_{j})^{l+1}} \Big( \sum_{1-K_{1}  \leq \widetilde{k} \leq k} Q_{\widetilde{k}} b_{k -
\widetilde{k}, j ,l} - \frac{z}{h} b_{k-1,j,l} \Big)  \nonumber  \\
\sim& i e^{i \varphi_{\infty} /h} \sum_{k} h^{k} \Big( \sum_{j < J_{1},l}
\frac{(l+1) !}{(S+\mu_{j})^{l+1}} b_{k-1,j,l+1} -  
\sum_{j < J_{1},l} \frac{l !}{(S+\mu_{j})^{l}} b_{k-1,j,l} \Big).\label{ram50}
\end{align}
Therefore $(P-z) A_{-} e^{i \varphi_{\infty} /h}=0$ microlocally in $\Omega$. One can also differentiate \eqref{ram50}, and obtain
the corresponding estimates. Then \eqref{ram52} follows from \eqref{ram51} and \eqref{ram50}.
Eventually, \eqref{ram53} follows from the fact that, for $x \in H_{-}$, $a$ has a compact support in $t$: 
 The formal stationary phase expansion (\ref{erp1}) can be given a meaning,  and gives  this last estimate.
\end{proof}

%-------------------------------------
%
%-------------------------------------
\section{The symbol of the transition operator} 
\label{seclocalscatt}

Now we finish the proof of  Theorem  \ref{explicit}. We compute  the principal symbol of  the operator $\CJ(z)$, defined after Theorem \ref{exis}, that is the microlocal value of the solution $u$ in Theorem \ref{repres} at some point $\rho_+\in \Lambda_+\setminus\widetilde{\Lambda_{+}}(\rho_-)$ (see the definition after (\ref{transpphi1})).

As in \cite[Section 5]{hsmw3}, we can assume $K_{1} =0$  (see (\ref{K1gh})) since the general case can be treated the same way. In that case, we recall that  the solution $u$ of the problem (\ref{zeproblem}) can be written as
\begin{equation}\label{recall1}
u(x,h)= \frac{1}{(2 \pi h)^{d-\frac{1}{2}}} \iint_{T^*\R^{d-1}} \int_{-1}^{+\infty} 
e^{i( \varphi (t,x,\eta ' )- y' \eta ' )/h} a(t,x, \eta ' ,z, h) u_{0}(\varepsilon , y' )
dt dy' d\eta ' ,
\end{equation}
 where $\varphi$ is defined in Section \ref{subsec:phase} and has the properties given in \eqref{pi1}--\eqref{pi2}, and $a$ is the symbol described in Section \ref{subsec:symbol}.
 
First of all, we  compute the principal  term $a_{0}$ of the symbol $a$ in (\ref{recall1}).
Performing  again a formal stationary phase with respect to $t$ in (\ref{solution}), we obtain, for $x=(\varepsilon ,x') \in H_{-}$,
\begin{equation}  \label{erp}
u (x,h) = \frac{1}{(2 \pi h)^{d-1}} \iint_{T^* \R^{d-1}} e^{i (x \cdot \eta '- y' \cdot \eta ' )/h} \widetilde{a} (x, \eta ',z,h) u_{0}(\varepsilon , y' ) dy' d\eta '  ,
\end{equation}
where $\widetilde{a}=1$ by  our choice in (\ref{condinit}). In particular for $x\in H_-$, we have
\begin{equation}
a_{0} (t(x, \eta '),x, \eta ',z) = e^{-i\pi /4} \vert
\partial_{tt}^{2} \varphi (t(x, \eta '), x, \eta ') \vert^{1/2} .
\end{equation}
Notice that we have done so that, microlocally near $\Lambda_-$,
\begin{equation}\label{}
\frac1{(2\pi h)^{d-1/2}}\int_{-1}^{+\infty} e^{i \varphi (t,x,\eta ' )/h} a(t,x, \eta ' ,h,z)  dt    =
b(x,\eta', h)e^{i\psi_{\eta'}(x)/h}
\end{equation}
where $\ds b(x,\eta',h)=(2 \pi h)^{-(d-1)} \sum_{j=0}^{\infty} h^{j} b_{j} (x,\eta'),$ is a symbol such that
\begin{equation}
\left\{ \begin{aligned}
&(P-z) \big( b(x,\eta',h) e^{i \psi_{\eta '} (x)/h} \big) = {\mathcal O} (h^{\infty})
&&\text{ near } \gamma_{\eta '}, \\
&b (x,\eta',h) =\tilde a (t(x, \eta '),x, \eta ',z)=1 &&\text{ on } H_{-},
\end{aligned} \right.
\end{equation}
%\begin{eqnarray}
%\nonumber
%&&
%b(x,\eta', h) = \\
%&&
%\qquad (2 \pi h)^{-\frac{d}{2}} e^{-i \psi_{\eta '} (x)/h}
%\int_{-1}^{+\infty} e^{i \varphi (t,x,\eta ' )/h} a(t,x, \eta ' ,h,z)  dt    
%= (2 \pi h)^{- \frac{d-1}{2}} \sum_{j=0}^{\infty} h^{j} b_{j} (x,\eta'),
%\end{eqnarray}
%where
%\begin{equation}  \label{tt3}
%b_{0} (x,\eta') = e^{i\pi /4} \frac{a_{0} (t(x, \eta '), x, \eta ',z,h)}{\vert \partial_{tt}^{2} \varphi (t(x, \eta '), x, \eta ') \vert^{1/2}}\cdotp
%\end{equation}
%Since $u$ is a solution of (\ref{zeproblem}), we have
%\begin{equation}
%\left\{ \begin{aligned}
%&(P-z) \big( b(x,\eta',h) e^{i \psi_{\eta '} (x)/h} \big) = {\mathcal O} (h^{\infty})
%&&\text{ near } \gamma_{\eta '}, \\
%&b_{0} (x,\eta',h) = a (t(x, \eta '),x, \eta ',z) &&\text{ on } H_{-},
%\end{aligned} \right.
%\end{equation}
The principal symbol $b_{0}$ of $b$ satisfies 
\begin{equation}  \label{tt3}
b_{0} (x,\eta') = e^{i\pi /4} \frac{a_{0} (t(x, \eta '), x, \eta ',z,h)}{\vert \partial_{tt}^{2} \varphi (t(x, \eta '), x, \eta ') \vert^{1/2}},
\end{equation}
and it is a solution of  the first transport equation
\begin{equation}  \label{op1}
\left\{ \begin{aligned}
&\partial_{\xi} p_{0} (x, \partial_{x} \psi_{\eta '}) \partial_{x} b_{0} + \Big( \frac{1}{2} \tr \big( \partial^{2}_{\xi , \xi} p_{0} (x, \partial_{x} \psi_{\eta '}) \partial^{2}_{x,x} \psi_{\eta '} \big) -i\frac{z}h \Big) b_{0} = 0 &&\text{ near } \gamma_{\eta '}, \\
&b_{0} (x,\eta') = \widetilde{a}_{0} (x', \eta ')=1 &&\text{ on } H_{-}.
\end{aligned} \right.
\end{equation}
In the Schr\"{o}dinger case, the first equation of \eqref{op1} can be written as
\begin{equation*}
2\nabla \psi_{\eta '} \cdot \nabla b_{0} + ( \Delta \psi_{\eta '} ) b_{0} -\frac{iz}{h}
b_{0} = 0 \text{ near } \gamma_{\eta '}.
\end{equation*}

We calculate $b_{0}$, starting with the computation of the trace in (\ref{op1}) (as e.g. in the book of V. Maslov and  M. Fedoryuk \cite{Ma}). Let $( x(t,x ' ,\eta'), \xi (t,x ',\eta'
))$ be the Hamiltonian curve with initial condition 
\begin{equation}
( x(0,x ' ,\eta'), \xi (0,x ',\eta'
))=(x_{1}^- , x', f_{-} ( x_{1}^-,x' , \eta ' ) , \eta ' ) \in \Lambda_{\psi_{\eta '}} .
\end{equation}
With the notations of Lemma \ref{intersectionpropre}, this curve is $\gamma_{\eta'}$ when $x'=x'(\eta')$.
As usual, we have
\begin{equation}
\partial_{t} \ln \det \frac{\partial x(t,x',\eta')}{\partial (t,x')} = \tr \big( \partial^{2}_{\xi , \xi} p_{0} (x, \partial_{x} \psi_{\eta '}) \partial^{2}_{x,x} \psi_{\eta '} \big) ,
\label{laplacepsi}
\end{equation}
and then (\ref{op1}) becomes
\begin{equation}
\partial_{t} \left( \sqrt{\det \frac{\partial x(t,x',\eta')}{\partial (t,x')}} b_{0} (x(t,x',\eta'),\eta') \right) 
= \frac{iz}{h} \sqrt{\det \frac{\partial x(t,x',\eta')}{\partial (t,x')}} b_{0} (x(t,x',\eta'),\eta'),
\end{equation}
which gives 
\begin{equation}  \label{tt1}
b_{0} (x(t,x',\eta'),\eta') = \frac{\sqrt{\partial_{\xi_{1}} p_{0} (x_{1}^- , x', f_{-} ( x_{1}^-,x' , \eta ' ) , \eta ' )}}{\sqrt{\det \frac{\partial x(t,x',\eta')}{\partial (t,x')}}} e^{itz/h}.
\end{equation}

We are interested in taking the limit $t\to +\infty$ in this expression. The point is that, as $t \rightarrow + \infty$,
\begin{equation}
\frac{1}{2} \partial_{t} \ln \det \frac{\partial x(t,x',\eta')}{\partial (t,x')}=(\sum \lambda_{j} /2 -\lambda_{1})t +o(1).
\end{equation}
Indeed, starting from (\ref{laplacepsi}), we have
\begin{equation}
\partial^{2}_{\xi , \xi} p_{0} (x, \partial_{x} \psi_{\eta '}) = 2 + {\mathcal O}
(e^{-\lambda_{1} t}) ,
\end{equation}
as a matrix, for $x \in \gamma_{\eta '} (t)$. On the other hand, writing $\psi_{\eta '} (x) = \varphi (t(x, \eta ')
, x, \eta ')$  and using the fact that $( \partial_{t} \varphi ) (t(x, \eta ')
, x, \eta ') =0$, we get
\begin{equation}
\partial_{x_{j}} \psi_{\eta'} = ( \partial_{x_{j}} \varphi ) (t(x, \eta ')
, x, \eta ') ,
\end{equation}
so that
\begin{equation}
\partial^{2}_{{x_{j}},{x_{k}}} \psi_{\eta '} = \partial^{2}_{{x_{j}},{x_{k}}} \varphi + ( \partial^{2}_{t,{x_{j}}} \varphi )  ( \partial_{x_{k}} t ).  \label{ee1}
\end{equation}
Now using (\ref{pi1}),  we have, 
\begin{equation}
\partial^{2}_{x,x} \varphi = \partial^{2}_{x,x}\varphi_{+} + {\mathcal O}
(e^{-\lambda_{1} t}) = \frac{L}{2} + {\mathcal O}
(e^{-\lambda_{1} t}),  \label{ee2}
\end{equation}
and
\begin{align}
\nonumber
\partial_{t, x_{j}}^{2} \varphi =& -\lambda_{1} ( \partial_{x_{j}}
\varphi_{1} ) e^{-\lambda_{1} t} + {\mathcal O} (e^{-\mu_{2} t})  \\
=&  \label{ee3}
\left\{ \begin{aligned}
&\lambda_{1}^{2} \vert g_{1} \vert e^{-\lambda_{1} t} + {\mathcal O}
(e^{-\mu_{2} t}) \quad &&\text{ if } j=1, \\
&{\mathcal O} (e^{-\mu_{2} t}) \quad &&\text{ if } j \neq 1.
\end{aligned} \right.  
\end{align}
Then using again the fact that $( \partial_{t} \varphi ) (t(x, \eta '), x, \eta ') =0$, we get
\begin{equation}  \label{ee4}
\partial_{x_{k}} t =
\left\{ \begin{aligned}
&- \vert g_{1} \vert^{-1} \lambda_{1}^{-1} e^{\lambda_{1} t} + {\mathcal O} (e^{(\lambda_{1}-\widehat{\mu}_{1} )t}) \quad &&\text{ if } k=1, \\
&{\mathcal O} (e^{(\lambda_{1}-\widehat{\mu}_{1} )t}) \quad &&\text{ if } k \neq 1.
\end{aligned} \right.
\end{equation}
Using (\ref{ee3}) and the estimates (\ref{ee2}), (\ref{ee3}), (\ref{ee4}), we obtain
\begin{equation}  \label{tt2}
\frac{1}{2} \tr \big( \partial^{2}_{\xi , \xi} p_{0} (x, \partial_{x} \psi_{\eta '}) \partial^{2}_{x,x} \psi_{\eta '} \big) = \sum_{j=1}^{d} \lambda_{j} /2 -\lambda_{1} + {\mathcal O} (e^{-\widehat{\mu}_{1}t}),
\end{equation}
on $\gamma_{\eta '} (t)$. Therefore, we shall write (\ref{tt1}) as
\begin{align}
b_{0} (x(t,x',\eta'),\eta') = \ds  e^{( \lambda_{1} -S(z/h))t} & e^{(\sum \lambda_{j} /2 -\lambda_{1})t - \frac{1}{2} \ln \det \frac{\partial x(t,x',\eta')}{\partial (t,x')}}  \nonumber   \\
&\times \sqrt{\partial_{\xi_{1}} p_{0} (x_{1}^- , x', f_{-} ( x_{1}^-,x' , \eta ' ) , \eta ' )} .   \label{tt7}
\end{align}

%\begin{align}
%\nonumber
%\partial^{2}_{t, x_{j}} \varphi =& -\lambda_{1} (\partial_{x_{j}} \varphi_{1}) e^{- \lambda_{1} t}
%+ {\mathcal O}(e^{-\mu_{2}t} )  \\
%=& \left\{ \begin{aligned}
%&\vert g_{1} \vert \lambda_{1}^{2} e^{-\lambda_{1} t} + {\mathcal O} (e^{-\mu_{2} t}) \quad
%&&\text{ if } j=1, \\
%&{\mathcal O} (e^{-\mu_{2} t}) \quad &&\text{ if } j \neq 1,
%\end{aligned} \right.
%\end{align}

Now  we compute $(\partial^{2}_{t,t} \varphi) (t(x, \eta '), x, \eta ')$ on the curve $\gamma_{\eta '}$. We have $\partial_{t} \varphi = - p(x, \partial_{x} \varphi)$, and 
\begin{equation}  \label{tt4}
\partial^{2}_{t,t} \varphi = -2 ( \partial_{x} \varphi ) \cdot ( \partial_{t,x} \varphi ).
\end{equation}
But, on $\gamma_{\eta '}$, we have \eqref{ee3} and
\begin{align}
\nonumber
\partial_{x_{j}} \varphi =& \partial_{x_{j}} \varphi_{+} + \partial_{x_{j}} \varphi_{1} e^{-\lambda_{1}t} + {\mathcal O} (e^{-\mu_{2} t}) \\
=& \left\{ \begin{aligned}
&- \frac{\vert g_{1} \vert \lambda_{1}}{2} e^{-\lambda_{1} t} + {\mathcal O} (e^{-\mu_{2} t}) \quad &&\text{ if } j=1, \\
&{\mathcal O} (e^{-\mu_{2} t}) \quad &&\text{ if } j \neq 1.
\end{aligned} \right.
\end{align}
Therefore, (\ref{tt4}) becomes
\begin{equation}  \label{tt5}
\partial^{2}_{t,t} \varphi (t, x(t,x',\eta'), \eta ') = \vert g_{1} \vert^{2} \lambda_{1}^{3} e^{-2 \lambda_{1}t} + {\mathcal O} (e^{-(\lambda_{1} + \mu_{2})t}).
\end{equation}

We recall that $a_{0}$ is expandible, namely
\begin{equation} \label{tt6}
a_{0}(t, x , \eta ') \sim \sum_{j=0}^{\infty} a_{0, j} (t, x, \eta ') e^{-(S(z/h) + \mu_{j})t} ,
\end{equation}
where $a_{0,j}$ are polynomials with respect to $t$, and $a_{0,0}$ does not depend on $t$. Using (\ref{tt3}), (\ref{tt7}), (\ref{tt5}) and (\ref{tt6}), we get
\begin{equation}
\begin{aligned}
a_{0,0} (x_{\eta '} (t), \eta ') = \vert g_{1} \vert \lambda_{1}^{3/2} e^{-i\pi /4} & \sqrt{\partial_{\xi_{1}} p_{0} (x_{1}^- , x', f_{-} ( x_{1}^-,x' , \eta ' ) , \eta ' )} \\
&\times e^{(\sum \lambda_{k} /2 -\lambda_{1})t - \frac{1}{2} \ln \det \frac{\partial x(t,x')}{\partial (t,x')}} b_{0} (x_{\eta '},\eta') + {\mathcal O} (e^{-\widehat{\mu}_{1} t}),
\end{aligned}
\end{equation}
where $\widehat{\mu}_1=\mu_2-\mu_1$, and then, since $x_{\eta'}\in H_-$,
\begin{equation}
a_{0,0} (0, \eta ') = \vert g_{1} \vert \lambda_{1}^{3/2} e^{-i\pi /4} \sqrt{\partial_{\xi_{1}} p_{0} (x_{1}^- , x', f_{-} ( x_{1}^-,x' , \eta ' ) , \eta ' )} \lim_{t\to + \infty} \frac{e^{(\sum \lambda_{k} /2 -\lambda_{1})t}}{{\sqrt{\det \frac{\partial x(t,x')}{\partial (t,x')}}}}.
\label{yyyyyyyyy}
\end{equation}
Notice that  the above  limit exists thanks to  \eqref{tt2}.

Finally we compute the solution $u(x,h)$ given by (\ref{solution}) microlocally near $\rho_{+}$. Since $\rho_{+} \in \Lambda_{+} \setminus \widetilde{\Lambda_{+}}(\rho_{-})$, we can use the calculus of \cite[Section 5]{hsmw3} and we get, microlocally near $\rho_{+}$,
\begin{equation}
\int_{-1}^{+\infty}e^{i\varphi(t,x,\eta')/h}a(t,x,\eta',z,h) dt = c(x,\eta',h) e^{i (\varphi_{+} (x)+\tilde \psi(\eta')) /h} .
\end{equation}
Here $c(x,\eta',h)$ is a symbol of class $\CS^{0}_{h}$ which satisfies
\begin{equation}
c(x,\eta',h) \sim \sum_{j=0}^{\infty} c_{j} (x,\eta', \ln h) h^{S(z/h)/\lambda_{1} + \widehat{\mu}_{j}/\lambda_{1}},
\end{equation}
where the $c_{j} (x,\eta', \ln h)$ are polynomial with respect to $\ln h$ and, in particular, 
\begin{equation}  \label{ttyy}
c_{0} (x,\eta') = \frac{1}{\lambda_{1}} (\varphi_{1} (x) /i)^{-S(z/h)/\lambda_{1}} \Gamma \left( S(z/h) /\lambda_{1} \right) a_{0,0} (x,\eta'),
\end{equation}
doesn't depend on $\ln h$. Here $\Gamma$ denotes Euler's Gamma function, and $(\widehat{\mu}_j)_{j\geq 0}$ is the increasing sequence  of the linear combinations over $\N$ of the $(\mu_{k}-\mu_1)$'s, $k\geq 2$.

On the other hand, since we want that the function $u(x,h)$, given by
\begin{equation}
u(x,h) = \frac1{(2\pi h)^{d-1/2}}\iint_{T^*\R^{d-1}} c(x,\eta',h) e^{i (\varphi_{+} (x) +\tilde \psi(\eta')-y'\cdot \eta')/h}u_0(y)dy'd\eta',
\label{zfinal}
\end{equation}
is a microlocal solution of $(P-z) u=0$ for any initial data $u_0$, the function $c_0$ should satisfy the usual transport equation:
\begin{equation}
\partial_{\xi} p_{0} (x, \partial_{x} \varphi_{+} ) \partial_{x} c_{0} + \Big( \frac{1}{2} \tr \big( \partial^{2}_{\xi , \xi} p_{0} (x, \partial_{x} \varphi_{+} ) \partial^{2}_{x,x} \varphi_{+} \big) -i\frac{z}h \Big) c_{0} = 0.
\end{equation}
Thus,  if $(x(t), \xi (t))$ is the  integral curve of $H_{p}$ in $\Lambda_{+}$ with initial condition $\rho = (x, \nabla \varphi_{+} (x))$, we have
\begin{equation}  \label{rfv}
c_{0} (x(t),\eta')=\ds e^{itz/h - \frac{1}{2} \int_{0}^{t} \tr ( \partial^{2}_{\xi , \xi} p_{0} ( \cdot , \partial_{x} \varphi_{+} ) \partial^{2}_{x,x} \varphi_{+} ) (x(s)) ds} c_{0} (x,\eta').
\end{equation}
Let us  compute $c_{0} (x(t),\eta')$ using (\ref{ttyy}). Since $\rho_{+} \notin \widetilde{\Lambda_{+}}(\rho_{-})$, we can assume that $\rho \notin \widetilde{\Lambda_{+}}(\rho_{-})$ for $\rho$ close enough to $\rho_{+}$. In particular $\rho \notin \widetilde{\Lambda_{+}}$ and then
\begin{equation}
x(t) \sim \sum_{j=1}^{\infty} g_{j}^{+} (t) e^{\mu_{j}t},
\end{equation}
as $t \to - \infty$, where the $g_{j}^{+} (t)$ are polynomials with respect to $t$ and $g_{1}^{+} (t)$ doesn't depend on $t$. 
%Notice that if  $p = \xi^{2}+ V(x)$, we have $g_{j}^{+} (t) = g_{j}^{-} (-t)$. Now
\begin{equation}
\varphi_{1}(x(t)) = -\lambda_{1} (g^-_{1} (\rho_{\eta '}) \cdot g^+_{1}(\rho)) e^{\lambda_{1} t} + {\mathcal O} (e^{(\mu_{2} -\varepsilon)t}),
\end{equation}
and the equations (\ref{ttyy}) and (\ref{rfv}) give
\begin{align}
\nonumber
c_{0} (x,\eta') =& \frac{1}{\lambda_{1}} e^{\frac{1}{2} \int_{0}^{t} ( \tr ( \partial^{2}_{\xi , \xi} p_{0} ( \cdot , \partial_{x} \varphi_{+} ) \partial^{2}_{x,x} \varphi_{+} ) (x(s)) - \sum \lambda_{l} ) ds} (i\lambda_{1} (g^-_{1} (\rho_{\eta '}) \cdot g^+_{1}(\rho)))^{-S(z/h)/\lambda_{1}}   \\
\nonumber
&\qquad \qquad \Gamma \left( S(z/h) /\lambda_{1} \right) a_{0,0} (x(t),\eta') + {\mathcal O} (e^{\varepsilon t})  \\
\nonumber
=& \frac{1}{\lambda_{1} } e^{\frac{1}{2} \int_{0}^{- \infty} ( \tr ( \partial^{2}_{\xi , \xi} p_{0} ( \cdot , \partial_{x} \varphi_{+} ) \partial^{2}_{x,x} \varphi_{+} ) (x(s)) - \sum \lambda_{j} ) ds} (i\lambda_{1} (g^-_{1} (\rho_{\eta '}) \cdot g^+_{1}(\rho)))^{-S(z/h)/\lambda_{1}}   \\
&\qquad \qquad \Gamma \left( S(z/h) /\lambda_{1} \right) a_{0,0} (0, \eta ').
\label{zzzzzzz}
\end{align}

At last,  we go back to (\ref{zfinal}) and we perform a stationary phase expansion with respect to $\eta'$ in that integral. Recalling (\ref{pi2}), we can write
\begin{equation}\label{zfinal2}
u(x,h)=\frac{1}{(2\pi h)^{d-1/2}}\iint_{T^*\R^{d-1}}e^{i\varphi(x,\eta',y')/h}c(x,\eta',h)u_0(y)dy' d\eta',
\end{equation}
where
\begin{equation}\label{zfinal3}
\varphi(x,\eta',y')=\varphi_+(x)+(x'(\eta')-y')\cdot\eta'-\varphi_-(x(\eta')).
\end{equation}
We have
\begin{equation}\label{zfinal4}
\nabla_{\eta'}\varphi(x,\eta',y')=
(x'(\eta')-y')+\nabla_{\eta'}x'(\eta')\cdot (\eta'-\nabla_{x'}\varphi_-(x(\eta')),
\end{equation}
since $x(\eta')=(x_1^-,x'(\eta'))$ where $x_1^-$ does not depend on $\eta'$. But $\rho(\eta')=(x(\eta'),\xi(\eta'))$ belongs to $\Lambda_-$ (see (\ref{rhoeta})), so that $\nabla\varphi_-(x(\eta'))=\xi(\eta')$, and in particular $\nabla_{x'}\varphi_-(x(\eta'))=\eta'$. Thus the last term in (\ref{zfinal4}) vanishes, and  $\eta'\mapsto \varphi(x,\eta',y')$ has a unique critical point $\eta'(y')$,  such that $y'=x'(\eta'(y'))$, with critical value
\begin{equation}\label{zfinal5}
\tilde\varphi(x,y')=\varphi(x,\eta'(y'),y')=\varphi_+(x)-\varphi_-(\varepsilon, y').
\end{equation}
Moreover, since $\nabla^2_{\eta' x'}\varphi_-(x(\eta'))=I$, we have
\begin{equation}\label{zfinal6}
\nabla^2_{\eta'\eta'}\varphi(x,\eta', y')=\nabla_{\eta'}x'(\eta')=
\left(
\nabla^2_{x'x'}\varphi_-(x(\eta'))
\right )^{-1}
.
\end{equation}
Thus, there exists a symbol $d(x,y',z,h) \sim \sum_{j \geq 0} d_{j} (x,y',z, \ln h) h^{\widehat{\mu}_{j} / \lambda_{1}} \in \CS_{h}^{0} (1)$, with $d_{j} (x,y',z, \ln h)$ polynomial with respect to  $\ln h$, such that
\begin{equation}\label{zfinal7ter}
\CJ (z) u_0(x,h)=\frac{h^{S(z/h)/\lambda_{1}}}{(2\pi h)^{d/2}}\int_{\R^{d-1}} d(x,y',z,h) e^{i(\varphi_+(x)-\varphi_-(\varepsilon,y'))/h} u_0(\varepsilon,y') dy' ,
\end{equation}
microlocally near $\rho_+$.
%\begin{equation}\label{zfinal7bis}
%\CJ(z)u_0(x,h)=\frac{h^{S(z/h)/\lambda_{1}}}{(2\pi h)^{d/2}}\int_{\R^{d-1}}
%d(x,y',h)e^{i(\varphi_+(x)-\varphi_-(\varepsilon,y'))/h} u_0(\varepsilon,y') dy'.
%\end{equation}
Moreover the principal symbol $d_0$ of $d$ is independent of $\ln h$, and can be written as
\begin{equation}\label{zfinal8}
d_0(x,y',z)={e^{-i(d-1)\pi/4}}\vert\det \nabla^2_{y'y'}\varphi_-(\varepsilon,y')\vert^{1/2}c_0(x,\eta'(y')),
\end{equation}
where $c_0(x,\eta'(y'))$ is given  by (\ref{zzzzzzz}) and (\ref{yyyyyyyyy}), and this finishes the proof Theorem \ref{explicit}.

\appendix
\section{A review of $h$-pseudodifferential calculus}
\label{pdo}

One of the main tool of this paper is the so-called $h$-pseudodifferential calculus, and we review here some basic facts. Since we deal with self-adjoint operators and spectral properties, we shall only use Weyl quantization. First we recall this calculus in standard classes of symbols, following closely \cite[Chapter 7]{disj}  (see also \cite{mabk}).

We say that $m : T^{*} \R^{d} \to [0, + \infty[$ is an order function when there are $C$, $N >0$ such that $m (x) \leq C \langle x-y \rangle^{N} m (y)$.

If $m (x,\xi)$ is an order function, and $\delta \geq 0$ a real number, we say that a function $a (x,\xi,h) \in C^{\infty} (T^{*} \R^{d} )$ is a symbol of class $\CS_{h}^{\delta} (m)$ when
\begin{equation}
\forall \alpha\in\N^{2d},\  \exists C_{\alpha}>0 ,\  \forall h\in ]0,1] , \   |\partial_{x,\xi}^\alpha a(x,\xi,h)| \leq C_{\alpha} h^{-\delta |\alpha |} m(x,\xi ). \label{shdelta1}
\end{equation}
If $e(h)$ is a function of $h$ only, sometimes we write  $\CS^\delta_{h}(e(h) m)$ instead of $e(h)\CS^\delta_{h}(m)$.

If $a (x,\xi ,h)$ is a symbol of class $\CS_{h}^{\delta} (m)$, we define the $h$-pseudodifferential operator  ${\rm Op}_{h} (a) $ with symbol $a$ by
\begin{equation}
\forall u \in \CC_{0}^{\infty} (\R^{d}),\  \left( {\rm Op}_{h} (a) u \right) (x) = \frac{1}{(2 \pi h)^{n}} \iint e^{i(x-y)\cdot \xi /h}a \Big( \frac{x+y} 2, \xi \Big) u(y) dy d \xi .
\label{oph}
\end{equation}
We also denote by $\Psi_{h}^{\delta} (m)$ the space of operators ${\rm Op}_{h} (\CS^\delta_{h}(m))$.

The composition rule between pseudodifferential operators in $\Psi_{h}^{\delta} (m)$ is given in the following proposition. It  is an easy  adaptation of Proposition 7.7 in \cite{disj}:
\begin{proposition}\sl
 If $a_{1}\in \CS^{\delta_{1}}_{h}(m_{1})$ and $a_{2}\in \CS^{\delta_{2}}_{h}(m_{2})$ with $0\leq \delta_{1},\delta_{2}\leq \frac12$ and $\delta_{1}+\delta_{2}<1$, then ${\rm Op}_{h} (a_{1}) \circ {\rm Op}_{h} (a_{2})$ belongs to $\Psi_{h}^{\max (\delta_{1}, \delta_{2})} (m_{1} m_{2})$, and, for any $N\in \N$,  its symbol $a_1\#a_{2}$ verifies 
\begin{eqnarray}
\nonumber
&&
(a_1\#a_{2})(x,\xi) =e^{\frac{ih}{2} \sigma(D_{x},D_{\xi}, D_{y},D_{\eta})} \big( a_1(x,\xi)a_{2}(y,\eta) \big) {\Big |_{y=x,\eta=\xi}} \\
\nonumber
&&
\qquad
=\sum_{k=0}^{N-1} \frac 1{k!} \Big( \big(\frac{ih}2 \sigma(D_{x},D_{\xi},D_{y},D_{\eta})\big)^k  a_1(x,\xi)a_{2}(y,\eta) \Big) \Big|_{y=x,\eta=\xi} \qquad\\
&&
\qquad\qquad 
+ h^{N(1-\delta_{1} - \delta_{2})} \CS^{\max(\delta_{1},\delta_{2})}_h (m_{1} m_{2}).
\end{eqnarray}
\label{a1a2}
\end{proposition}

Notice that in this theorem and below, we use the  standard notation $\sigma(D_{x},D_{\xi}, D_{y},D_{\eta})=D_{\xi}D_{y}-D_{x}D_{\eta}$.

To control the norm of a pseudodifferential operator in ${\mathcal L} (L^{2} (\R^{d}))$, we use the following classical result:
\begin{theorem}\sl (Calder\`on--Vaillancourt)
Let $a \in \CS_{h}^{\delta} (1)$ with $0\leq \delta \leq 1/2$. Then there exists $C>0$ such that
\begin{equation}
\forall u \in L^{2} (\R^{d}),\  \Vert {\rm Op}_{h} (a) u \Vert_{L^{2} (\R^{d})} \leq C \| u \|_{L^{2} (\R^{d})}.
\end{equation}
Furthermore, $C$ is bounded by a semi-norm of $a \in \CS_{h}^{\delta} (1)$.
\label{dss}
\end{theorem}

We now recall the semiclassical sharp G\aa rding inequality and  Fefferman--Phong's inequality:
\begin{theorem}\sl (G\aa rding's inequality)
Let $a (x,\xi ,h)$ be a real valued symbol in $\CS_{h}^{0} (1)$. If $a(x,\xi ,h) \geq 0$ for all $(x, \xi ,h) \in T^{*} \R^{d} \times [0,1]$, then there exists $C >0$ such that
\begin{equation}
\forall u \in L^{2} (\R^{d}),\  \big( {\rm Op}_{h} (a) u , u \big)_{L^{2} (\R^{d})} \geq -C h \| u \|^{2}_{L^{2} (\R^{d})}.
\end{equation}
Furthermore, $C$ is bounded by a semi-norm of $a \in \CS_{h}^{0} (1)$.
\label{gar}
\end{theorem}

\begin{theorem}\sl (Fefferman-Phong's inequality)
Let $a (x,\xi ,h)$ be a real valued symbol in $\CS_{h}^{0} (1)$. If $a(x,\xi ,h) \geq 0$ for all $(x, \xi ,h) \in T^{*} \R^{d} \times [0,1]$, then there exists $C >0$ such that
\begin{equation}
\forall u \in L^{2} (\R^{d}),\  \big( {\rm Op}_{h} (a) u , u \big)_{L^{2} (\R^{d})} \geq -C h^{2} \| u \|^{2}_{L^{2} (\R^{d})}.
\end{equation}
Furthermore, $C$ is bounded by a semi-norm of $a \in \CS_{h}^{0} (1)$.
\label{fff}
\end{theorem}

We now give the composition rule in the class $\widetilde{\CS}_{h}$ we use in Section 4, which can be seen as a particular case of the  semiclassical Weyl--H\"{o}rmander calculus. Let $m(x,\xi)$ be an order function. We say that a function $a (x,\xi ,h)$ is a symbol of class $\widetilde{\CS}_{h} (m)$ if $\forall \alpha , \beta \in \N^{d}$,   $\exists C_{\alpha , \beta}>0$ such that , $\forall h \ll 1$,
\begin{equation}
 \ |\partial_{x}^\alpha \partial_{\xi}^{\beta} a (x,\xi,h)| \leq C_{\alpha , \beta} m(x,\xi ) \langle x \rangle^{- \vert \alpha \vert /2} \langle \xi \rangle^{- \vert \beta \vert /2}. \label{shdelta2}
\end{equation}

Concerning the product rule,  we have the following result, which is similar to Proposition \ref{a1a2}:

\begin{proposition}\sl
If $a_{1} \in \widetilde{\CS}_{h} ( m_{1} )$ and $a_{2}\in \widetilde{\CS}_{h} (m_{2})$, then ${\rm Op}_{h} (a_{1}) \circ {\rm Op}_{h} (a_{2})$ is a pseudodifferential operator of class $\widetilde{\CS}_{h} ( m_{1} m_{2} )$ and its symbol is given by
\begin{align}
a\#b(x,\xi) =& e^{\frac{ih}{2} \sigma(D_{x},D_{\xi}, D_{y},D_{\eta})} \big( a(x,\xi)b(y,\eta)  \big) {\Big |_{y=x,\eta=\xi}}   \label{rez}  \\
=& \sum_{k=0}^{N-1} \frac 1{k!} \Big( \big(\frac{ih}2 \sigma(D_{x},D_{\xi},D_{y},D_{\eta})\big)^k a(x,\xi)b(y,\eta) \Big) \Big|_{y=x,\eta=\xi} \\
&+ h^{N} \widetilde{\CS}_{h} \big( \langle x \rangle^{-N/2} \langle \xi \rangle^{-N/2} m_{1} m_{2} \big) .
\end{align}
\end{proposition}

\begin{proof}
We follow the proof of \cite[Proposition 7.7]{disj}. Since $a_{j} \in \CS_{h}^{0} (m_{j})$, we now that ${\rm Op}_{h} (a_{1}) \circ {\rm Op}_{h} (a_{2})$ is a pseudodifferential operator whose symbol in $\CS_{h}^{0} (m_{1} m_{2})$ is given by (\ref{rez}). Let $X = (x,y, \xi , \eta)$, $\widetilde{X} = ( \widetilde{x} , \widetilde{y} , \widetilde{\xi} , \widetilde{\eta} )$ and $\chi \in C_{0}^{\infty} (\R^{4d})$ be equal to $1$ near $0$. Using  Fourier's inversion formula, one can show that, if we set
\begin{equation}
I (X) = e^{\frac{ih}{2} \sigma(D_{\widetilde{X}})} \Big( \chi \Big( \frac{X- \widetilde{X} }{\langle X \rangle^{\mu}} \Big) a_{1} a_{2} (\widetilde{X}) \Big) (X) - \sum_{j \leq N-1}\frac{1}{ j !} \left(\frac{ ih}{2} \sigma (D_{X})  \right)^{j} a_{1} a_{2} (X),
\end{equation}
we have
\begin{equation}
\vert I (X) \vert \lesssim h^{N} \sum_{\vert \alpha \vert \leq d/2+1} \Big\Vert D_{\widetilde{X}}^{\alpha} \sigma ( D_{\widetilde{X}} )^{N} \Big( \chi \Big( \frac{X- \widetilde{X} }{\langle X \rangle^{\mu}} \Big) a_{1} a_{2} (\widetilde{X}) \Big) \Big\Vert_{L^{2}_{\widetilde{X}}}.  \label{opo}
\end{equation}

Now, using the estimate (\ref{shdelta2}), we have
\begin{align}
\nonumber
\Big\vert D_{\widetilde{X}}^{\alpha} \sigma ( D_{\widetilde{X}} )^{N} & \Big( \chi \Big( \frac{X- \widetilde{X} }{\langle X \rangle^{\mu}} \Big) a_{1} a_{2} (\widetilde{X}) \Big) \Big\vert    \\
\nonumber
\lesssim & \sum_{\vert \beta \vert + \vert \gamma \vert = N} \Big\vert D_{\widetilde{X}}^{\alpha} \partial_{\widetilde{x}}^{\beta} \partial_{\widetilde{\eta}}^{\beta} \partial_{\widetilde{y}}^{\gamma} \partial_{\widetilde{\xi}}^{\gamma} \Big( \chi \Big( \frac{X- \widetilde{X} }{\langle X \rangle^{\mu}} \Big) a_{1} a_{2} (\widetilde{X}) \Big) \Big\vert   \\
\nonumber
\lesssim &   \sum_{\vert \beta \vert + \vert \gamma \vert = N} \sum_{\fract{0 \leq j,k \leq \vert \beta \vert}{0 \leq l,m \leq \vert \gamma \vert}} \langle X \rangle^{-j \mu} \langle \widetilde{x} \rangle^{- (\vert \beta \vert -j) /2} \langle X \rangle^{-k \mu} \langle \widetilde{\eta} \rangle^{- (\vert \beta \vert -k) /2}   \\
&\langle X \rangle^{-l \mu} \langle \widetilde{y} \rangle^{- (\vert \gamma \vert -l) /2} \langle X \rangle^{-m \mu} \langle \widetilde{\xi} \rangle^{- (\vert \gamma \vert -m) /2} m_{1} (\widetilde{X}) m_{2} (\widetilde{X}).
\end{align}
But since $\widetilde{X}$ is in the support of $\chi \big( (X-\widetilde{X}) \langle X \rangle^{- \mu} \big)$, we also have
\begin{equation}
\langle X \rangle^{-j \mu} \langle \widetilde{x} \rangle^{- (\vert \beta \vert -j) /2} \lesssim \left\{ \begin{aligned}
&\langle X \rangle^{- 2 \vert \beta \vert \mu /3} \quad &&\text{ for } j \geq 2 \vert \beta \vert /3 \\
& \langle x \rangle^{- \vert \beta \vert /6} \langle X \rangle^{\vert \beta \vert \mu /2} \quad &&\text{ for } j \leq 2 \vert \beta \vert /3
\end{aligned}
\right.
\end{equation}
Therefore, using the fact that $m_{j}$ are order functions, we obtain the estimate
\begin{align}
\nonumber
\vert I (X) \vert \lesssim & h^{N} \langle X \rangle^{N_{0}} m_{1} m_{2} (X) \sum_{\vert \beta \vert + \vert \gamma \vert = N} \big( \langle X \rangle^{- 2 \vert \beta \vert \mu /3} + \langle x \rangle^{- \vert \beta \vert /6} \langle X \rangle^{\vert \beta \vert \mu /2} \big)  \\
\nonumber
&\big( \langle X \rangle^{- 2 \vert \beta \vert \mu /3} + \langle \eta \rangle^{- \vert \beta \vert /6} \langle X \rangle^{\vert \beta \vert \mu /2} \big) \big( \langle X \rangle^{- 2 \vert \gamma \vert \mu /3} + \langle y \rangle^{- \vert \gamma \vert /6} \langle X \rangle^{\vert \beta \vert \mu /2} \big) \\
&\big( \langle X \rangle^{- 2 \vert \gamma \vert \mu /3} + \langle \xi \rangle^{- \vert \gamma \vert /3} \langle X \rangle^{\vert \beta \vert \mu /6} \big), 
 \label{popopo}
\end{align}
where $N_{0}$ is independent of $N$. 

Now, if we assume $y = x$ and $\eta = \xi$, we have,
\begin{align}
\nonumber
\big( \langle X \rangle^{- 2 \vert \beta \vert \mu /3} + \langle x \rangle^{- \vert \beta \vert /6} & \langle X \rangle^{\vert \beta \vert \mu /2} \big) \big( \langle X \rangle^{- 2 \vert \beta \vert \mu /3} + \langle \eta \rangle^{- \vert \beta \vert /6} \langle X \rangle^{\vert \beta \vert \mu /2} \big)   \\
\nonumber
&\lesssim \langle x, \xi \rangle^{-  \vert \beta \vert \mu /6} + \langle x \rangle^{- \vert \beta \vert /6} \langle \xi  \rangle^{- \vert \beta \vert /6} \langle x, \xi \rangle^{\vert \beta \vert \mu}  \\
&\lesssim \langle x, \xi \rangle^{- \vert \beta \vert (1/6 - \mu )},
\end{align}
so that (\ref{popopo}) gives
\begin{equation}  \label{rtt1}
\vert I( x, \xi , x, \xi ) \vert \lesssim  h^{N} \langle x, \xi \rangle^{- (1/6 - \mu) N + N_{0}}m_{1} (x, \xi ) m_{2} (x, \xi ).
\end{equation}

One obtains the same way  the same estimate for the derivatives of $I$, and we are left with the estimate of
\begin{align}
J (X) :=& e^{\frac{ih}{2} \sigma(D_{\widetilde{X}})} \Big( \Big( 1 - \chi \Big( \frac{X- \widetilde{X} }{\langle X \rangle^{\mu}} \Big) \Big) a_{1} a_{2} (\widetilde{X}) \Big) (X)  \\
=& \frac{1}{(2 \pi h)^{2d}} \int e^{-2i \sigma (X - \widetilde{X} )/h} \big( 1 - \chi \big( (X - \widetilde{X} ) \langle X \rangle^{- \mu} \big) \big) a_{1} a_{2} (\widetilde{X}) d \widetilde{X}.
\end{align}
We make integrations by parts, using the operator
\begin{equation}
{}^{t} L = \big( \partial_{\widetilde{X}} \sigma (X - \widetilde{X}) \big)^{-2} \big( \partial_{\widetilde{X}} \sigma (X - \widetilde{X}) \big) . \partial_{\widetilde{X}}.
\end{equation}
At each integration, we gain a factor $h$ and an  $|X - \widetilde{X}|^{-1}$, which is lower than $\langle X \rangle^{- \mu}$ on the support of $1- \chi$. Then, for each $M \gg 1$,
\begin{align}
\nonumber
\vert J (X) \vert \lesssim& h^{M - 2d} \langle X \rangle^{- (M - M_{0}) \mu} \sum_{\vert \alpha \vert \leq M}  \big\Vert \langle X - \widetilde{X} \rangle^{M_{0}} \partial_{\widetilde{X}}^{\alpha} a_{1} a_{2} \big\Vert_{L^{1}_{\widetilde{X}}}  \\
\lesssim & h^{M - 2d} \langle X \rangle^{- (M - M_{0}) \mu} m_{1} (X) m_{2} (X) \big\Vert \langle X - \widetilde{X} \rangle^{M_{0}} \langle X - \widetilde{X} \rangle^{M_{0}} \big\Vert_{L^{1}_{\widetilde{X}}}.
\end{align}
Eventually, if $y=x$ and $\eta = \xi$, we get, for all $M \in \N$,
\begin{equation}  \label{rtt2}
\vert J( x, \xi , x, \xi ) \vert \lesssim  h^{M} \langle x, \xi \rangle^{-M} m_{1} (x, \xi ) m_{2} (x, \xi ).
\end{equation}

We can prove also the same estimates for the derivatives of $J$, and  the proposition follows from (\ref{rtt1}) and (\ref{rtt2}).
\end{proof}

% \bibliographystyle{amsplain}
% \bibliography{reflq}
% 
% \end{document}

\end{document}